\documentclass[11pt,reqno]{amsart} 
\usepackage{amsthm,amssymb,amsfonts,amsmath}
\usepackage{amsrefs}
\usepackage{fullpage}
\usepackage{bm, bbm, mathrsfs}
\usepackage{graphics, graphicx}
\usepackage[mathscr]{euscript}
\usepackage{stmaryrd}

\theoremstyle{plain}

\theoremstyle{definition}

\theoremstyle{remark}
\newtheorem{remark}{Remark}
\newtheorem{example}{Example}

\newcommand{\mi}{\ensuremath{\mathrm{i}}}
\newcommand{\dif}{\ensuremath{\mathrm{d}}}

\newcommand{\CC}{\ensuremath{\mathbb{C}}}
\newcommand{\RR}{\ensuremath{\mathbb{R}}}

\newcommand{\lam}{\ensuremath{\lambda}}

\newcommand{\Ea}{\ensuremath{\mathcal{E}_A}}
\newcommand{\Ti}{\ensuremath{T_\mathrm{ig}}}
\newcommand{\ei}{\ensuremath{e_\mathrm{ig}}}

\newcommand{\f}{\ensuremath{\mathsf{f}}}

\newcommand{\calW}{\ensuremath{\mathcal{W}}}
\renewcommand{\t}{\ensuremath{\tau}}
\renewcommand{\l}{\ensuremath{\lambda}}
\newcommand{\tr}{\ensuremath{\top}}

\newcommand{\spm}{{\ensuremath{{\scriptscriptstyle\pm}}}}
\renewcommand{\sp}{{\ensuremath{{\scriptscriptstyle +}}}}
\newcommand{\sm}{{\ensuremath{{\scriptscriptstyle -}}}}

\newcommand{\bp}{\begin{pmatrix}}
\newcommand{\ep}{\end{pmatrix}}
\newcommand{\beq}{\begin{equation}}
\newcommand{\eeq}{\end{equation}}
\newcommand{\br}{\begin{remark}}
\newcommand{\er}{\end{remark}}
\newcommand\bt{\begin{todo}}
\newcommand\et{\end{todo}}

\numberwithin{equation}{section}

\title[Viscous Hyperstability of Detonation Waves]{Viscous hyperstabilization of detonation waves in one space dimension}
\author{Blake Barker}
\address{Department of Mathematics, Indiana University, Bloomington, IN 47405}
\email{bhbarker@indiana.edu}
\thanks{B.B. was partially supported by the National Science Foundation under grant DMS-0801745.} 
\author{Jeffrey Humpherys}
\address{Department of Mathematics, Brigham Young University, Provo, UT 84602}
\email{jeffh@math.byu.edu}
\thanks{J.H. was partially supported by the National Science Foundation under grant DMS-0847074 (CAREER)}
\author{Gregory Lyng}
\address{Department of Mathematics, University of Wyoming, Laramie, WY 82071}
\email{glyng@uwyo.edu}
\thanks{G.L. was partially supported by the National Science Foundation under grant DMS-0845127 (CAREER)}
\author{Kevin Zumbrun}
\address{Department of Mathematics, Indiana University, Bloomington, IN 47405}
\email{kzumbrun@indiana.edu}
\thanks{K.Z. was partially supported by the National Science Foundation under grant DMS-0801745.}
\date{Last Updated:  \today}
\begin{document}
\begin{abstract} 
It has long been a standard practice to neglect diffusive effects in stability analyses of detonation waves. 
	Here, with the principal aim of quantifying the impact of these oft-neglected effects on the stability characteristics of such waves, we use numerical Evans-function techniques to study the (spectral) stability of viscous strong detonation waves---particular traveling-wave solutions of the Navier--Stokes equations modeling a mixture of reacting gases. Remarkably, our results show a surprising synergy between the high-activation-energy limit typically studied in stability analyses of detonation waves and the presence of small but nonzero diffusive effects. While our calculations do show a modest delay in the onset of instability in agreement with recently reported calculations by direct numerical simulation of the physical equations, our approach, based on the Evans function, also provides additional spectral information. In particular, for each of the families of detonation waves in our computational domain, we find a completely unexpected kind of hysteresis in the limit of increasing activation energy; that is, \emph{our calculations suggest that, whenever diffusive effects are present, there is a return to stability as unstable eigenvalues return to the stable complex half plane}.   
%
%
%
\end{abstract}
\maketitle

\section{Introduction}\label{sec:intro}

\subsection{Viscous detonation waves \& stability}\label{ssec:vdw}
A detonation wave is a particular and dramatic kind of combustion wave that arises in the dynamics of mixtures of reacting gases. These distinctive waves are characterized by their high speed, substantial energy conversion, and shock-like structure \cite{FD}. Indeed, the traditional viewpoint---going back at least to Zel'dovich, von Neumann, and D\"oring in the 1940s---is that these waves are initiated by a gas-dynamical shock wave (which compresses the gas mixture, heating it and inducing chemical reactions) and their structure is that of a shock wave coupled to and supported by a burning tail \cites{CF,FD}. Importantly, the exothermic energy release due to chemical reaction plays an important role in the dynamics of detonation waves. That is, in order to correctly describe the dynamics, it is necessary to properly model the interaction of the combustion processes and the nonlinear behavior of the gas mixture. A standard tenet of the Zel'dovich--von Neumann--D\"oring (ZND) theory is that diffusive effects\footnote{These are the effects of viscosity, heat conductivity, and species diffusion. By an abuse of language, we frequently refer to the collection of all of these effects as the \emph{viscous effects} or simply the \emph{viscosity}.} can be safely neglected as negligible relative to the more significant reaction and convection effects which are expected to dominate the behavior \cites{CF,FD,W}.    
Indeed, it is a standard practice to model detonations using the (inviscid)  reactive Euler/ZND equations \cite{BS_ARFM07}.
In recent years, however, there has been a growing interest in the impact of these oft-neglected diffusive effects on the structure and stability of detonation waves. For example, Gardner \cite{G_TAMS83}, Gasser \& Szmolyan \cite{GS_SIMA93}, and Williams \cite{W_IUMJ10} have all established---under various structural assumptions---the existence of steady 
detonation waves with viscosity, and notable mathematical progress towards a stability theory for these waves can be found in the work of   
Tan \& Tesei \cite{TT_N97}; 
Lyng \& Zumbrun \cite{LZ_ARMA04}; 
Jenssen, Lyng, \& Williams \cite{JLW_IUMJ05}; 
Lyng et al.\ \cite{LRTZ_JDE07}; 
Texier \& Zumbrun \cite{TZ_CMP11}; 
and Zumbrun \cite{Z_ARMA11}. 
Complementing this theoretical work, we mention also the important recent numerical calculations of Romick, Aslan, \& Powers \cites{RAP_AIAA11,RAP_JFM12} and also the numerical calculations of Colella, Majda, \& Roytburd \cite{CMR_SJSSC86} and by Bourlioux, Majda, \& Roytburd \cite{BMR}; 
as discussed in \S\ref{ssec:vs} these calculations, especially the recent ones, form an important foundation for the present work.  

The first substantive investigation of the linear stability of detonation waves was initiated by Erpenbeck \cites{E_PF62,E_PF64} at the Los Alamos National Laboratory during the 1960s; Erpenbeck's program for analyzing the stability of ZND detonation waves is described in detail in the text of Fickett \& Davis \cite{FD}, and it forms the foundation on which subsequent work in this area has been built.  
For example, Lee \& Stewart \cite{LS_JFM90}, still in the ZND setting, considerably improved the state of affairs by providing a far more complete description of the neutral stability boundaries in parameter space\footnote{This problem was recently revisited, inside the Evans-function framework, by Humpherys \& Zumbrun \cite{HZ_QAM12}.}. Here, continuing the Evans-function program initiated by Lyng \& Zumbrun \cite{LZ_ARMA04}, we examine, by numerically computing the relevant Evans function, the stability of strong-detonation-wave solutions of the Navier--Stokes equations modeling a mixture of reacting gases. The Evans function, an analytic function of a complex variable,
is a ``stability function'' analogous to the stability function introduced by Erpenbeck \cites{E_PF62,E_PF64} for the ZND system. In particular, its zeros in the unstable complex half plane signal the presence of perturbations that grow exponentially in time. That is, the Evans function is constructed so that its zeros coincide in location and multiplicity with the eigenvalues of the linearized operator about the wave.
There is now a well-developed and growing literature describing the numerical approximation of the Evans function in the related case of viscous shock profiles \cites{BHLZ-el,BHLZ-flux,B_MC01,BZ_MC02,HSZ_NM06,HZ_PD06}. These techniques have been successfully used to analyze the stability of viscous shock profiles in gas dynamics \cites{BHRZ_CMP08,HLZ_ARMA09} and for detonation waves in simplified, scalar models \cites{BZ_Majda-ZND,HLZ_Majda,HHLZ_Majda}. Here, we adapt these techniques to the physically relevant and computationally challenging case of viscous detonation waves in the setting of the Navier--Stokes equations, system \eqref{eq:rns} below.

\subsection{Description of results: viscous hyperstabilization}\label{ssec:results}
Here, briefly and in an informal fashion, we describe the principal results of our computations. To start, we outline the steps required to approximate the Evans function for a given detonation wave. The first step is to approximate the wave itself. As a traveling-wave solution of the Navier--Stokes equations (system \eqref{eq:rns} below) the viscous structure of the wave is encoded in a nonlinear two-point boundary-value problem posed on the real line. Said differently, the wave is a heteroclinic orbit in the phase space $\RR^4$. To approximate it, we truncate the spatial domain, we apply appropriate projective boundary conditions at the end points of our (now finite) computational domain, and we employ a collocation method to solve the differential equation; a detailed description of this step can be found in \S\ref{ssec:profiles} below. Once the traveling-wave profile is in hand, we may turn our attention to the computation of the Evans function. 
This computation requires the solution of a nonautonomous linear system of differential equations whose coefficients depend on the structure of the wave in question and on the spectral parameter. Here, we use \textsc{StabLab}, a computational package developed precisely for this purpose \cite{STABLAB}, to carry out these calculations. One obstacle to the use of the Evans function for large physical systems like \eqref{eq:rns} below has been the inherent stiffness in the problem that arises from the need to resolve modes of differing exponential growth (decay) rates to obtain the subspaces of growing (decaying) solutions from which the Evans function is built. Importantly, the computational routines implemented by \textsc{StabLab} allow the user to deal with this issue either by the traditional method of ``exterior products'' or by the ``analytic orthogonalization''  method used here. In either case, care is taken to maintain the analyticity of the Evans function with respect to the spectral parameter; this facilitates the use of winding number calculations to search for zeros (unstable eigenvalues).
A more complete description of the methodology for computing the Evans function can be found below in \S\ref{sec:methodology} and in \cites{BHLZ-el,BHLZ-flux,B_MC01,BZ_MC02,HSZ_NM06,HZ_PD06}.

The key point is that, with the Evans function in hand, we may systematically search for unstable eigenvalues of the linearized operator. The first task, then, is to determine the shape of the neutral stability boundary and to compare it to the neutral stability boundary, as computed by Erpenbeck and his successors, for the corresponding ZND profiles. For example, Romick et al.\ \cites{RAP_AIAA11,RAP_JFM12}, using direct numerical simulation of the temporal evolution of the physical system of equations, found that the presence of viscosity delayed the onset of instability by about 10\%. We recover this result, but our approach yields more information. 
In particular, for each of the families of detonation waves in our computational domain, we find a completely unexpected kind of hysteresis in the limit of increasing activation energy; that is, our calculations suggest that, whenever diffusive effects are present, there is a return to stability as unstable eigenvalues return to the stable complex half plane. \emph{We call this behavior \textbf{viscous hyperstabilization}, and we conjecture that it persists for all nonzero viscosities}. We believe that the further investigation of this phenomenon via singular perturbation theory is a fundamental open problem in the theory.

\begin{remark}[Growth of the upper stability boundary]
A particularly interesting feature of these calculations is the growth of the upper-stability boundary in the limit of vanishing viscosity. Zumbrun has shown that this boundary must escape to infinity as the viscosity tends to zero \cite{Z_ARMA11}, but our computational results suggest a very slow, perhaps logarithmic, growth, in the stabilizing activation energy as a function of viscosity. See Figure \ref{fig:sb} below. We believe that further investigation of this phenomenon is warranted.  
\end{remark}

\subsection{About the computations}
We see several noteworthy aspects of these computations which we now describe in some detail.

\subsubsection{Computational tractability}\label{sssec:compute} 
First and foremost, we see the computation of the Evans function in this physically relevant setting as a substantial step forward. Indeed, we view this as an important, practical validation of the power of the polar-coordinate (``analytic orthogonalization'') technique \cite{HZ_PD06} for approximating the Evans function. For context, we recall that the computation of the Evans--Lopatinski\u\i\ function (essentially Erpenbeck's stability function)
for the ZND system is well known to be a computationally intense problem. For more details about the difficulties, we refer the reader to the foundational paper of Lee \& Stewart \cite{LS_JFM90} and also to more recent work by Humpherys \& Zumbrun \cite{HZ_QAM12}. Briefly, however, we can see by a dimension count that the ZND calculation requires the (stiff) computation of a single mode of a $4\times 4$ system while the computation for the equations with viscosity requires one to track 3 growing and 3 decaying solutions of a $7\times 7$ system. Thus, the (assumed small) effects of viscosity are responsible for a considerable increase in the size and complexity of the linearized system.
We also note that the use of Lagrangian coordinates (cf.\ Remark \ref{rem:el}) and the use of flux variables (cf.\ equation \eqref{eq:Ydef} in \S\ref{ssec:flux}) are important ingredients in our computational formulation; full details of these devices and their virtues can be found in \cites{BHLZ-el,BHLZ-flux}.

\subsubsection{Virtues of the Evans-function approach}\label{sssec:evansapproach} 
A particular benefit of our approach, based on the Evans function, is that it reveals information about the eigenstructure of the linearized operator. For example, as noted above, one of our experiments utilizes the same parameter values as used in the recent studies of Romick et al.\ \cites{RAP_AIAA11,RAP_JFM12}.
Their computations, however, are based on the direct numerical simulation of the temporal evolution of the governing physical system of partial differential equations.
%
Here, using the same parameters, we find the same 10\% delay in the onset of instability as seen by Romick et al.\ \cites{RAP_AIAA11,RAP_JFM12}. At first glance, this supports the conventional wisdom that diffusive effects are essentially perturbative. However, as noted above, continuing to follow roots as the activation energy \Ea\ is increased, we find a return to stability---viscous hyperstabilization; see Figure \ref{fig:PowersEA} below. Moreover, rather than the well-known cascade of Hopf bifurcations associated with the stability problem in the ZND setting, \emph{we see only two pairs of complex conjugate roots crossing the imaginary axis into the unstable complex half plane.} Thus, the inclusion of diffusive effects in the model results in striking changes at the spectral level. The mechanism for this behavior is a kind of  synergy between viscous effects and activation energy.  The inviscid limit is highly nonuniform with respect to $\Ea$, and this nonuniformity manifests itself in both the structure of the viscous profiles (see the discussion of ``benches'' in \S\ref{ssec:benches} below) and subsequently in eigenvalues\footnote{We recall the paradigm of Henry \cite{H}: essential spectrum of the linearized operator corresponds to far-field behavior while the point spectrum (located by the Evans function) is associated with the detailed structure of the shock layer and nearby reaction zone.}. We note that the benching phenomenon appears to occur for a rather  pedestrian reason. That is, its origin is due to the different powers of $\tau$, the specific volume, that appear in the viscous terms on the right-hand side in the system \eqref{eq:rns} below.  

\subsubsection{Implications \& future work}\label{sssec:implications} 
Finally, and perhaps most importantly in the long run, these findings may be significant in various applications in certain parameter regimes. For example, we note that instability characteristics of detonation waves are critically relevant in the design of any viable detonation-based propulsion system, a problem which remains largely open \cite{SK_JPP06}.
While the implications for applications remain to be explored, we note that, from the point of view of applied mathematics and modeling, these important behaviors are ``right on the surface,'' and our findings certainly warrant a more detailed analysis of the viscous system. The Evans-based techniques utilized here appear to offer an effective starting point for this endeavor. 


\subsection{Plan}
In \S\ref{sec:model} we recall the form of the Navier--Stokes equations describing a reacting mixture of gases. We fix notation, and we describe the forms of the constitutive relations (equations of state) that we use in our computations. In \S\ref{sec:profiles} we discuss the traveling-wave equation and a particularly useful scaling that we use. In this section we also discuss our approach to approximating the solutions of the traveling-wave equation. Indeed, even at this initial stage, we find some interesting and unexpected behavior in the small-viscosity regime. The linearization and the subsequent computation of the Evans function are discussed in \S\ref{sec:linearize} and \S\ref{sec:methodology} although a number of the details of the tedious computation of the linearization are deferred to Appendix \ref{sec:details}. In \S\ref{sec:results}, we present the principal result of the paper. Finally, in \S\ref{sec:discuss} we discuss our findings and suggest some possible future directions for research.

\section{Equations for reacting mixtures of gases}\label{sec:model}
\subsection{Navier--Stokes equations}\label{ssec:rns}
In a single spatial dimension (Lagrangian coordinates), the Navier--Stokes equations for a reacting two-species gas mixture take the form 
\begin{subequations}\label{eq:rns}
\begin{align}
\tau_t&-u_x=0\,,\label{eq:mass}\\
u_t&+p_x=\left(\frac{\nu u_x}{\tau}\right)_x\,,\label{eq:momentum}\\
E_t&+(pu)_x=\left(\frac{\nu uu_x}{\tau}+\frac{\kappa T_x}{\tau}\right)_x+qk\phi(T)z\,,\label{eq:energy}\\
z_t&=-k\phi(T)z+\left(\frac{dz_x}{\tau^2}\right)_x\,.\label{eq:reaction}
\end{align}
\end{subequations}
Here, subscripts denote partial derivatives ($t$ is time and $x$ is a Lagrangian marker), $\tau$ is the specific volume (length), $u$ is the velocity, $p$ is the pressure, and $z\in[0,1]$ is the mass fraction of reactant. The specific total energy $E$ is made up of the specific internal energy $e$ and the kinetic energy:
\begin{equation}
E=e+\frac{u^2}{2}\,.\label{eq:internal_kinetic}
\end{equation}
The positive constants $\nu$, $\kappa$, and $d$ represent viscosity, heat conductivity, and species diffusivity respectively. The constant $k>0$ is the reaction rate, and $q$ measures the difference in the heats of formation in the reactant and the product. That is, it is a measure of the amount of energy released during the reaction process. We assume $q>0$ which corresponds to an exothermic reaction. The function $\phi$ is the ignition function; it serves as an on/off switch for the reaction. Finally, $T$ is the temperature, and we assume initially that the internal energy $e$ and the pressure $p$ are known functions of the specific volume, the temperature, and the mixture of the product and reactant gases:
\begin{equation*}
p=p_0(\tau,T,z),\quad e=e_0(\tau,T,z).
\end{equation*}
We describe the functions $p_0$ and $e_0$ in \S\ref{ssec:eosi} immediately below.  
A detailed derivation of the system \eqref{eq:rns} can be found in the text by Williams~\cite{W}. For brevity, we shall refer to the system \eqref{eq:rns} of partial differential equations as the RNS system or as the the RNS equations. 

\begin{remark}[Eulerian vs. Lagrangian]\label{rem:el}
Here, as in previous studies of the stability of viscous shocks in gas dynamics \cites{BHRZ_CMP08,HLZ_ARMA09}, we work in Lagrangian rather than Eulerian coordinates. Now, the same viscous profile (see \S\ref{sec:profiles} below) and the associated linearization (see \S\ref{sec:linearize} and Appendix \ref{sec:details} below) can be expressed in either of the two different coordinate systems, and, evidently, the corresponding Evans functions must carry the same stability information (zero set). However, this does not mean that the two different analytic functions provide equally convenient platforms from which to attempt a (numerical) stability analysis. In particular, a cornerstone of our approach is to use the analyticity of the Evans function to count zeros by computing winding numbers. In this framework, for example, it is clearly of interest to eliminate excessive winding and unwinding about the origin. Remarkably, our various experiments have shown that, for viscous shock profiles in gas dynamics, the Evans function arising from the Lagrangian system is considerably better behaved than the one that arises from reaction-convection-diffusion system \eqref{eq:rns} expressed in Eulerian coordinates. See \cite{BHLZ-el} for a further, detailed discussion of this point in the setting on gas dynamics and for computations illustrating the transformation from one Evans function to the other. 
\end{remark}

\subsection{Equations of state \& ignition}\label{ssec:eosi}
\subsubsection{Ideal gas}\label{ssec:ideal}
The simplest case of interest occurs when we consider an ideal,  
polytropic gas. In this case the energy and pressure functions take  
the specific form
\begin{equation}\label{eq:ideal_gas}
p_0(\tau,T,z)=\frac{RT}{\tau},\quad e_0(\tau,T,z)=c_vT,
\end{equation}
where $R\ge 0$ and $c_v>0$ are constants that characterize the gas.
Alternatively, we can write the pressure in terms of the internal energy and the specific volume as
\begin{equation}\label{Gammaeq}
p=\frac{\Gamma e}{\tau},
\end{equation}
where $\Gamma =\gamma -1= \frac{R}{c_v }\ge 0$,
$\gamma \ge 1$ the gas constant. We note that due to \eqref{eq:ideal_gas}, the specific internal energy is simply a multiple of the temperature, and these two quantities can effectively be used interchangeably.  

\subsubsection{Ignition function}\label{sssec:ignition}
We make the standard assumption that the ignition function $\phi$ is of the Arrhenius form with an ignition temperature cut-off at a distinguished  threshold temperature \Ti. That is, we assume
\beq
\phi(T)=\begin{cases}
\exp\left(-\Ea/(T-\Ti)\right)\,, & \text{if}\;T\geq\Ti \\
0\,, & \text{otherwise}
\end{cases}
\eeq
with $\Ea$ the activation energy and $\Ti$ the ignition temperature. Using the relationship $e=c_vT$ in \eqref{eq:ideal_gas} between the temperature and internal energy, we sometimes write $\check\phi(e)=\phi(T)$. For our numerical experiments we choose the ignition value of the specific internal energy, $\ei$, to be a convex combination of the end state $e_\sp$ and $e_\mathrm{mid}$---defined to be the value of the specific internal energy at $x=0$ (at the jump) in the corresponding ZND solution. The convex combination is weighted heavily to favor the contribution of the unburned state $e_\sp$. 

\section{Traveling-wave profiles \& scaling}\label{sec:profiles}
\subsection{Traveling waves}\label{ssec:tw}
The viscous strong detonation waves that are the focus of our analysis are traveling-wave solutions of \eqref{eq:rns}. 
Suppose that $(\tau,u,e,z)(t,x)=(\hat\tau,\hat u,\hat e,\hat z)(x-st)$ is a traveling-wave solution connecting constant states $(\tau_\spm,u_\spm,e_\spm,z_\spm)$. Our interest is in waves connecting an unburnt state ($z_\sp=1$) to a completely burnt state ($z_\sm=0$).

\subsection{A convenient scaling}\label{ssec:scaling}
Evidently, a traveling-wave solution is a stationary solution of the system 
\begin{subequations}\label{eq:rns2}
\begin{align}
\tau_t&-s\tau_x-u_x=0\,,\label{eq:mass2}\\
u_t&-su_x+p_x=\left(\frac{\nu u_x}{\tau}\right)_x\,,\label{eq:momentum2}\\
E_t&-sE_x+(pu)_x=\left(\frac{\nu uu_x}{\tau}+\frac{\kappa_\nu e_x}{\tau}\right)_x+qk\check\phi(e)z\,,\label{eq:energy2}\\
z_t&-sz_x=-k\check\phi(e)z+\left(\frac{dz_x}{\tau^2}\right)_x\,.\label{eq:reaction2}
\end{align}
\end{subequations}
We scale the independent variables using $L>0$ as 
\beq\label{eq:scale1}
x\to\frac{\tau_\sp sx}{L}=:\bar x\,,\quad t\to\frac{\tau_\sp s^2t}{L}=:\bar t\,,
\eeq
and we scale the dependent variables as 
\beq
\tau\to\frac{\tau L}{\tau_\sp}=:\bar\tau\,, \quad
u\to\frac{uL}{\tau_\sp s}=:\bar u\,,\quad
e\to\frac{L^2e}{\tau_\sp^2s^2}=:\bar e\,.
\eeq
Then, with $\bar u(\bar t(t),\bar x(x))=\frac{Lu(t,x)}{\tau_\sp s}$ and so forth, the equations in barred variables read as 
\begin{subequations}\label{eq:rns3}
\begin{align}
\bar\tau_{\bar t}&-\bar\tau_{\bar x}-\bar u_{\bar x}=0\,,\label{eq:mass3}\\
\bar u_{\bar t}&-\bar u_{\bar x}+{\bar p}_{\bar x}=\left(\frac{\nu \bar u_{\bar x}}{\bar\tau}\right)_{\bar x}\,,\label{eq:momentum3}\\
\bar E_{\bar t}&-\bar E_{\bar x}+(\bar p\bar u)_{\bar x}=\left(\frac{\nu \bar u\bar u_{\bar x}}{\bar\tau}+\frac{\kappa_\nu \bar e_{\bar x}}{\bar\tau}\right)_x+\bar q\bar k\bar{\check\phi}(\bar e)\bar z\,,\label{eq:energy3}\\
\bar z_{\bar t}&-\bar z_{\bar x}=-\bar k\bar{\check\phi}(\bar e)\bar z+\left(\frac{d\bar z_{\bar x}}{\bar\tau^2}\right)_{\bar x}\,.\label{eq:reaction3}
\end{align}
\end{subequations}
In \eqref{eq:rns3}, we have used
\[
\bar{\check\phi}(\bar e)=\check\phi\left(\frac{\tau_\sp s}{L^2}\bar e\right)\,,\quad
\bar k =\frac{kL}{\tau_\sp s^2}\,,\quad
\bar q =\frac{qL^2}{\tau_\sp s^2}\,.
\]
Remarkably, the pressure law is \emph{unchanged} in this new formulation; for further discussion see \cite{HLZ_ARMA09}.
Henceforth, we drop the bars. 
From \eqref{eq:rns3}, the profile system is thus the stationary system of ordinary differential equations ($'=\dif/\dif x$)
\begin{subequations}\label{eq:twode}
\begin{align}
-\tau'&-u'=0\,, \label{eq:mass4} \\
-u'&+p'=\left(\frac{\nu u'}{\tau}\right)'\,,\label{eq:momentum4}\\
-E'&+(pu')=\left(\frac{\nu uu'}{\tau}+\frac{\kappa_\nu e'}{\tau}\right)'+ qk\check\phi(e) z\,,\label{eq:energy4}\\
-z'&=- k\check\phi(e)z+\left(\frac{d z'}{\tau^2}\right)'\,.\label{eq:reaction4}
\end{align}
\end{subequations}
Then using (from the scaling and translation invariance)
\[
s=1\,,\tau_\sp=1,\,u_\sp=0\,,
\]
and the ideal gas pressure law, we integrate the system \eqref{eq:twode} once, and we rewrite the result as a first-order system with 
\beq\label{eq:ydef}
y:=\frac{dz'}{\tau^2}\,.
\eeq
The result is
\begin{subequations}\label{eq:tw}
\begin{align}
\tau'&=-\frac{1}{\nu}\big[\tau(\tau-1)+\Gamma(e-e_\sp\tau)\big]\,, \label{eq:tautw}\\
e'&=-\frac{\tau}{\kappa_v}\left[\frac{(\tau-1)^2}{2}+(e-e_\sp)+\Gamma e_\sp(\tau-1)+q(y+z-z_\sp)\right]\,, \label{eq:etw}\\
z'&=d^{-1}\tau^2y\,, \label{eq:ztw}\\
y'&=k\check\phi(e)z-d^{-1}\tau^2y\,. \label{eq:ytw}
\end{align}
\end{subequations}
In system \eqref{eq:tw}, we have used the conservation of mass equation to eliminate the velocity $u$ from the system;
of course, the velocity $u$ and the specific volume $\tau$ are related by the simple relation 
\beq
u=1-\tau\,,
\eeq
and we have used the following convenient notation
\beq
\kappa_v:=\frac{\kappa}{c_v}\,.
\eeq

\begin{remark}[Existence of solutions]
As noted above, there are several proofs of the existence of such profile solutions. For example, Gasser \& Szmolyan \cites{GS_SIMA93} (and more recently by a distinct methodology, Williams \cite{W_IUMJ10}) have established the existence of solutions for sufficiently small viscosity. In both cases, the arguments rely on the use of the inviscid ZND solution as a template for the viscous profile, and they both build the layer solution using a multi-scale approach.
We also note that Gardner \cite{G_TAMS83} has established the existence of profiles using the topological techniques based on the Conley Index.  
\end{remark}

\subsection{Parametrization}\label{ssec:parametrization}
Summing up, we may take without loss of generality $s=1$, $\tau_\sp=1$
and (by translation invariance in $u$), $u_\sp=0$,
leaving $e_\sp>0$ as the parameter determining the wave.
By explicit computation using the 
Rankine--Hugoniot relations (\cite{Z_ARMA11}, Appendix C), we have then that the burned end state is given by
\begin{align}
\tau_\sm&= 
\frac{
(\Gamma+1) (\Gamma e_\sp+1)-
\sqrt{ 
(\Gamma+1)^2 (\Gamma e_\sp+1)^2
- \Gamma (\Gamma +2) ( 1+2(\Gamma +1)e_\sp +2q))  }
}
{\Gamma +2}\,,\\
u_\sm& = 1-\tau_\sm\,, \\
e_\sm& = 
\frac{\tau_-(\Gamma e_\sp+1-\tau_\sm)}{\Gamma},
\end{align}
where $q$ satisfies $0\leq q\leq q_\mathrm{CJ}$ and 
\begin{align}\label{cjlim}
q_{\mathrm{CJ}} & :=
\frac{ (\Gamma+1)^2(\Gamma e_+ + 1)^2- \Gamma (\Gamma+2) ( 1+2(\Gamma +1)e_\sp ) }
{2\Gamma(\Gamma+2)}\,.
\end{align}
Thus, we can parametrize all possible profiles by 
\beq\label{param}
(e_\sp, q, \Ea, \Gamma, \nu, d, \kappa_v, k)\,, 
\eeq
where 
$0\le e_\sp \le \frac{1}{\Gamma(\Gamma+1)}$,
$0\le q\le q_{\mathrm{CJ}}(e_\sp)$,
$0\le \Ea < \infty$, and $0<\Gamma<\infty$.
By a further common rescaling of space and time, 
we can fix $k$ at any desired value.
Following standard
practice, we shall choose the value of $k$ so that the 
reaction zone of the associated inviscid (ZND) profile is
of roughly constant length; more precisely, to match corresponding
ZND computations done previously in \cite{BZ_ZND}. In particular, we choose $k$ so that the corresponding ZND (inviscid) solution with this ignition function has the value $z=1/2$ at $x=-10$ in the spatial domain; the jump is at $x=0$. 
See \cite{Z_ARMA11} for further details.

\subsection{The high-overdrive limit and the scaling of Erpenbeck}\label{ssec:f}

A similar scaling was used by Erpenbeck in
\cites{E_PF62,E_PF64}, but with $e_\sp$ held fixed instead of wave speed $s$.
As pointed out previously \cites{HLZ_ARMA09,Z_ARMA11}, the advantage of 
our choice is that the coefficients of the linearized eigenvalue equations remain bounded in the limit of increasing strength of the  Neumann shock. 
Converting from Erpenbeck's to our scaling amounts to rescaling the
wave speed, so that $T\to T/s^2$ and $\Ea\to 
\Ea/s^2$, and $t\to ts^2$
($u$ is translation invariant, so irrelevant).
Thus, as noted by Zumbrun \cite{Z_ARMA11},
the high-overdrive limit discussed by Erpenbeck \cite{E_PF62}, in which
$s\to \infty$ with $u_\sp$ held fixed, corresponds in our scaling
to taking $\Ea= \mathcal{E}_0 e_\sp$, 
$q= q_0 e_\sp$, and varying $e_\sp$ from $e_\sp=e_\mathrm{cj}(q_0) $
 ($s$ minimum) to $ 0$ ($s=\infty$). Here, the value of $e_\mathrm{cj} $
is may be obtained from the relation $q_\mathrm{cj}(e_\mathrm{cj})=q_0 e_\mathrm{cj}$.
This is the simultaneous \emph{zero heat release},
{\it zero activation energy}, and {\it strong shock limit}
($e_\sp\to 0$, 
$\Ea\to 0$,
and $q\to0$).
To complete the translation from Erpenbeck's
coordinatization to ours, and back, we must compute the overdrive,
which Erpenbeck uses in place of our $e_\sp$ to complete the description
of the wave (recall that he holds $e_\sp$ fixed at value $1$).
%
The overdrive factor is defined for fixed $q_0$ and
$\tau_\sp=e_\sp=1$, as the ratio $\f:=M^2/M_{\mathrm{CJ}}^2$, where
$M$ is the Mach number and $M_{\mathrm{CJ}}$ is the Mach number
of the detonation with same right end state and $q_0$,
but traveling with CJ speed.
In the present setting (the scaling introduced above), we must therefore compute the Mach number
for a detonation with $q=q_0 e_{\sp}^\mathrm{CJ}$, $\tau_\sp=1$,
and $q_\mathrm{CJ}(e_{\sp}^\mathrm{CJ})=q_0 e_{\sp}^\mathrm{CJ}$, or
\beq\label{impcj}
q_0 e_{\sp}^\mathrm{CJ}=
\frac{ (\Gamma+1)^2(\Gamma e_{\sp}^\mathrm{CJ} + 1)^2- \Gamma (\Gamma+2) ( 1+2(\Gamma +1)e_\sp^\mathrm{CJ} ) }
{2\Gamma(\Gamma+2)}\,.
\eeq
Evidently, for $e_\sp^\mathrm{CJ}$ small, this gives asymptotically
$$
e_\sp^\mathrm{CJ}\sim 
\frac{1}{2\Gamma (\Gamma+2)q_0}\,,
$$
hence, noting that $q_0=q/e_\sp$, we find that
\beq\label{fform}
\f=\frac {e_\sp^\mathrm{CJ}} {e_\sp}
\sim \frac{1}{ 2\Gamma(\Gamma+2)q} ,
\eeq
where $q\lesssim \frac{1}{ 2\Gamma(\Gamma+2)}.  $
Taking the limit all the way to $e_\sp=0$, we may thus estimate
the overdrive simply in terms of $q$, via \eqref{fform}.

\begin{example}
For example, $e_\sp=0$, $q=0.1$, $\Gamma=0.2$ 
gives $\f\approx 11.3$, but $q_0$, $\mathcal{E}_0\to \infty$.
\end{example}

\begin{example}[Fickett \& Wood \cite{FW_PF66}]
To duplicate Fickett and Wood's benchmark problem
$q_0=\mathcal{E}_0=50$, $\Gamma=0.2$, we take
$e_\sp^\mathrm{CJ}\sim
 \frac{1}{2\Gamma (\Gamma+2)q_0} \approx
.023 ,  $
and vary $\mathcal{E}=\mathcal{E}_0e_\sp=50e_\sp$, $q=q_0 e_\sp=50e_\sp$
as $e_\sp$ varies between $e_{+}^\mathrm{CJ}$ and $0$, with
$1\leq \f= e_\sp^\mathrm{CJ}/e_\sp<+\infty$.
Here, we expect (and find) the onset of instability around $\f\approx 1.73$,
or $e_\sp\approx .012$, $q\approx (.012)(50)= .6$, $\Ea=.6$.
\end{example}
\begin{remark}[Fickett \& Wood \cite{FW_PF66}]
We observe that Fickett \& Wood's benchmark problem has now been duplicated by direct numerical simulation, evaluation of the Evans--Lopatinski\u\i\ function (ZND) by Lee \& Stewart \cite{LS_JFM90}, and now by numerical computation of the Evans function (RNS) here. Thus, there are now multiple independently computed sources confirming this behavior. 
\end{remark}

\subsection{Viscous Scaling}\label{ssec:vs}
To determine appropriate ranges for the viscous parameters in our model, we use---as in related studies \cites{RAP_AIAA11,RAP_JFM12}---the ratio of width of the viscous shock layer to the length of the reaction zone. That is, we choose the reaction rate $k$ (as discussed above in \S\ref{ssec:parametrization}) to hold the length of reaction zone constant, and then we vary $\nu$, $\kappa$, and $D$ to adjust the width of the shock layer. The key observation is that the ratio of shock to reaction width is independent of the rescalings we have performed above. As a starting point, we adopt the ratio of Romick et al. \cite{RAP_JFM12} in which the ratio of length scales is $1/10$.  As noted there, this ratio is an order of magnitude too big to be physical; indeed, Powers \& Paolucci have shown that a large range of length scales are important in realistic models (including detailed chemical kinetics) \cite{PP_AIAAJ05}. Thus, our experiments include a range of ratios including $1/10$, but we also compute closer to the inviscid limit; see Figure \ref{fig:sb}. For example, visual inspection of Figure \ref{fig:bench} shows a ratio of approximately $1/10$.

\subsection{Approximation of profiles and phenomenology}\label{ssec:profiles}
\subsubsection{Computational protocol}\label{ssec:protocol}
From the point of view of dynamical systems, a strong detonation is a heteroclinic orbit of \eqref{eq:tw} in $\RR^4$ connecting equilibria $(\t_\spm,e_\spm,z_\spm,y_\spm)$. Said differently, to find the detonation profile, we must solve a nonlinear two-point boundary-value problem posed on the whole real line. 
The numerical approximation of solutions of such problems is well understood. 
First, we truncate the problem to a finite computational domain $[-M_\sm,M_\sp]$ for $M_\spm>0$, and we supply appropriate projective boundary conditions at $\pm M_\spm$. Finally, we use \textsc{MatLab}'s boundary-value solver, an adaptive Lobatto quadrature scheme, to compute the approximate profile. A frequent challenge associated with this approach is the need to supply an adequate initial guess to the boundary-value solver to seed a Newton iteration. Because of this, the ability to continue solutions as parameter values change is an important art.
Additionally, to ensure that the the profile is completely resolved in the domain, the values for plus and minus spatial computational infinity, $M_\spm$, must be chosen with some care. Writing the traveling-wave equation \eqref{eq:tw} as $\hat U'=F(\hat U)$ together with the condition that $\hat U\to U_\spm$ as $\xi\to\pm\infty$, the typical requirement is that $M_\spm$ should be chosen so that $ |\hat U(\pm M_\spm)-U_\spm|$ is within a prescribed tolerance of 1e-4.  We also set the relative and absolute error tolerances for the boundary value solver to be 
1e-6 and 1e-8 respectively.

\subsubsection{``Benches''}\label{ssec:benches}
\begin{figure}[ht] 
   \centering
   \includegraphics[width=4in]{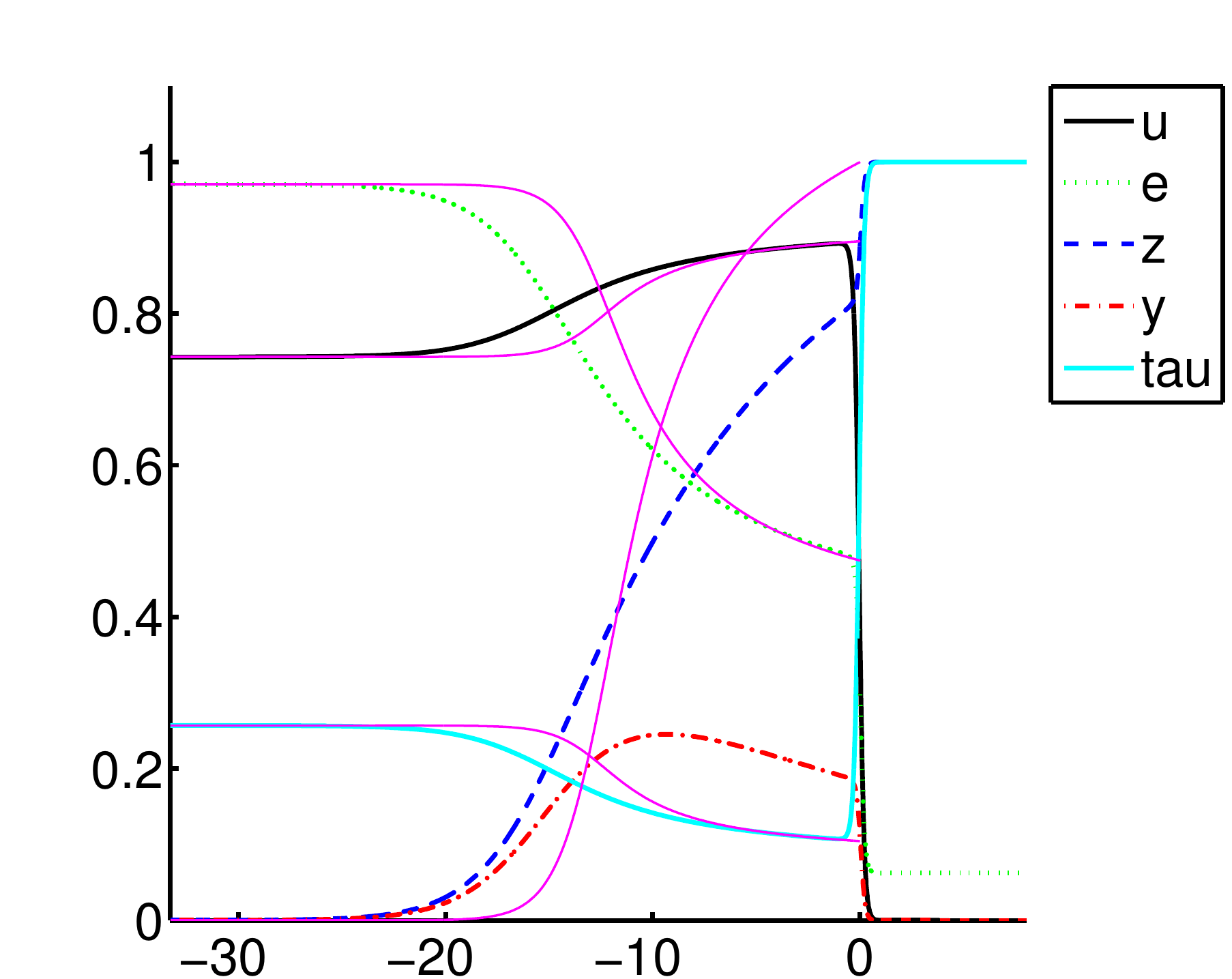} 
   \caption{A comparison of the viscous profile to the inviscid ZND solution. The legend shows the values of the velocity ($u$), specific internal energy ($e$), mass fraction of reactant ($z$), $y$, and specific volume ($\tau$). The thin (magenta) lines denote the components of the corresponding ZND solution. Note the discrepancy between the $z$ profiles for the two systems; the discrepancy is especially pronounced near the shock layer (close to $x=0$). Here, $e_\sp=6.23$e-2, $k=$2.71e-1, $d$=0.1, $\nu$=0.1, $\kappa$=0.1, $q$=6.23e-1, $\Ea$=3.1, $\Gamma$=0.2, $c_v$=1.}
   \label{fig:bench}
\end{figure}

An example illustrating the numerical solution, via the method of \S\ref{ssec:protocol}, of the traveling-wave equation \eqref{eq:tw} is shown in Figure \ref{fig:bench}. The figure compares the viscous solution to the ZND solution connecting the same end states. One interesting result of these computations regards the behavior of the $z$ component. In particular, visually, the deviation of the RNS $z$ profile from the limiting ZND profile is quite striking relative to that of the the other components of the solution. This discrepancy is particularly pronounced near the shock layer. As shown in Figure \ref{fig:bench}, the $z$ profile hits the shock layer at a height below $z=1$ and makes a near vertical adjustment in the shock layer to achieve the correct limiting value. 
We call this feature of the traveling wave a ``bench,'' and we saw benches of various sizes for all of the positive viscosity values tested. Further investigation reveals that the bench formation appears to be, more or less, a result of the distinct diffusive scalings in the traveling-wave equation. We note that the diffusivities scale like $1$ in \eqref{eq:tautw}, $\tau$ in \eqref{eq:etw}, and $\tau^2$ in \eqref{eq:ztw} and \eqref{eq:ytw} and that the value of $\tau$ is actually rather small for the profiles involved in our computations (roughly 0.15 in Figure \ref{fig:bench}). 
Thus, the diffusivities range over several orders of magnitude.



\section{Linearization \& Evans function}\label{sec:linearize}

Here, we outline the basics steps in the construction of the Evans function. The first step is to linearize the system \eqref{eq:rns2} about the steady strong-detonation profile. Taking the Laplace transform with respect to time, we obtain the eigenvalue problem for the linearized operator.  
For the convenience of the reader, we provide a detailed account of these computations in Appendix \ref{sec:details}. The upshot is the we obtain a linear, first-order system of ODEs with coefficients depending on the spectral parameter $\lam$. We write this system as 
\beq\label{eq:evanssystem}
\mathcal{W}'=G(x;\lam)\mathcal{W}\,,\quad '=\dif/\dif x\,,
\eeq
and we note that $\mathcal{W}\in\CC^7$ and the coefficient matrix $G$ is a $7\times 7$ matrix. The forty-nine entries of $G$ are given in \S\ref{ssec:system} in Appendix \ref{sec:details}. The Evans function $\lambda\mapsto D(\lambda)$ is computed from the system \eqref{eq:evanssystem}.


\section{Methodology}\label{sec:methodology}
There is now an established and growing collection of literature describing various techniques for the numerical approximation of the Evans function; see, e.g., \cites{B_MC01,HSZ_NM06,HZ_PD06,STABLAB}. We describe the basic steps below.
\begin{description}
\item[Step \#1. Approximation of the viscous profile] This initial step is described in detail in \S\ref{ssec:protocol} above.
\item[Step \#2. Approximation of the Evans function]
The evaluation of the Evans function is accomplished by means of the \textsc{StabLab} package, a \textsc{MatLab}-based package developed for this purpose \cite{STABLAB}. This package allows the user to choose to approximate the Evans function either via exterior products, as in \cites{AB_NM02,B_MC01,AS_NW95}, or by a polar-coordinate (``analytic orthogonalization'') method \cite{HZ_PD06}, which is used here. Because the analyticity of the Evans function is one of its essential properties, we use Kato's method \cite{K}*{p. 99} to analytically determine the relevant initializing eigenvectors; see \cites{BZ_MC02,BDG_PD02,HSZ_NM06} for details.  Throughout our study, we set the relative and absolute error tolerances on \textsc{MatLab}'s ODE solver {\tt ode45} to be 1e-6 and 1e-8 respectively.
\item[Step \#3. Winding \& Roots]
Finally, we describe how we find zeros of the Evans function $D$; our approach is based on the method of moments as in, e.g., \cites{OMB_PD86,B_PD96}. Being fundamentally based on Rouch\'e's theorem, this technique takes advantage of the analyticity of the Evans function. 
Suppose that $D$ is holomorphic on and inside a simple closed positively oriented contour $\gamma$, and denote by $\lam_1, \lam_2,\ldots, \lam_m$ the zeros of $D$ inside $\gamma$. Then, the $p$th moment of $D$ about $\hat\lam$ is, by definition,
\beq\label{eq:moment}
M_p(\hat\lam)=\frac{1}{2\pi\mi}\int_\gamma \frac{(\lam-\hat\lam)^pD'(\lambda)}{D(\lam)}\,\dif \lam
=\sum_{k=1}^m(\lambda_k-\hat\lam)^p\,.
\eeq
In particular, $M_0(0)=m$ gives the number of zeros inside $\gamma$ while $M_1(0)=\sum_{k=1}^m\lambda_k$ returns the sum of the zeros. Thus, we may determine the number and location of roots by approximating the integral in \eqref{eq:moment} via Simpson's rule. 
We note that unlike Newton-based, iterative methods for root finding, this approach requires neither a good initial guess for the root location nor excessive evaluations of the function $D$. These features make this approach attractive, and they provide justification for the care taken above to ensure analyticity. For the RNS system, we use the method of moments in conjunction with a two-dimensional bisection method to locate zeros of $D$. 
\end{description}

\section{Results: viscous hyperstabilization}\label{sec:results}

We now describe our principal experiments. For various families of viscous strong detonation waves we compute the zeros of the associated Evans functions in the large-\Ea\ limit. Plots of this spectral information, shown together with the spectral data from the corresponding inviscid ZND profiles gives a direct way to visualize the impact of the viscous terms in the RNS model on the stability of these waves. 
This is shown, for example, in Figures \ref{fig:PowersEA}, \ref{fig:025}, and \ref{fig:025-2}. Initally, we find a crossing of a complex conjugate pairs of eigenvalues into the unstable half plane; this signal of a Hopf-type bifurcation is  completely expected. However,  continuing to follow the eigenvalues as \Ea\ increases, we find an unexpected \emph{restabilization}. We describe our observations of the behavior of the spectra in more detail below. 

\subsection{The viscosity of Romick et al.\ \cite{RAP_JFM12}}
\begin{figure}[p] 
   \centering
   \begin{tabular}{ccc}
   (a) \includegraphics[width=1.6in]{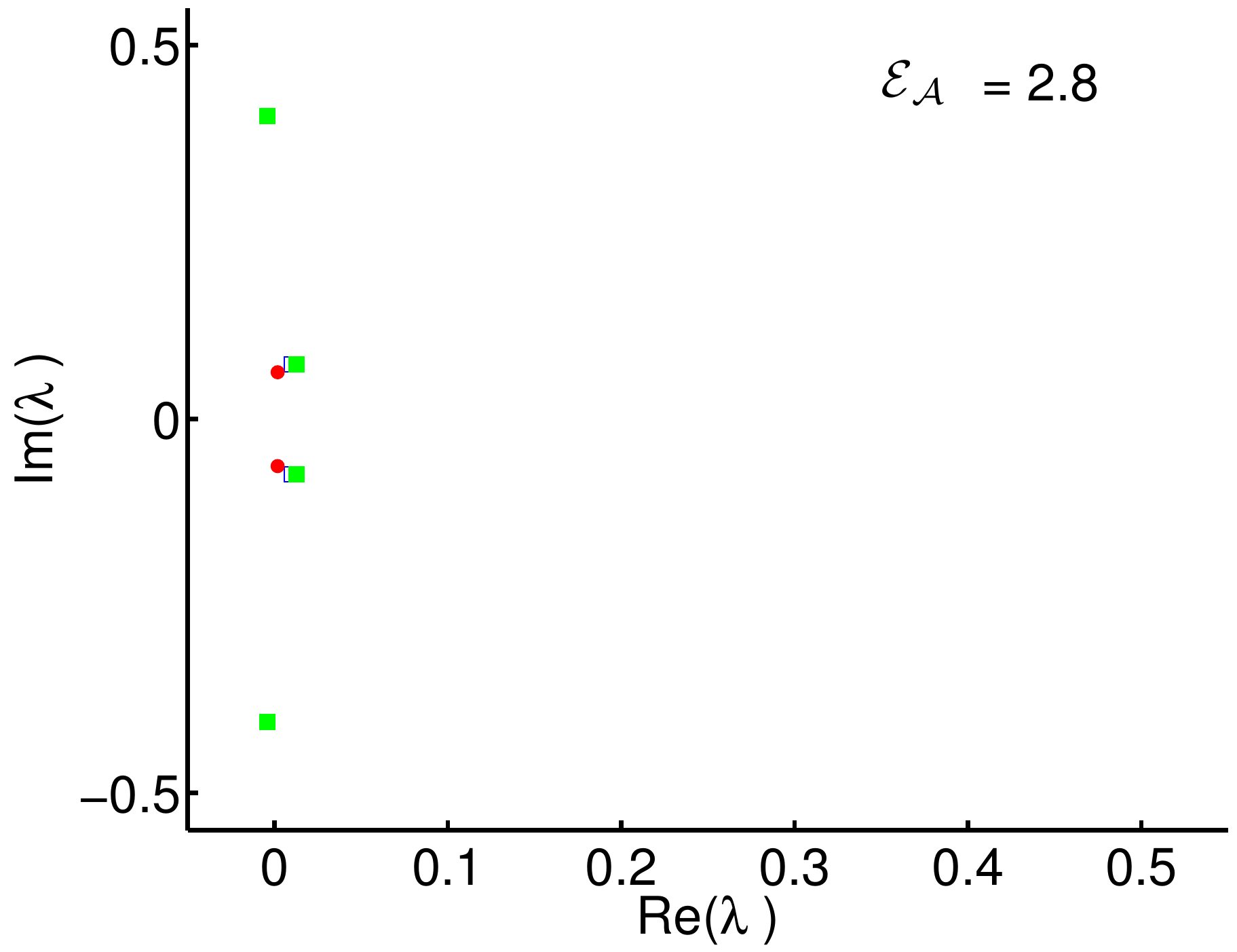} &  
   (b) \includegraphics[width=1.6in]{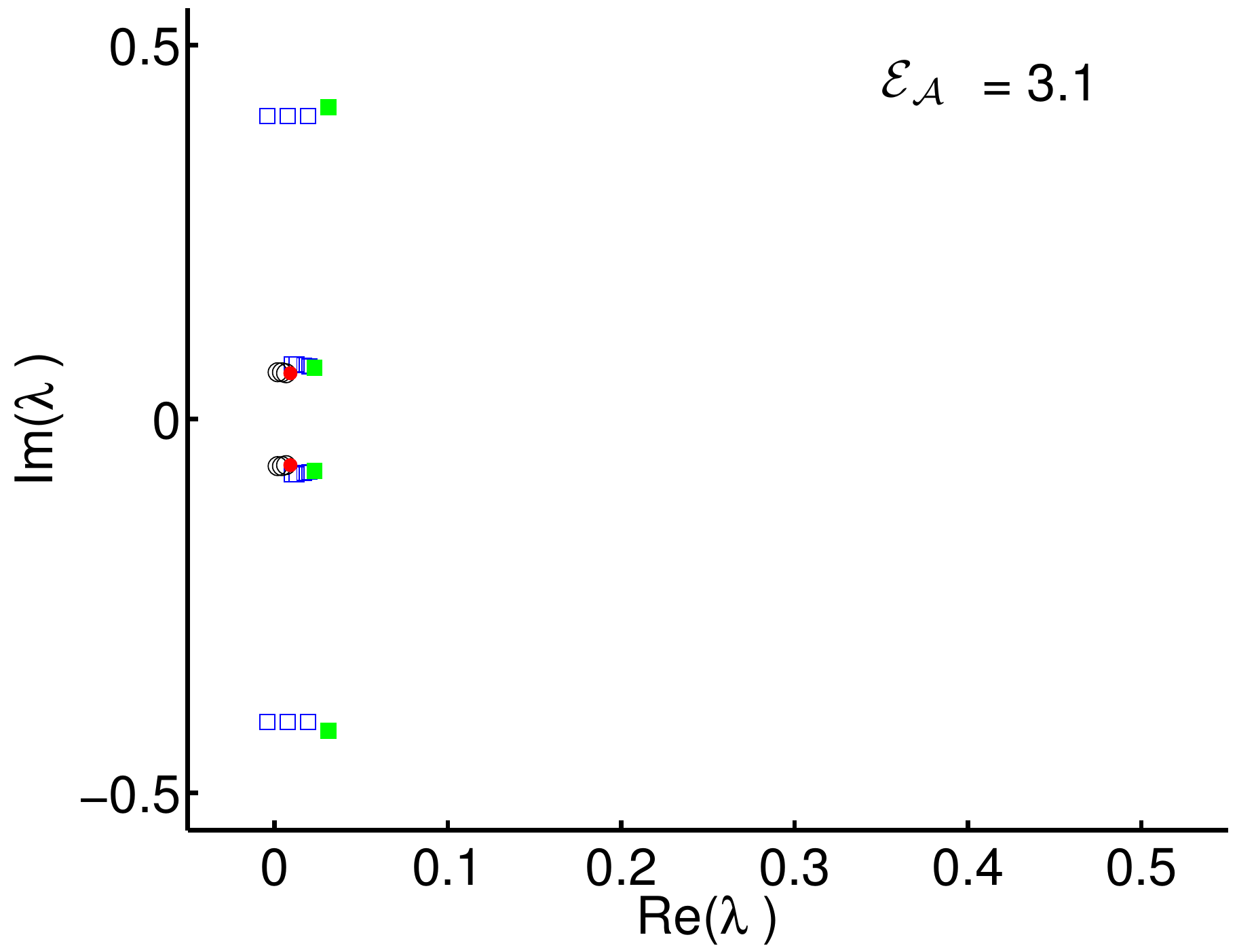} &  
   (c) \includegraphics[width=1.6in]{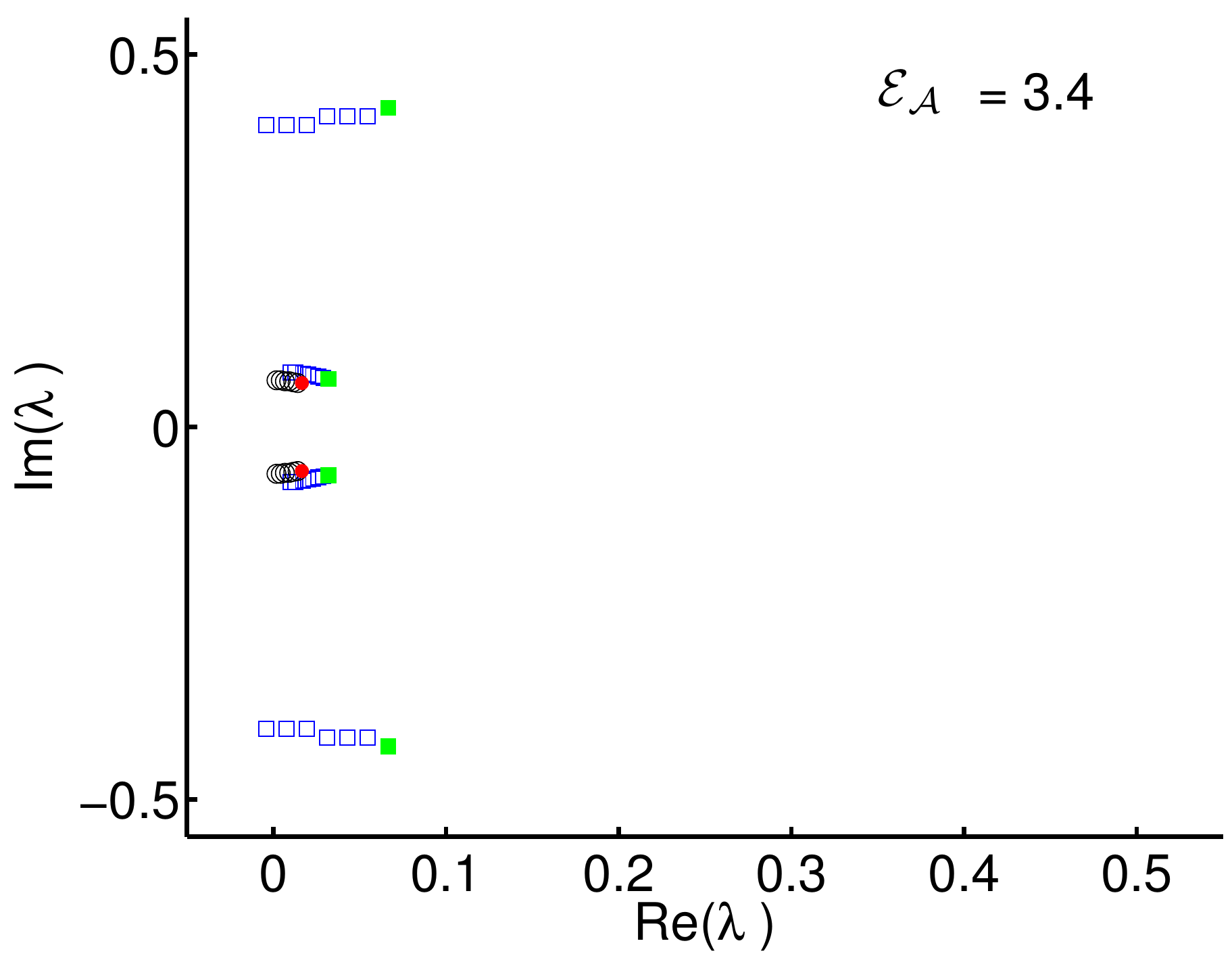} \\
   (d) \includegraphics[width=1.6in]{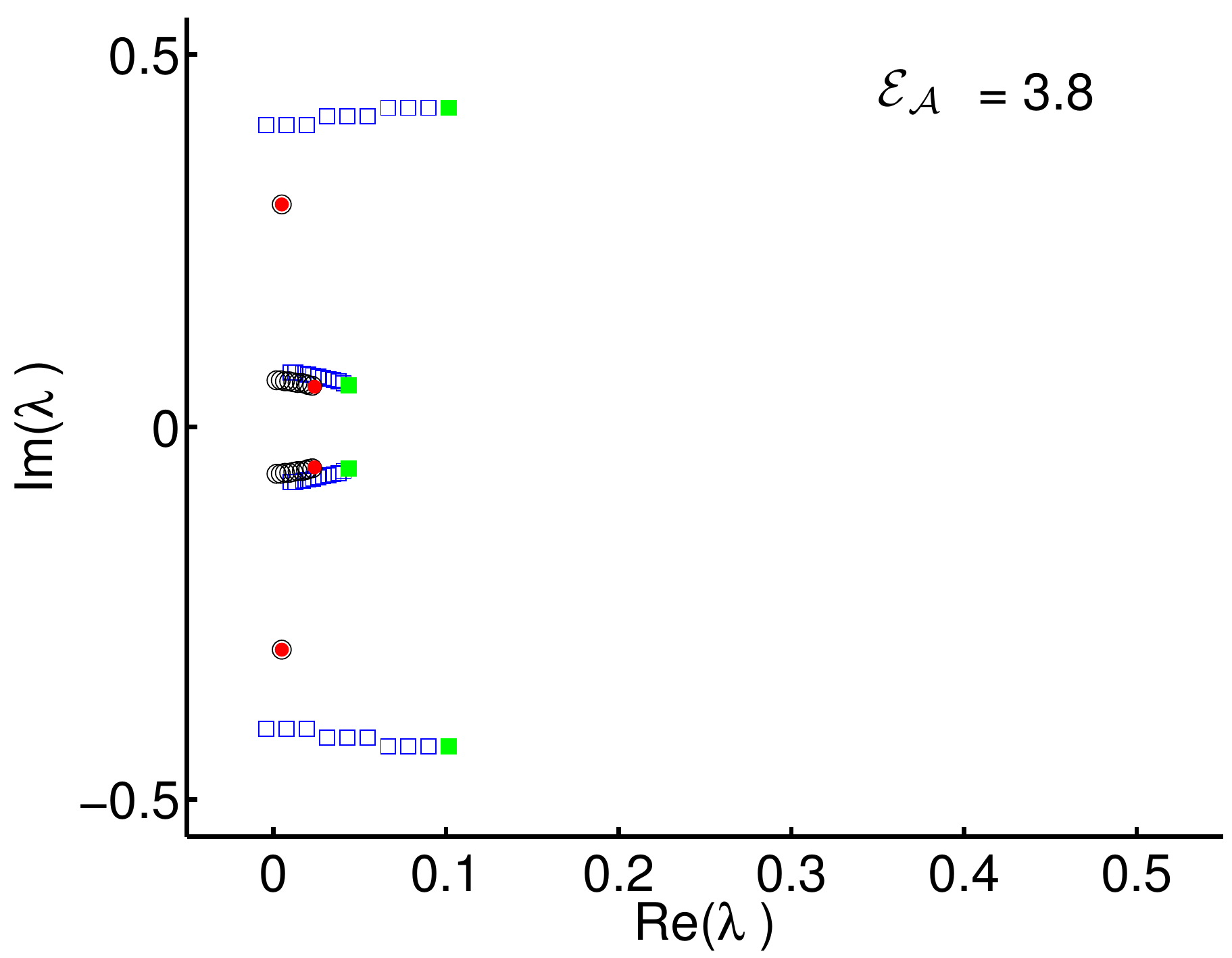} &  
   (e) \includegraphics[width=1.6in]{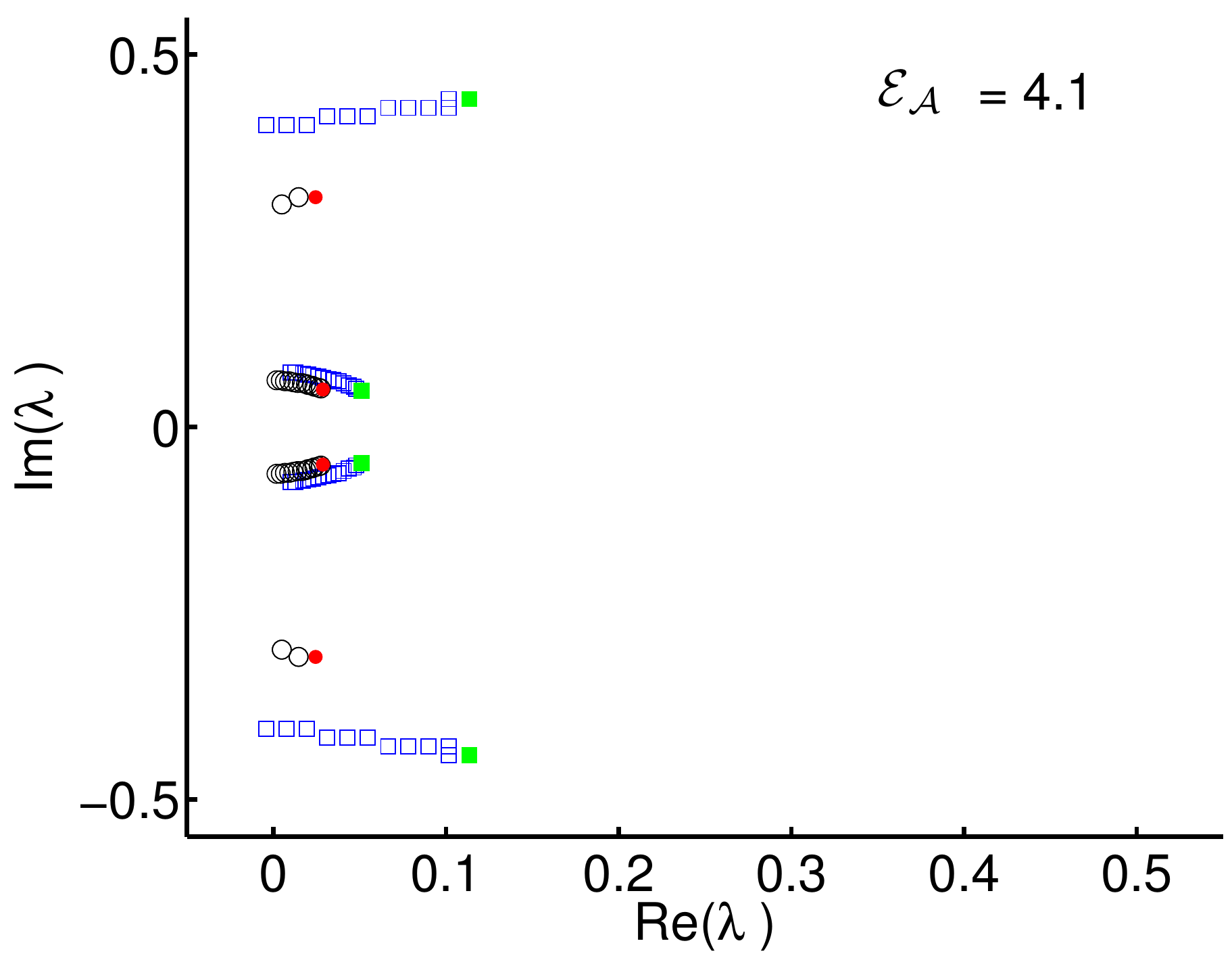} &  
   (f) \includegraphics[width=1.6in]{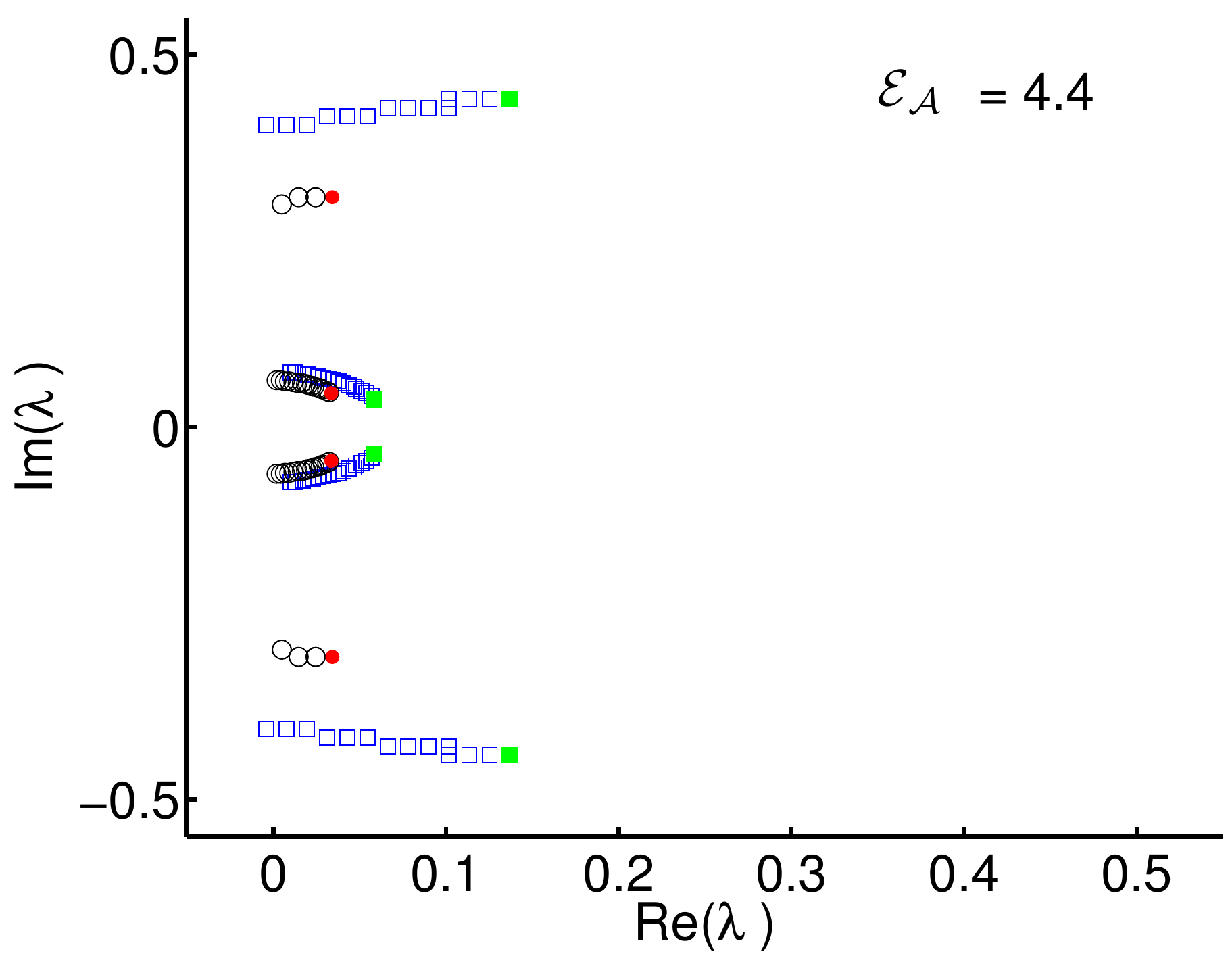} \\
   (g) \includegraphics[width=1.6in]{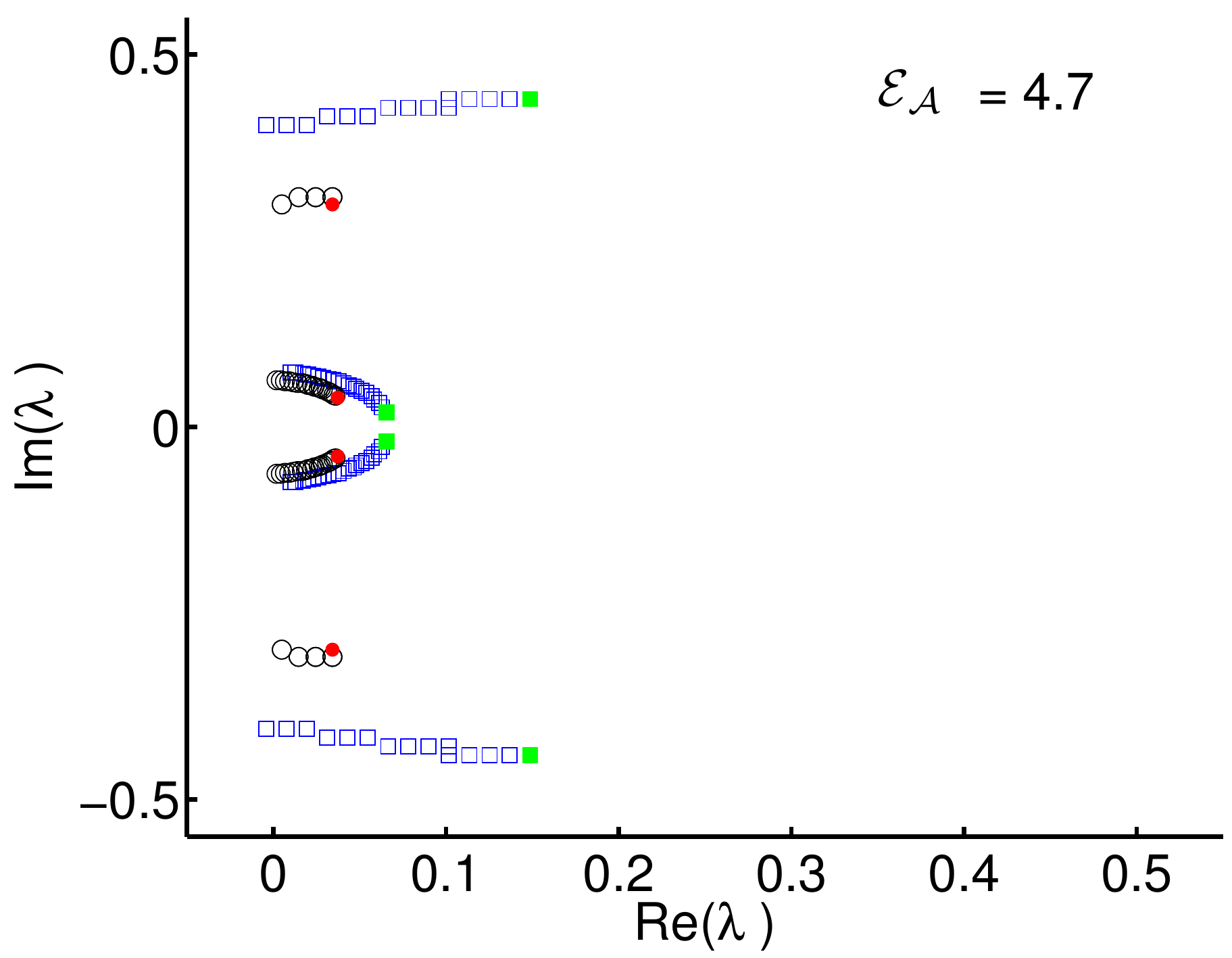} &  
   (h) \includegraphics[width=1.6in]{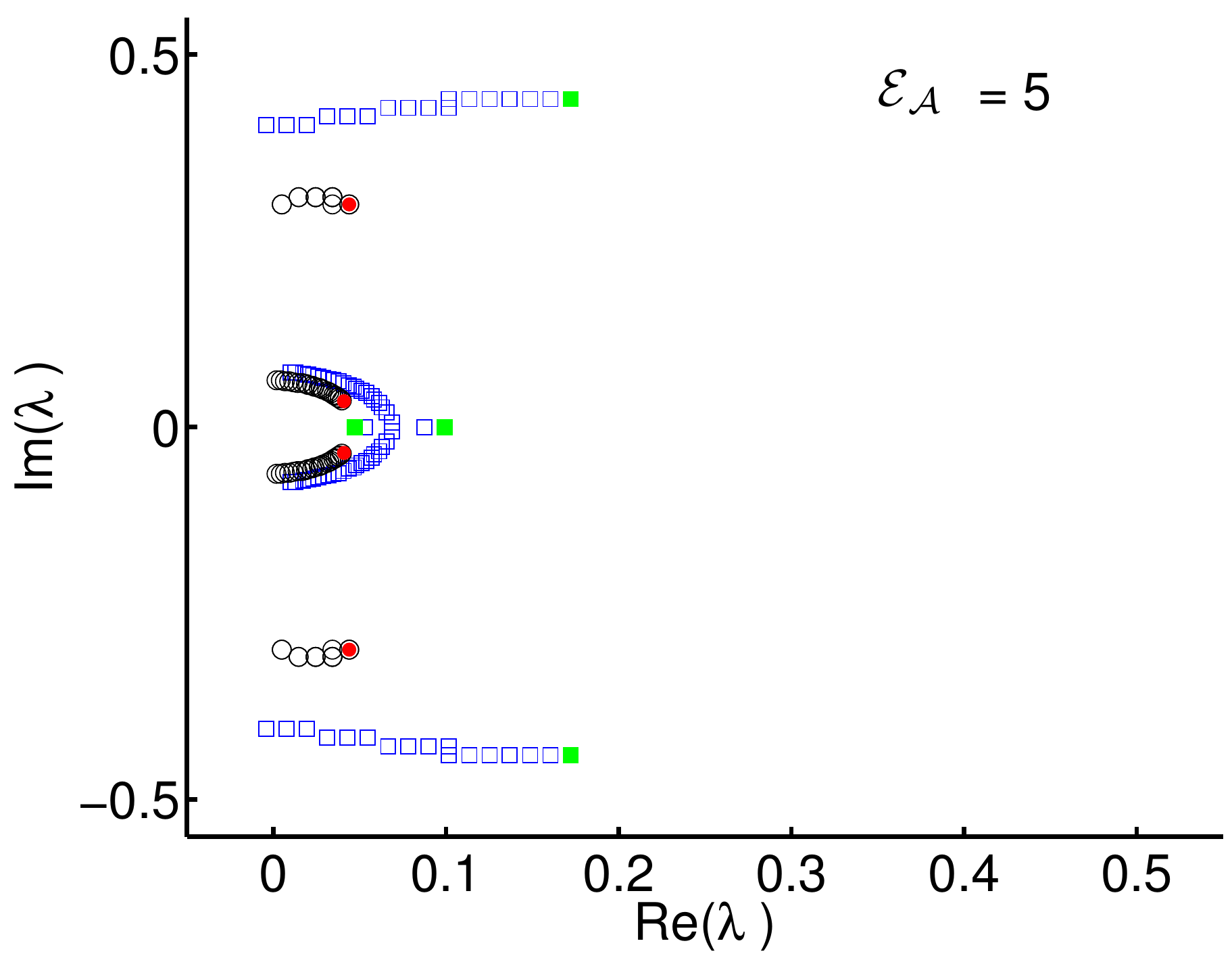} &  
   (i) \includegraphics[width=1.6in]{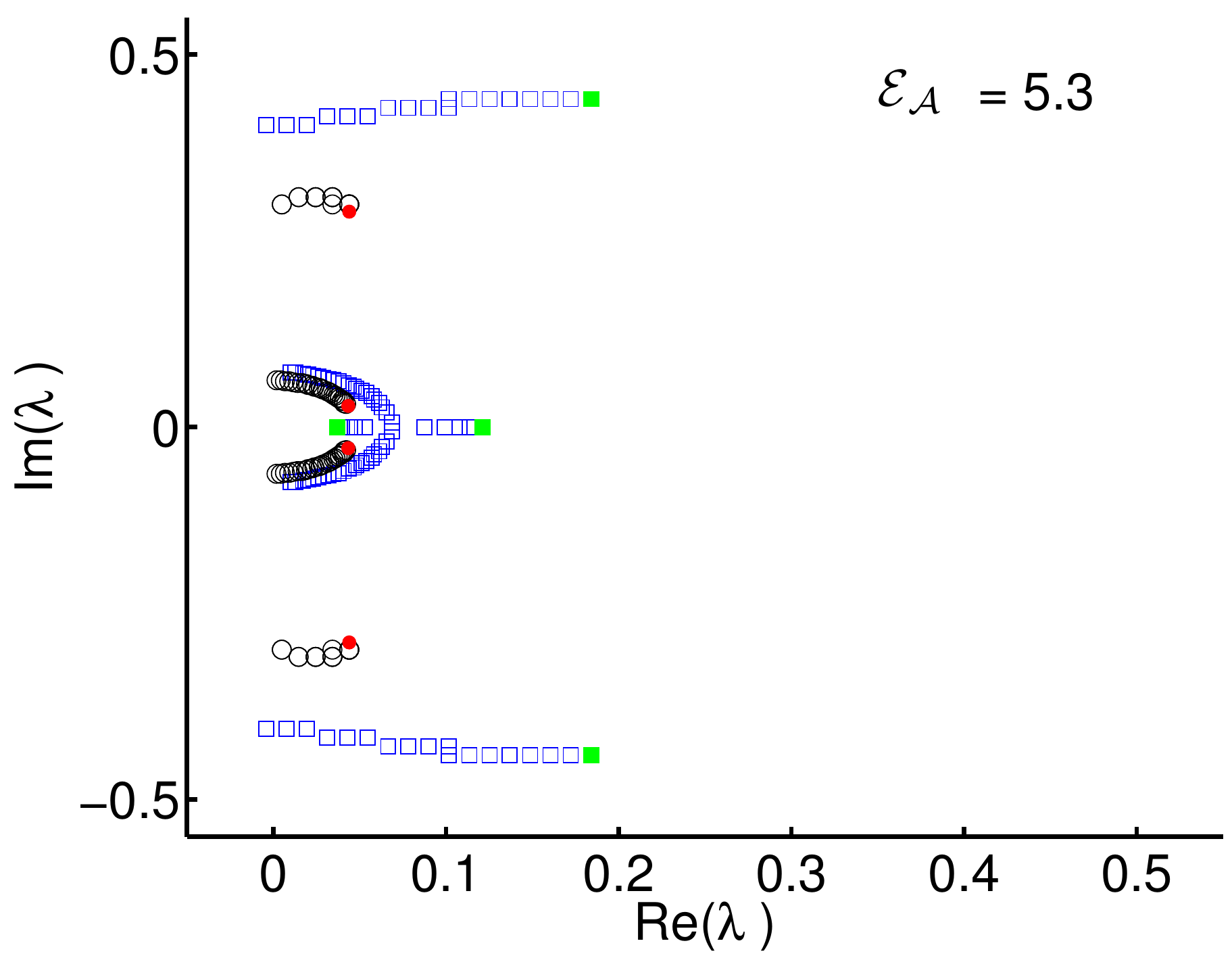} \\
   (j) \includegraphics[width=1.6in]{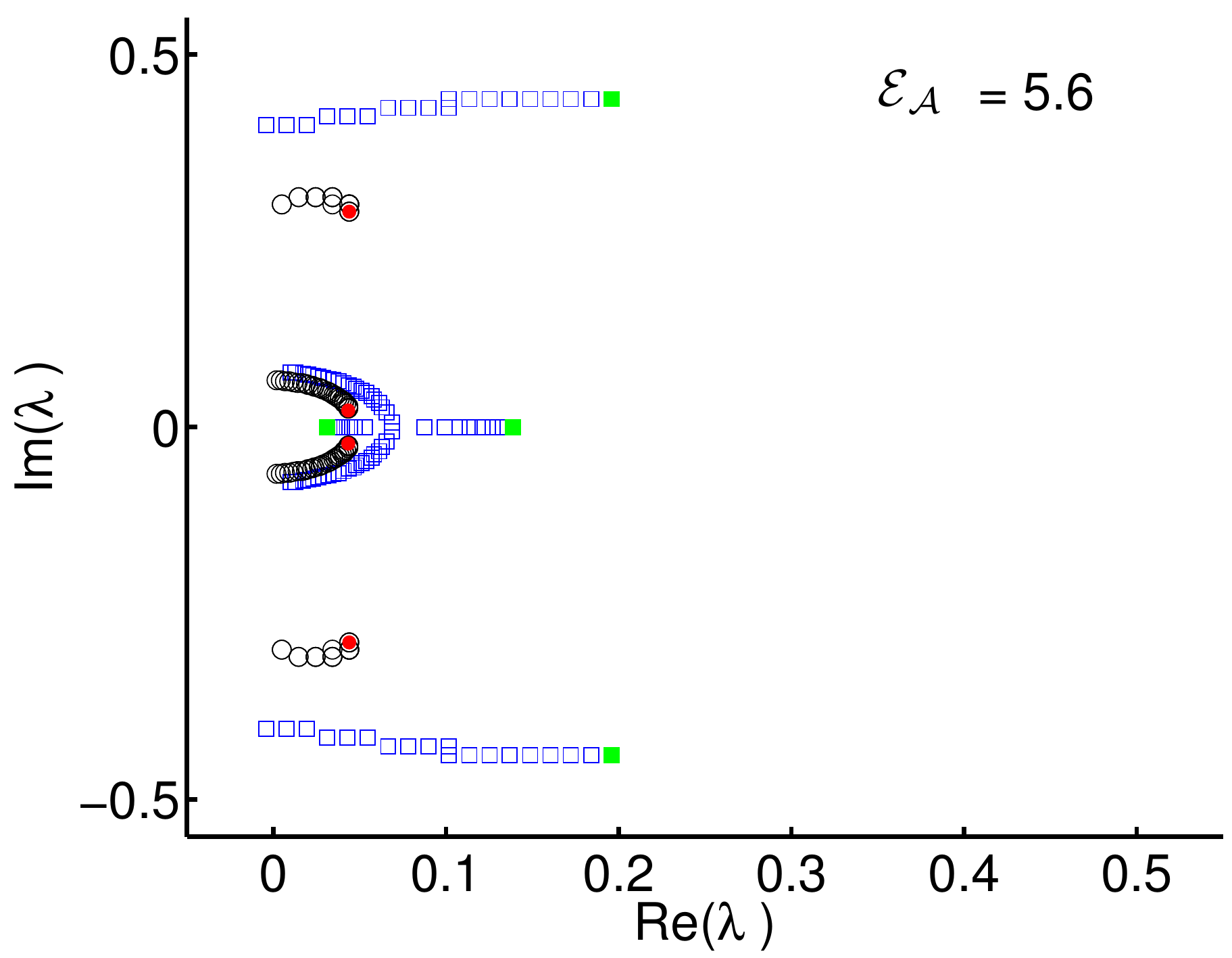} &  
   (k) \includegraphics[width=1.6in]{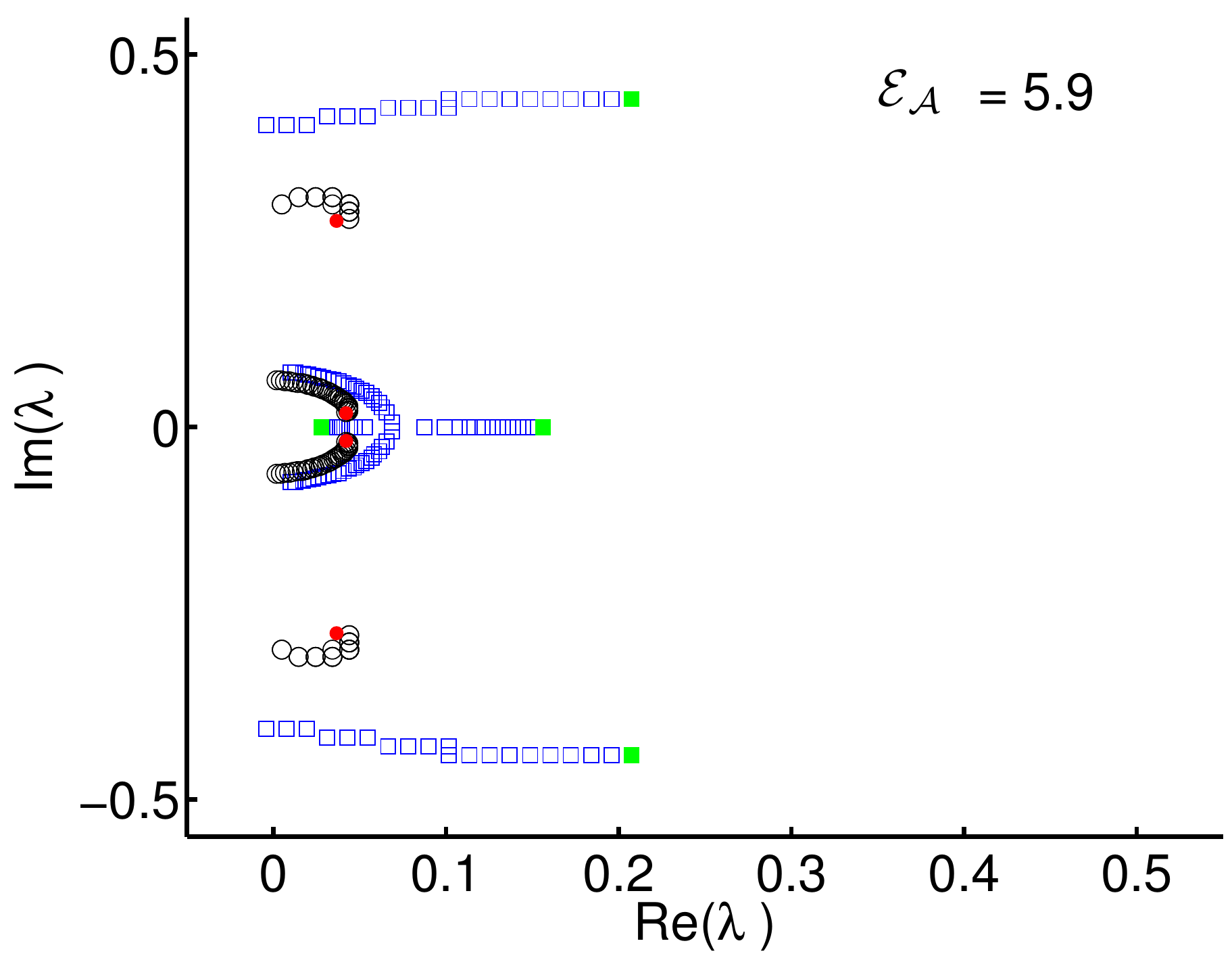} &  
   (l) \includegraphics[width=1.6in]{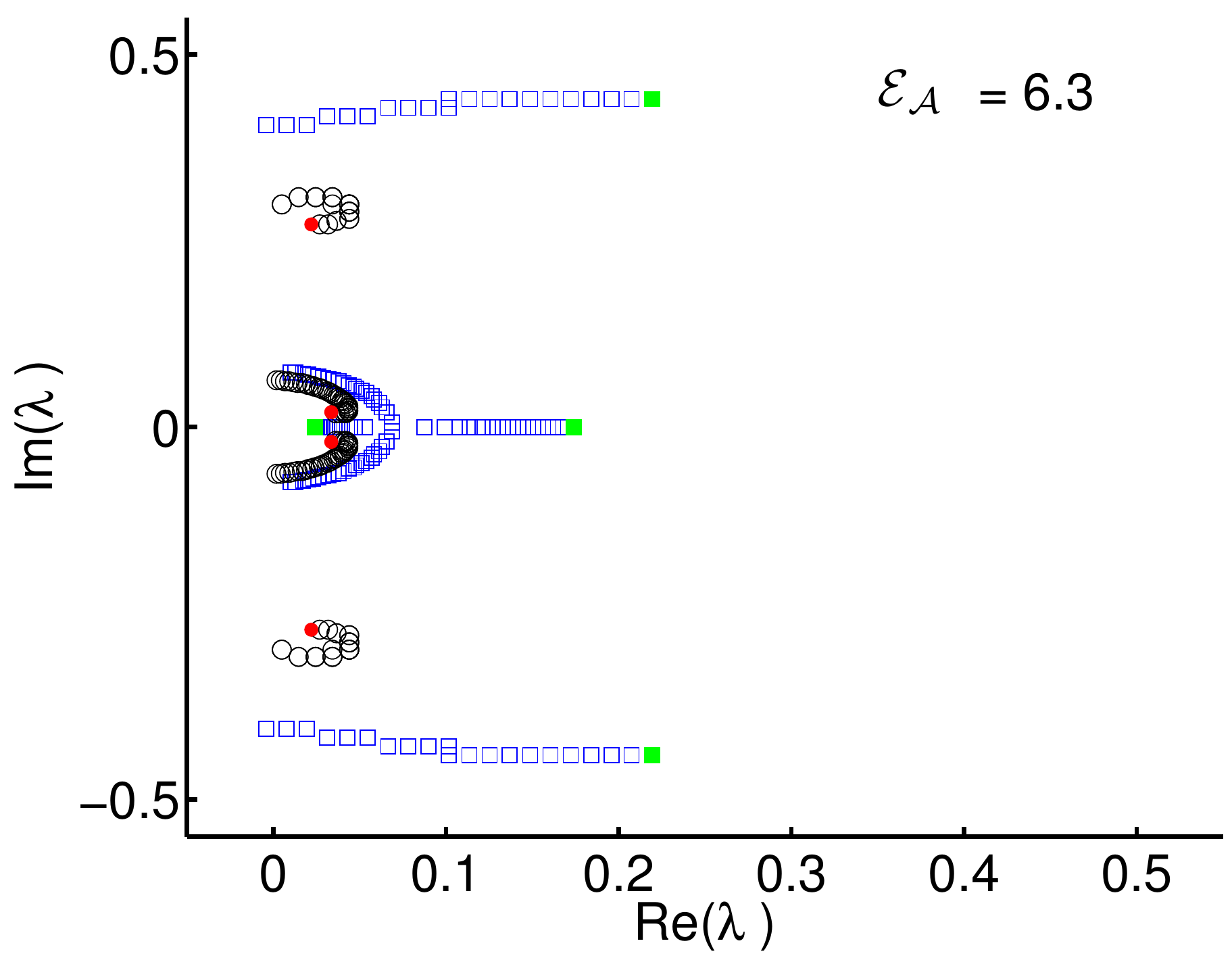} \\
   (m) \includegraphics[width=1.6in]{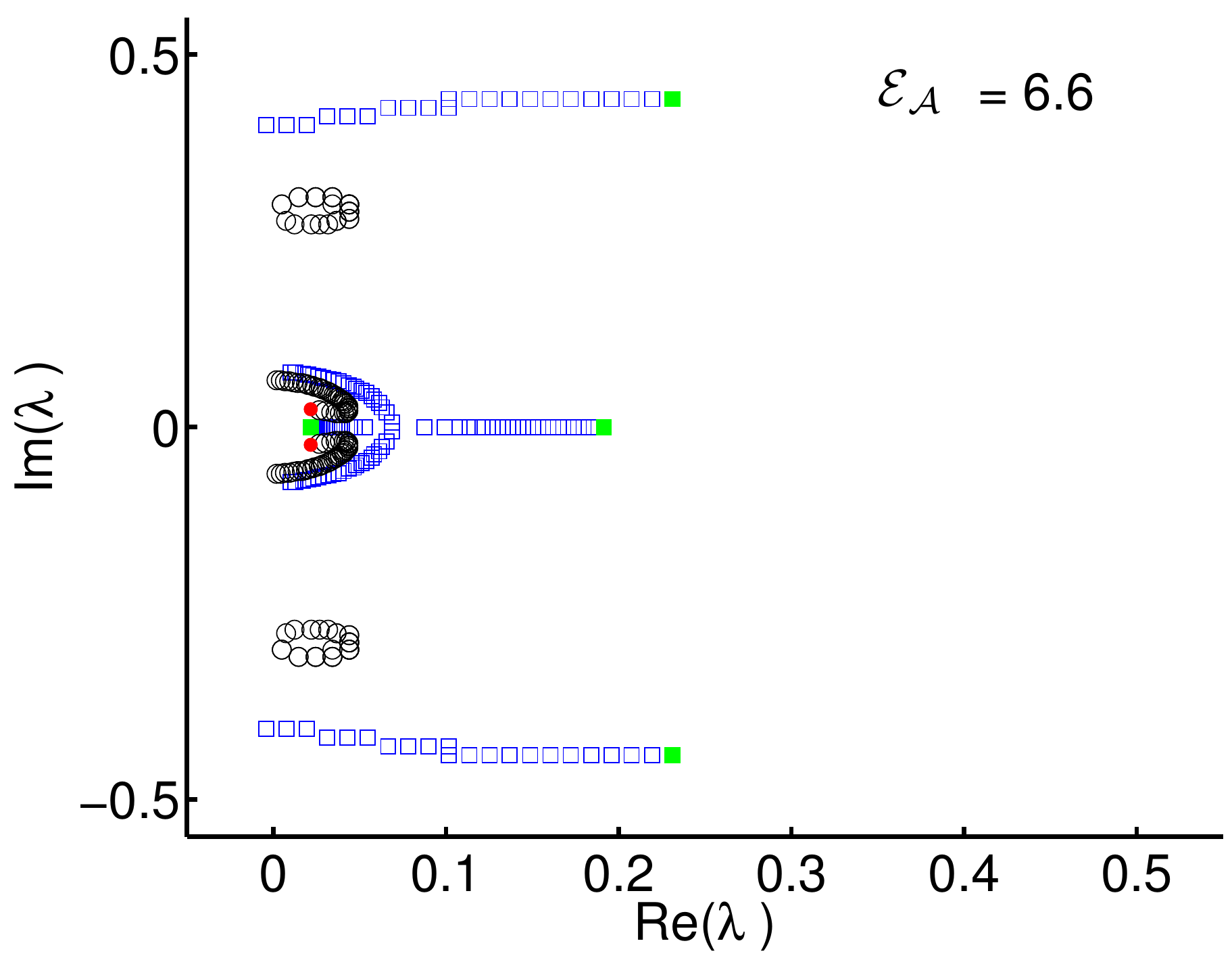} &
   (n) \includegraphics[width=1.6in]{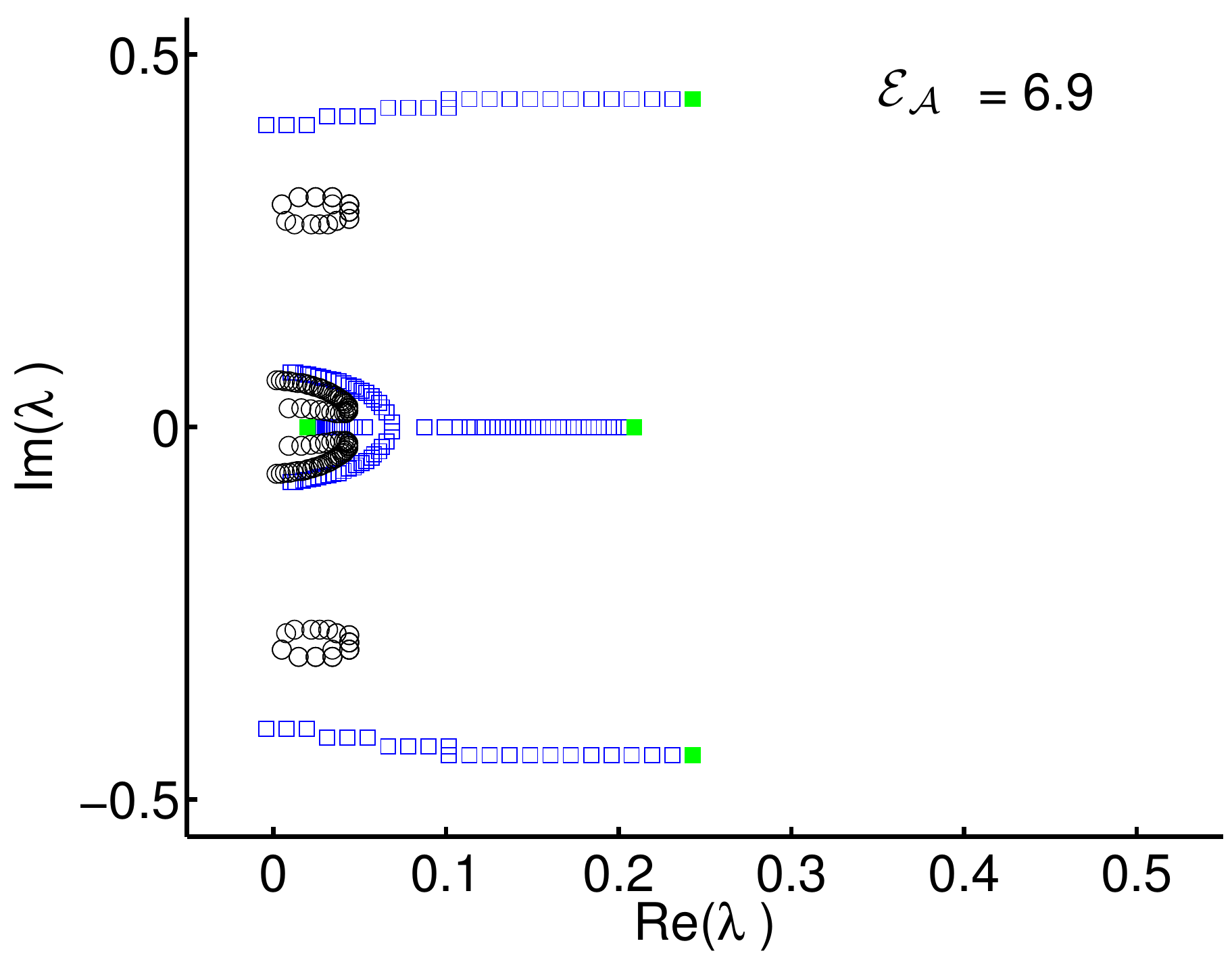} & 
   (o) \includegraphics[width=1.6in]{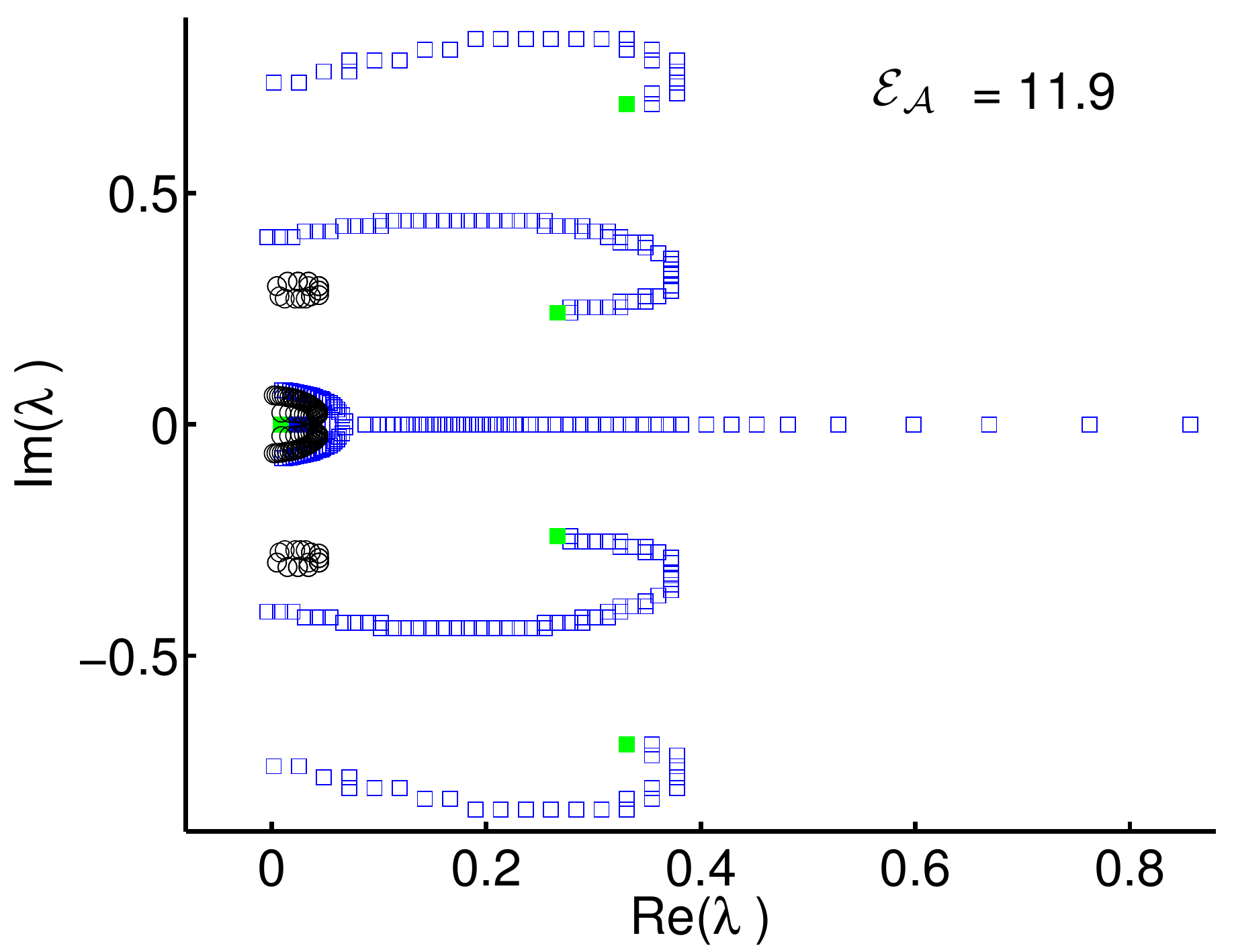} \\
   \end{tabular}
  \caption{The movement of unstable roots in the complex plane as $\Ea$ increases. Circles (solid dot for current value of $\Ea$) mark the roots corresponding to the RNS model, while open squares (solid square for the current value of $\Ea$) correspond to the ZND model.     The smaller modulus roots enters for $\Ea \approx 0.29$ (Panel (a)) and the higher modulus roots for $\Ea \approx 3.7$ (Panel (d)). The high modulus roots have a turning point at about $\Ea \approx 5.2$, and the smaller modulus roots have theirs at $\Ea \approx5.5$ (Panels (i) and (j)). The large modulus roots leave at $\Ea\approx 6.4$ and the small modulus roots leave at approximately $\Ea\approx6.85$ (Panel (n)). The other parameters are  $\Gamma=0.2$; $e_\sp$=6.23e-2; $q$= 6.2e-1; $d=0.1$; $\kappa= 0.1$; $\nu= 0.1$; $k$ chosen as described in  \S\ref{ssec:parametrization}; $\tau_\sm$= 2.57e-1; $u_\sm$= 7.43e-1; $e_\sm$=9.71e-1; $c_v=1$; $z_\sp=1$; $z_\sm= 0$; $\tau_\sp=1$; $u_\sp=0$; $s=1$; $y_\spm=0$; $T_\sp$=6.2e-2; $T_\sm$= 9.71e-1; and $\Ti$=6.64e-2.} 
\label{fig:PowersEA}
\end{figure}

\begin{figure}[ht] 
   \centering
   \includegraphics[width=3in]{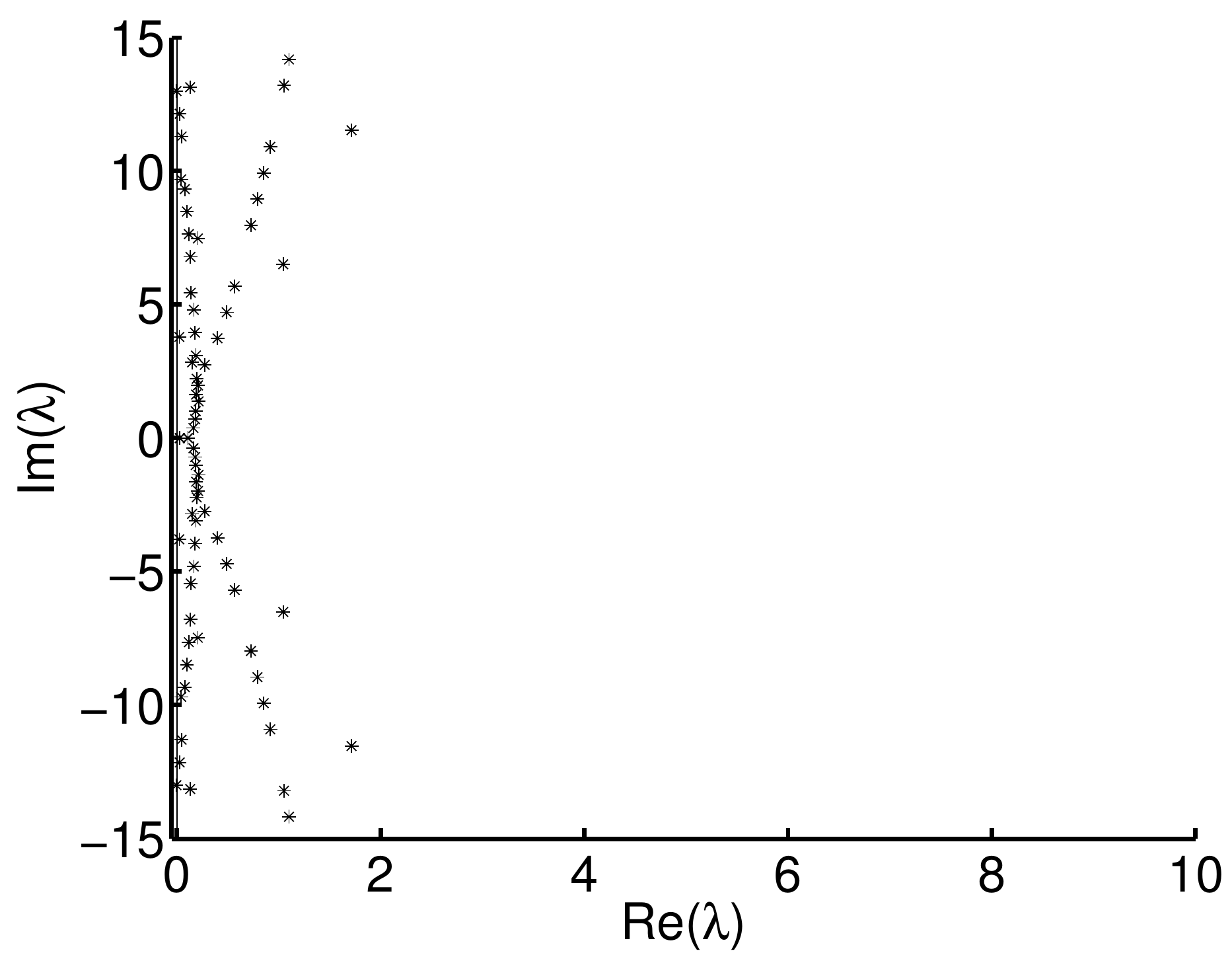} 
   \caption{Unstable eigenvalues in the complex plane for the inviscid ZND problem. The activation energy is $\Ea=7.1$ and all of the other parameter are as in Figure \ref{fig:PowersEA}. There are 49 roots in the picture.  This figure is taken from \cite{BZ_ZND} (with permission of the authors); details of this computation can be found in \cite{BZ_ZND}.}
   \label{fig:ZNDCompare}
\end{figure}

In this first experiment, we use the viscous parameter values utilized in the recent numerical study of Romick et al.\ \cite{RAP_JFM12}. Figure \ref{fig:PowersEA} shows that the smaller modulus roots enters for $\Ea \approx 0.29$ (Panel (a)) and the higher modulus roots for $\Ea \approx 3.7$ (Panel (d)). The high modulus roots have a turning point at about $\Ea \approx 5.2$, and the smaller modulus roots have theirs at $\Ea \approx5.5$ (Panels (i) and (j)). The large modulus roots leave at $\Ea\approx 6.4$ and the small modulus roots leave at approximately $\Ea\approx6.85$ (Panel (n)). 
Indeed, at the spectral level the contrast with ZND is striking; see Figure \ref{fig:ZNDCompare} which shows the ZND spectra in  a much larger portion of the unstable half plane than is shown in Figure \ref{fig:PowersEA}. In Figure \ref{fig:ZNDCompare} the value of \Ea\ is 7.1 which corresponds to just after panel (n) in the computation with viscosity. Notably, at this point in the computation with viscosity, the hyperstabilization has already taken place. 

\subsection{Smaller viscosity}
\begin{figure}[ht] 
   \centering
   \begin{tabular}{ccc}
   (a) \includegraphics[width=1.7in]{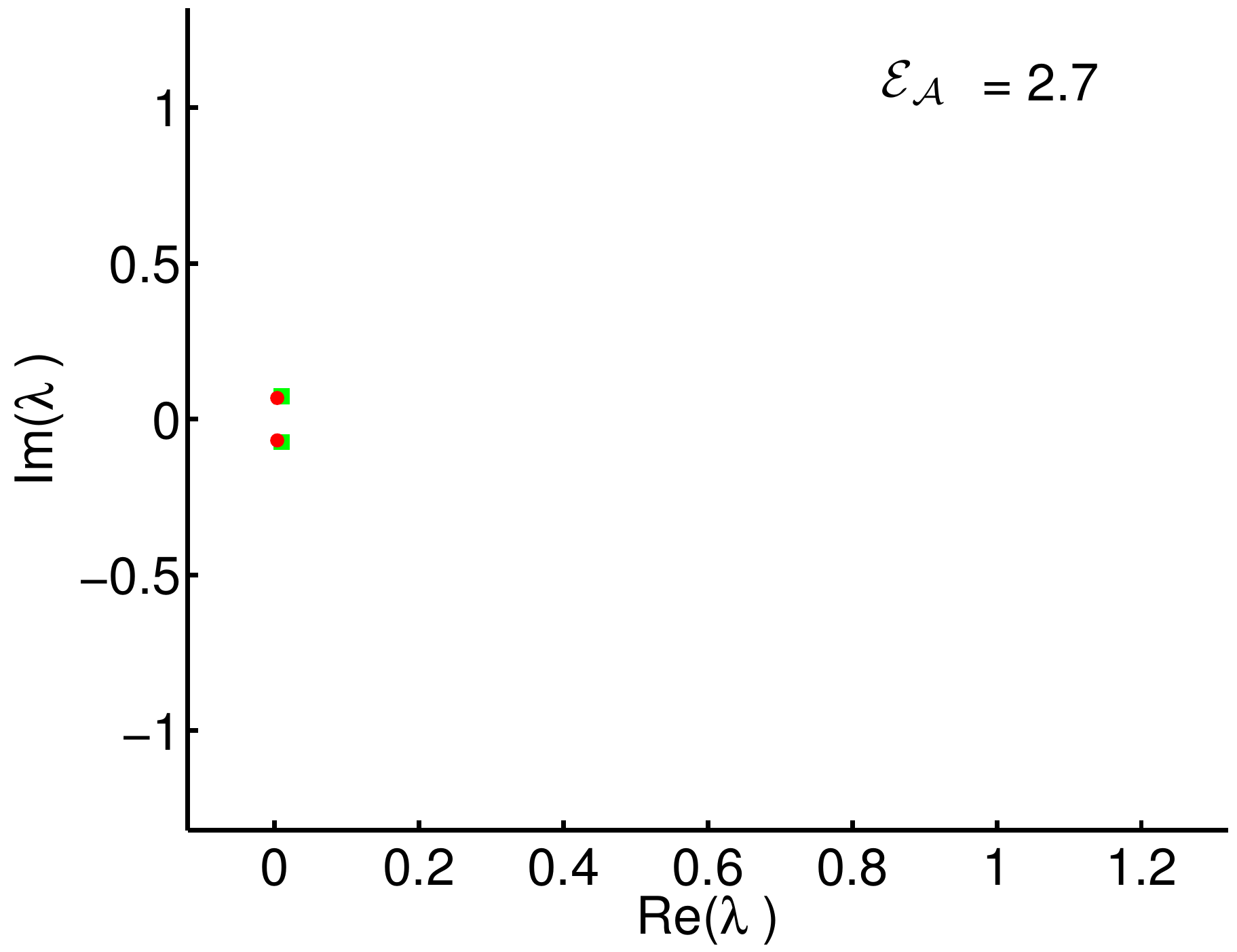} &  
   (b) \includegraphics[width=1.7in]{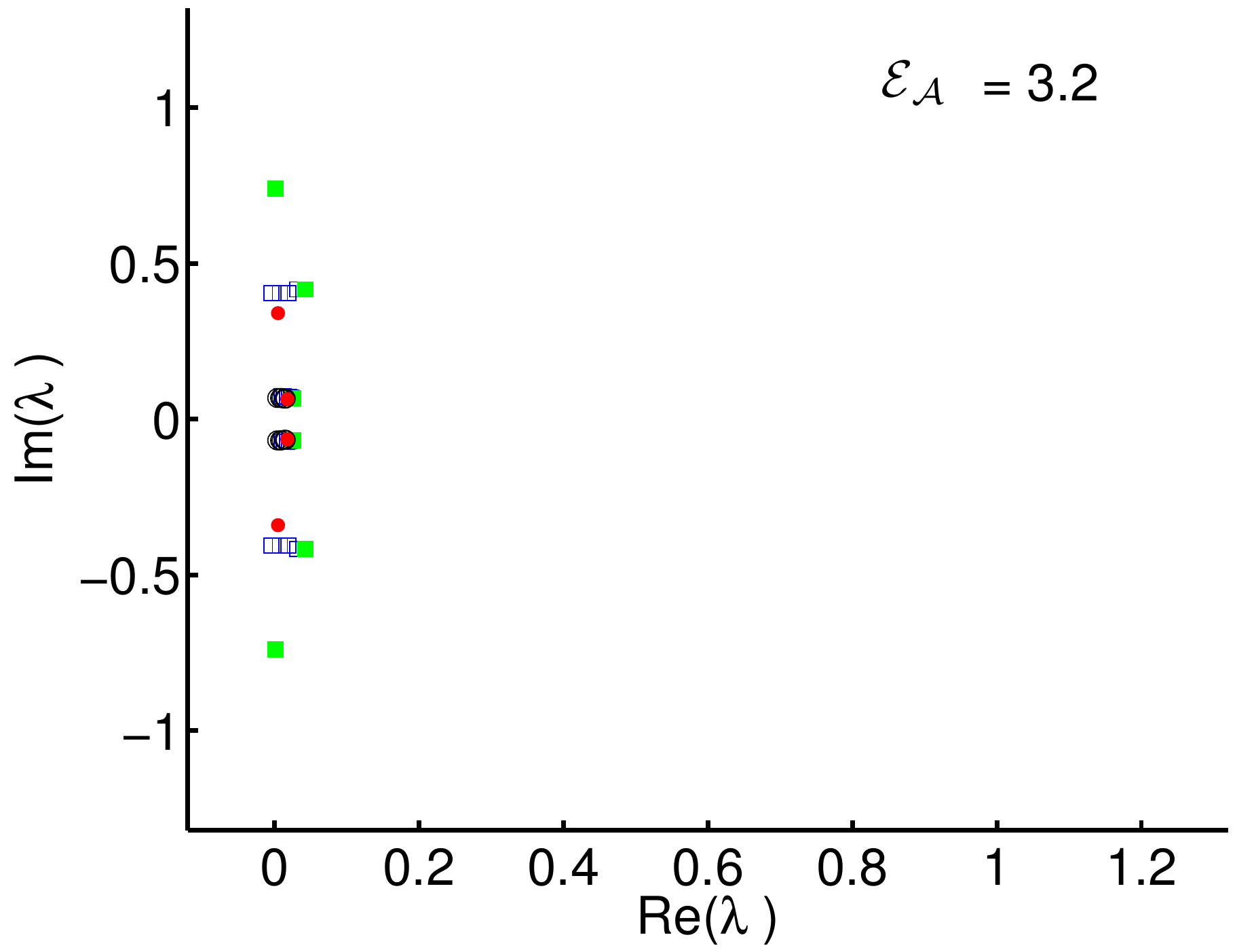} &  
   (c) \includegraphics[width=1.7in]{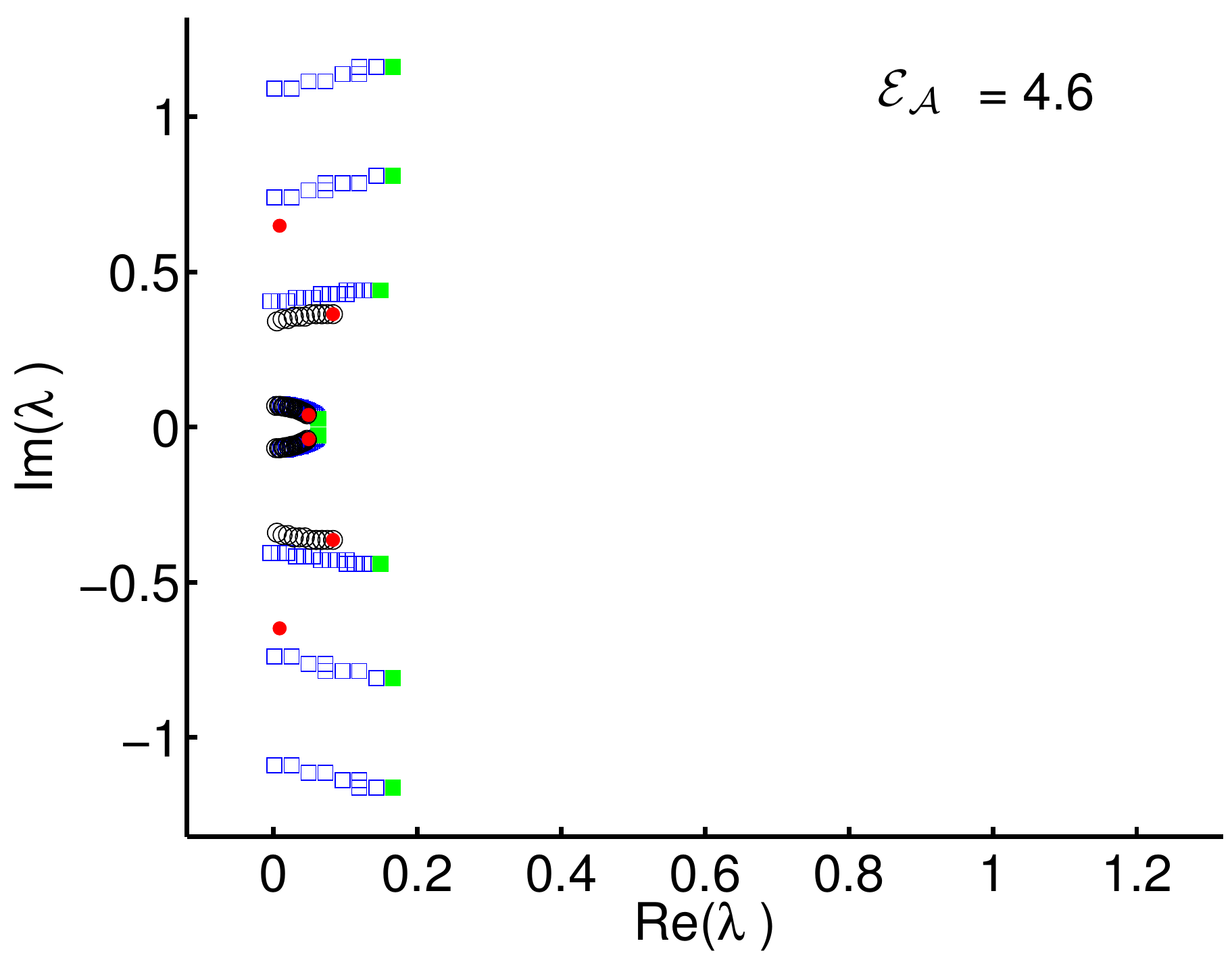} \\
   (d) \includegraphics[width=1.7in]{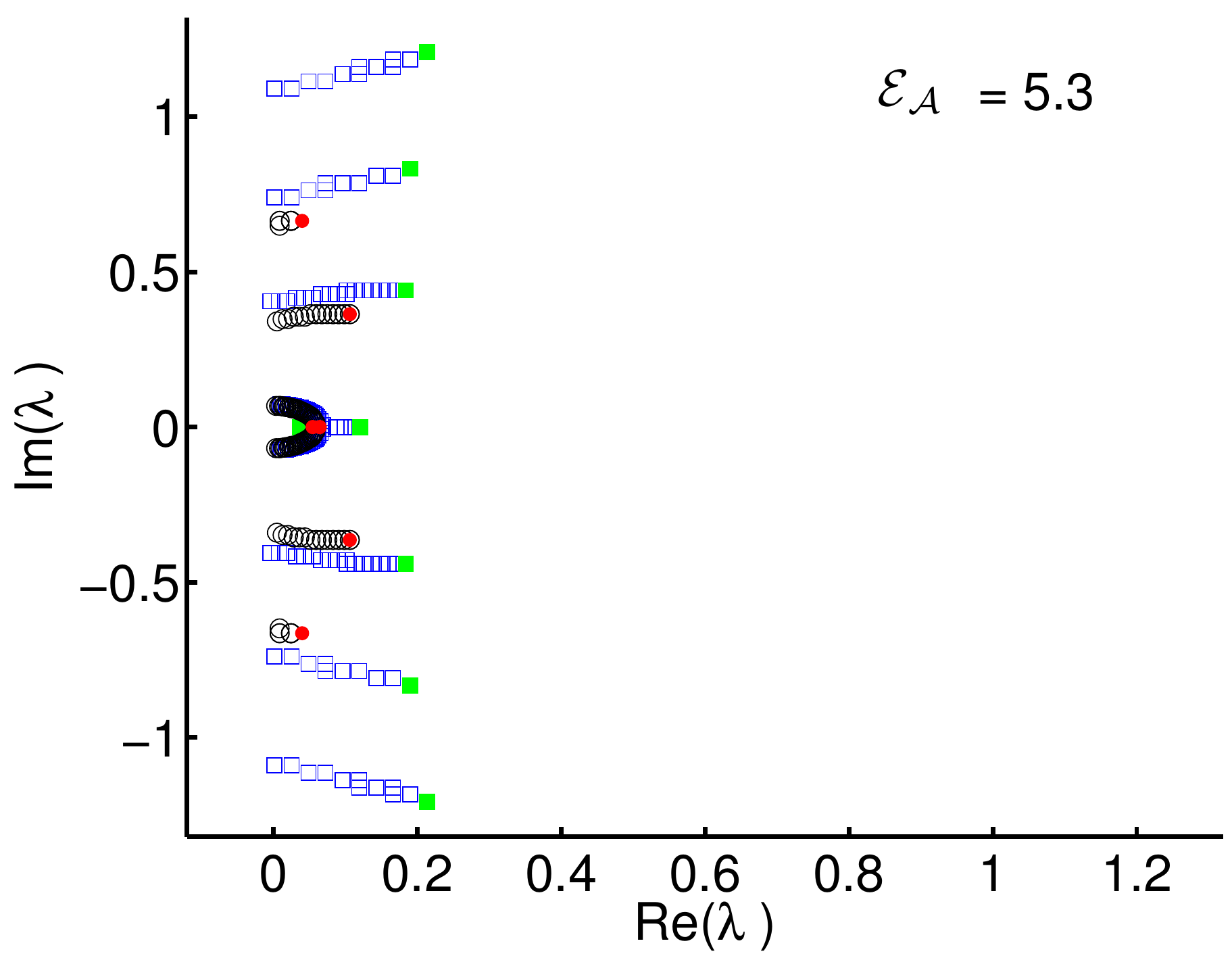} &  
   (e) \includegraphics[width=1.7in]{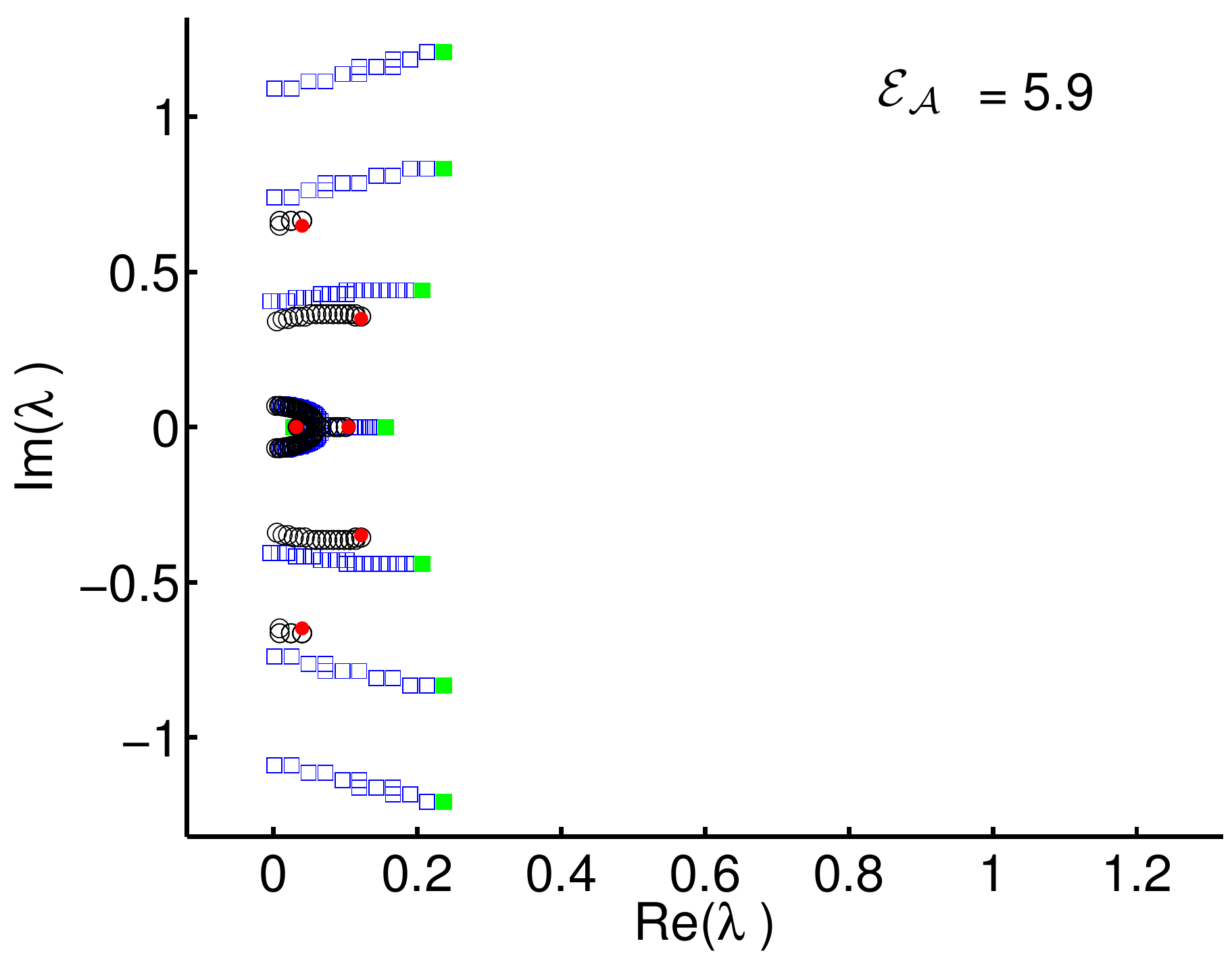} &  
   (f) \includegraphics[width=1.7in]{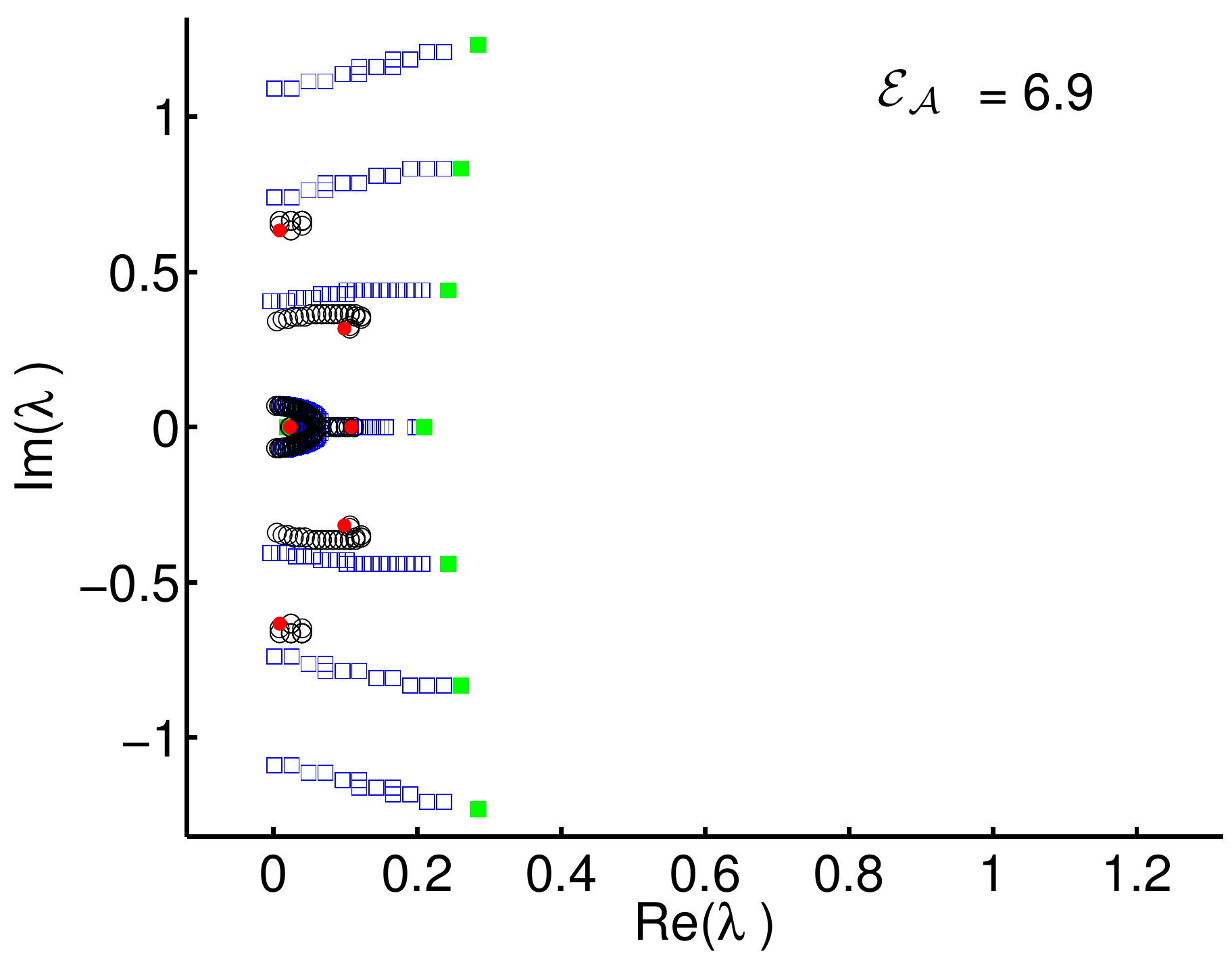} \\
   (g) \includegraphics[width=1.7in]{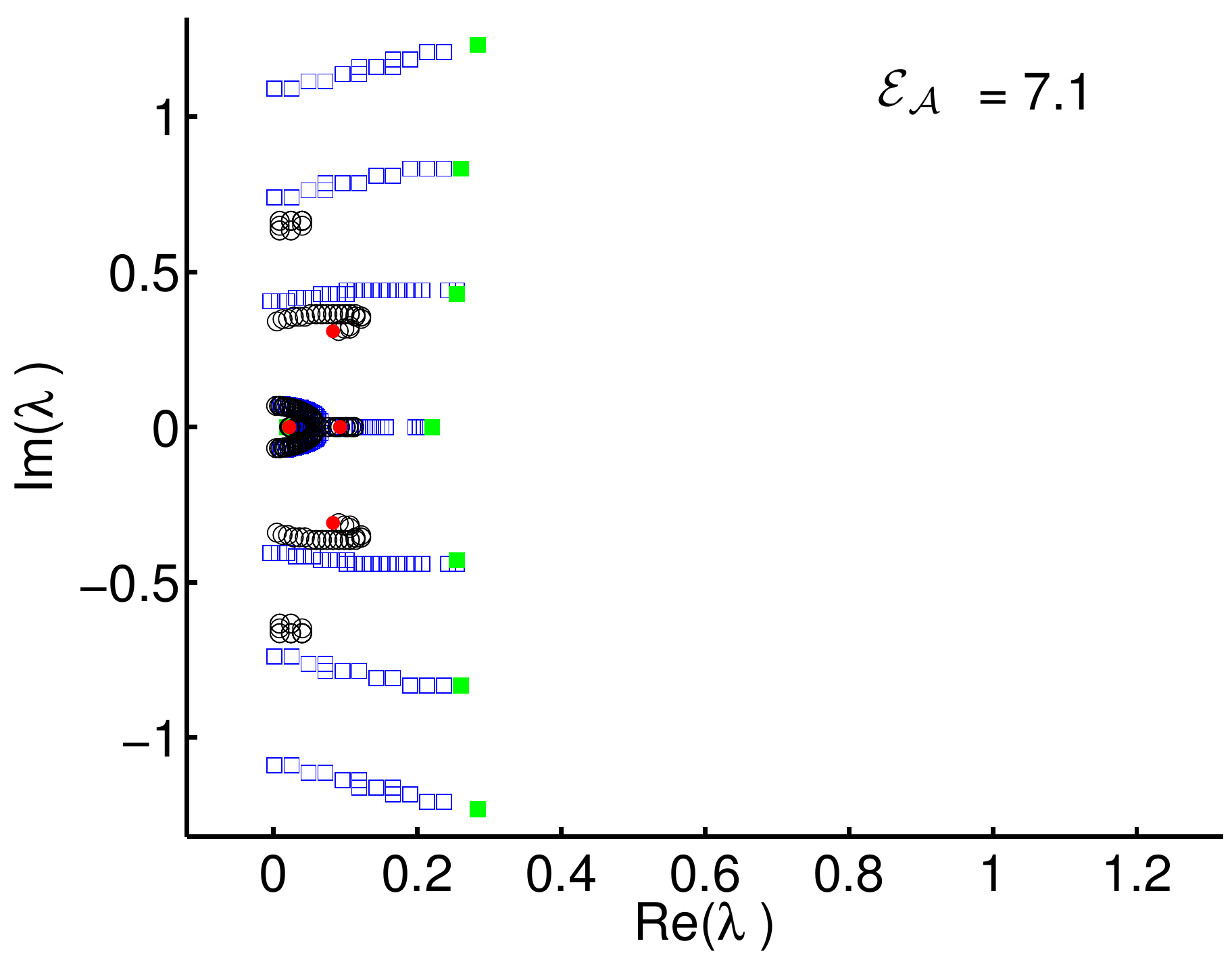} &  
   (h) \includegraphics[width=1.7in]{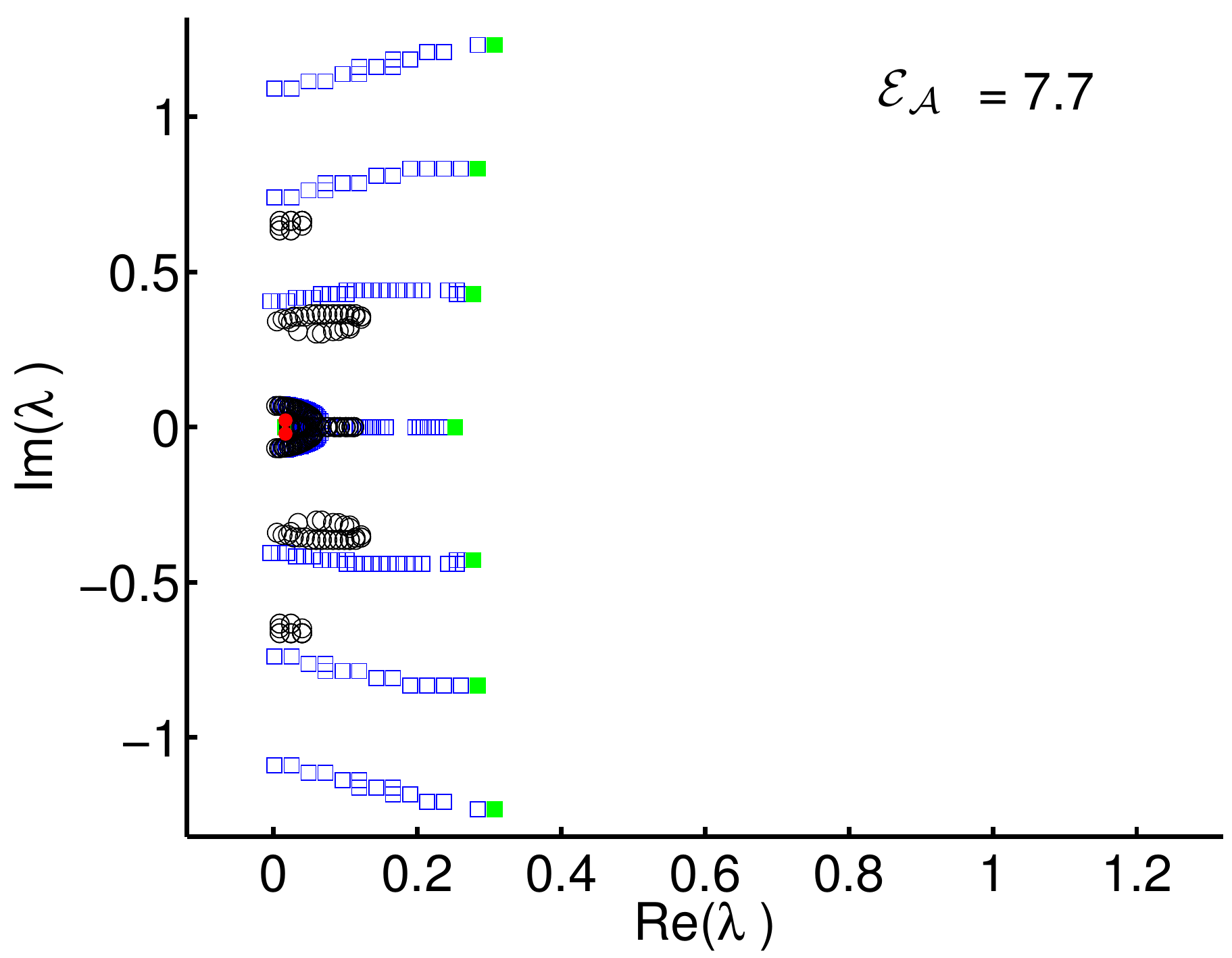} &  
   (i) \includegraphics[width=1.7in]{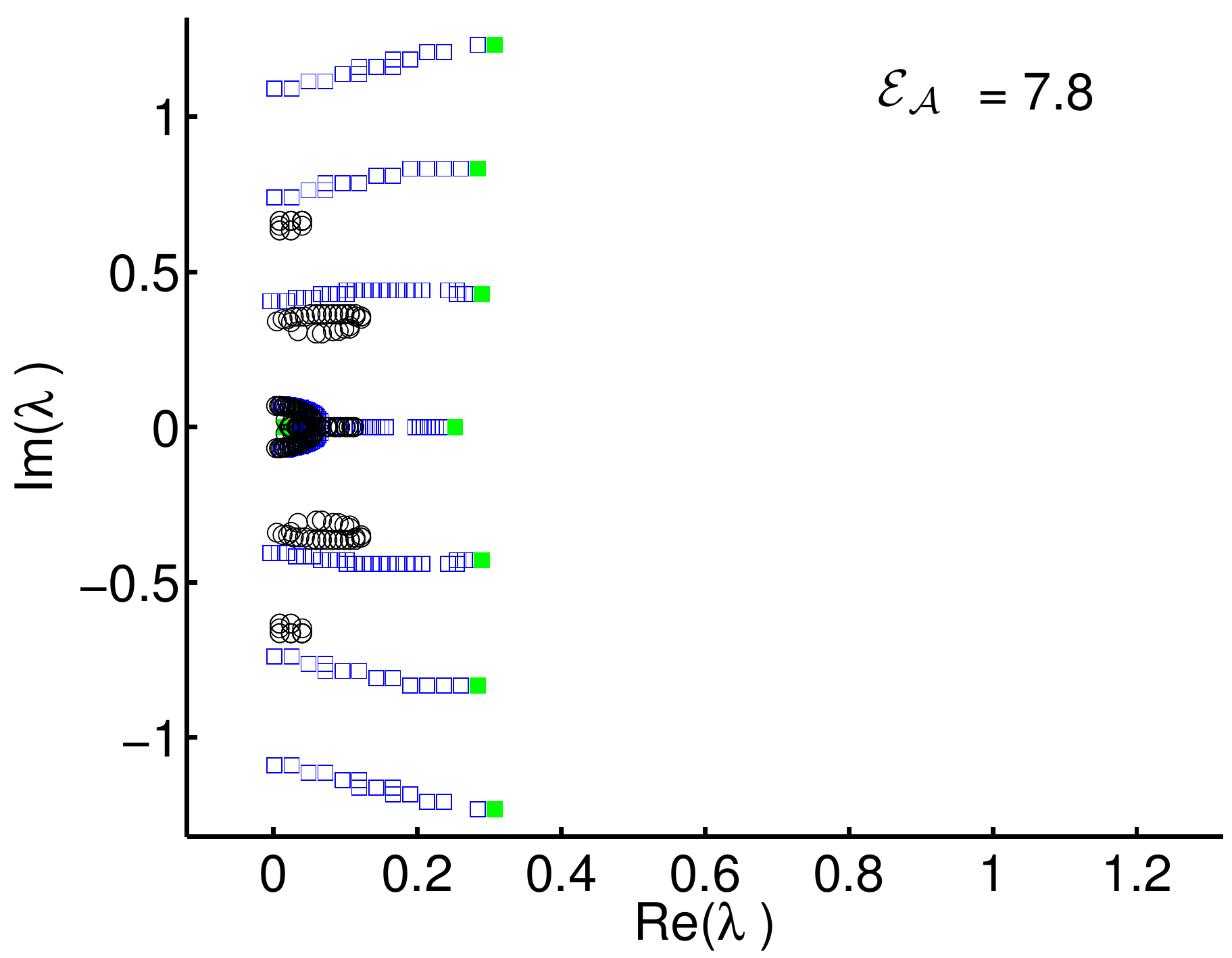} \\
     \end{tabular}
   \caption{The movement of unstable roots in the complex plane as $\Ea$ increases. Open circles (solid dot for current value of $\Ea$) mark the roots corresponding to the RNS, while open squares (solid square for current value of $\Ea$) correspond to the ZND model. The roots with the smallest modulus when they initially cross the imaginary axis cross at $\Ea \approx 2.7$, then collide at $\Ea \approx 5.25$, then at $\Ea \approx 6.55$ the largest real root turns back, at $\Ea \approx 7.5$ the two real roots collide again, and at $\Ea \approx 7.75$ they cross back to the stable side of the imaginary axis. The second largest modulus roots at time of initial crossing cross at $\Ea \approx 3.15$, turn back at $\Ea \approx 5.85$, and cross the imaginary axis again at $\Ea \approx 7.65$. Finally the largest modulus roots at time of crossing cross at  $\Ea \approx 4.55$, turn around at $\Ea \approx 5.85$, and cross back to the stable side of the imaginary axis at $\Ea \approx 7.05$. 
The other parameters are  $\Gamma=0.2$; $e_\sp$=6.23e-2; $q$= 6.2e-1; $d=0.05$; $\kappa= 0.05$; $\nu= 0.05$; $k$ chosen as described in  \S\ref{ssec:parametrization}; $\tau_\sm$= 2.57e-1; $u_\sm$= 7.43e-1; $e_\sm$=9.71e-1; $c_v=1$; $z_\sp=1$; $z_\sm= 0$; $\tau_\sp=1$; $u_\sp=0$; $s=1$; $y_\spm=0$; $T_\sp$=6.2e-2; $T_\sm$= 9.71e-1; and $\Ti$=6.64e-2.}    
\label{fig:025}
\end{figure}

\begin{figure}[ht] 
   \centering
   \begin{tabular}{ccc}
   (a) \includegraphics[width=1.6in]{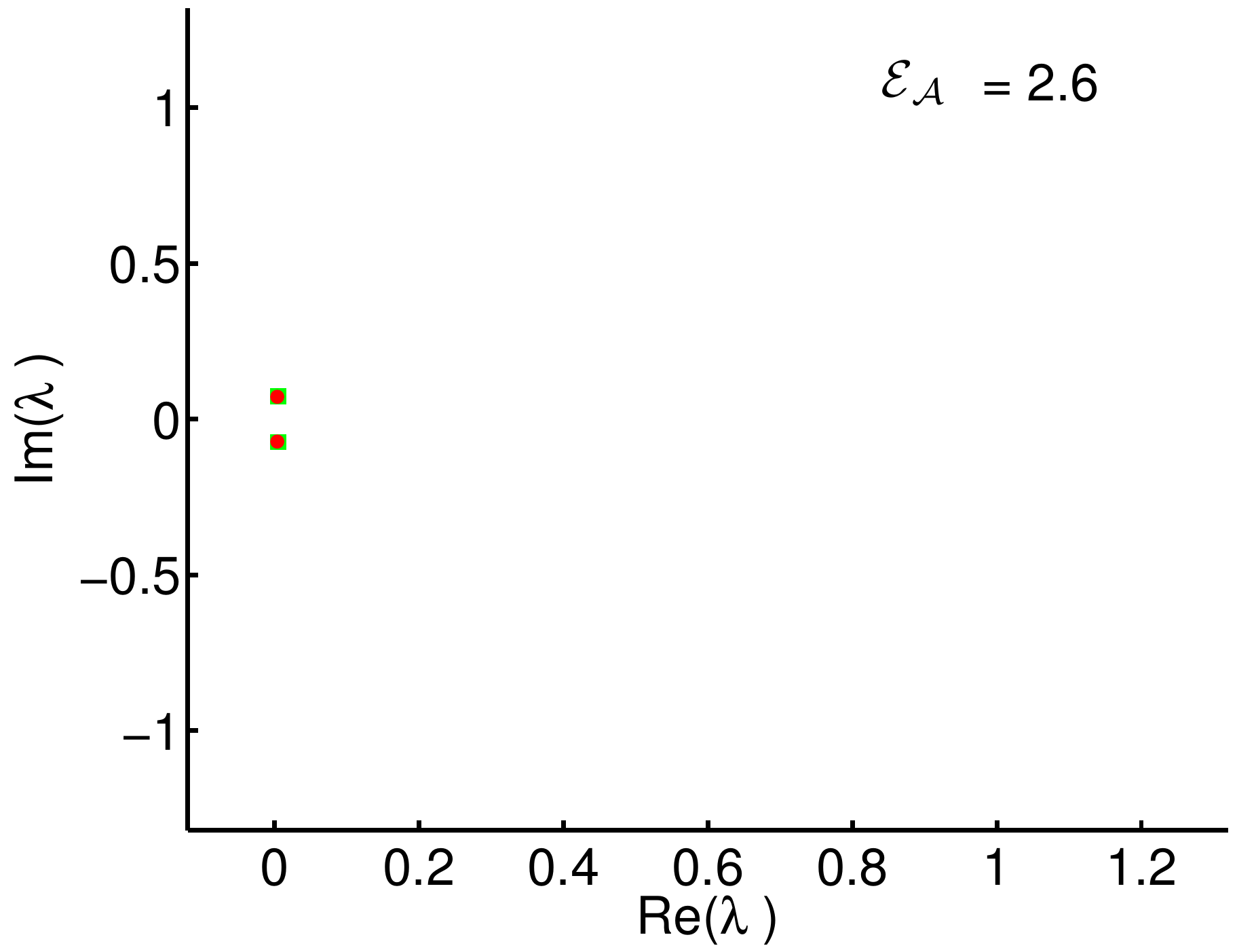} &  
   (b) \includegraphics[width=1.6in]{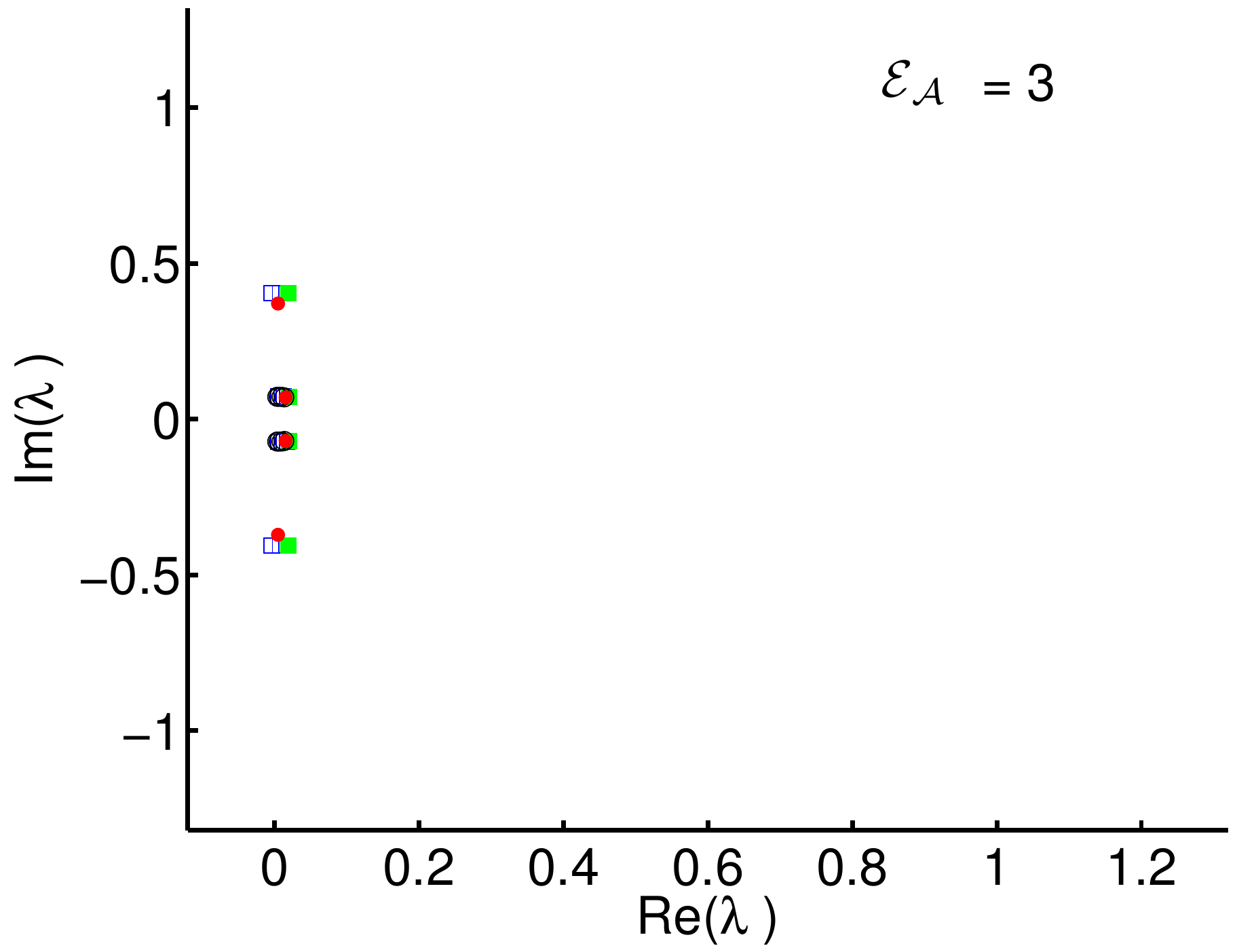} &  
   (c) \includegraphics[width=1.6in]{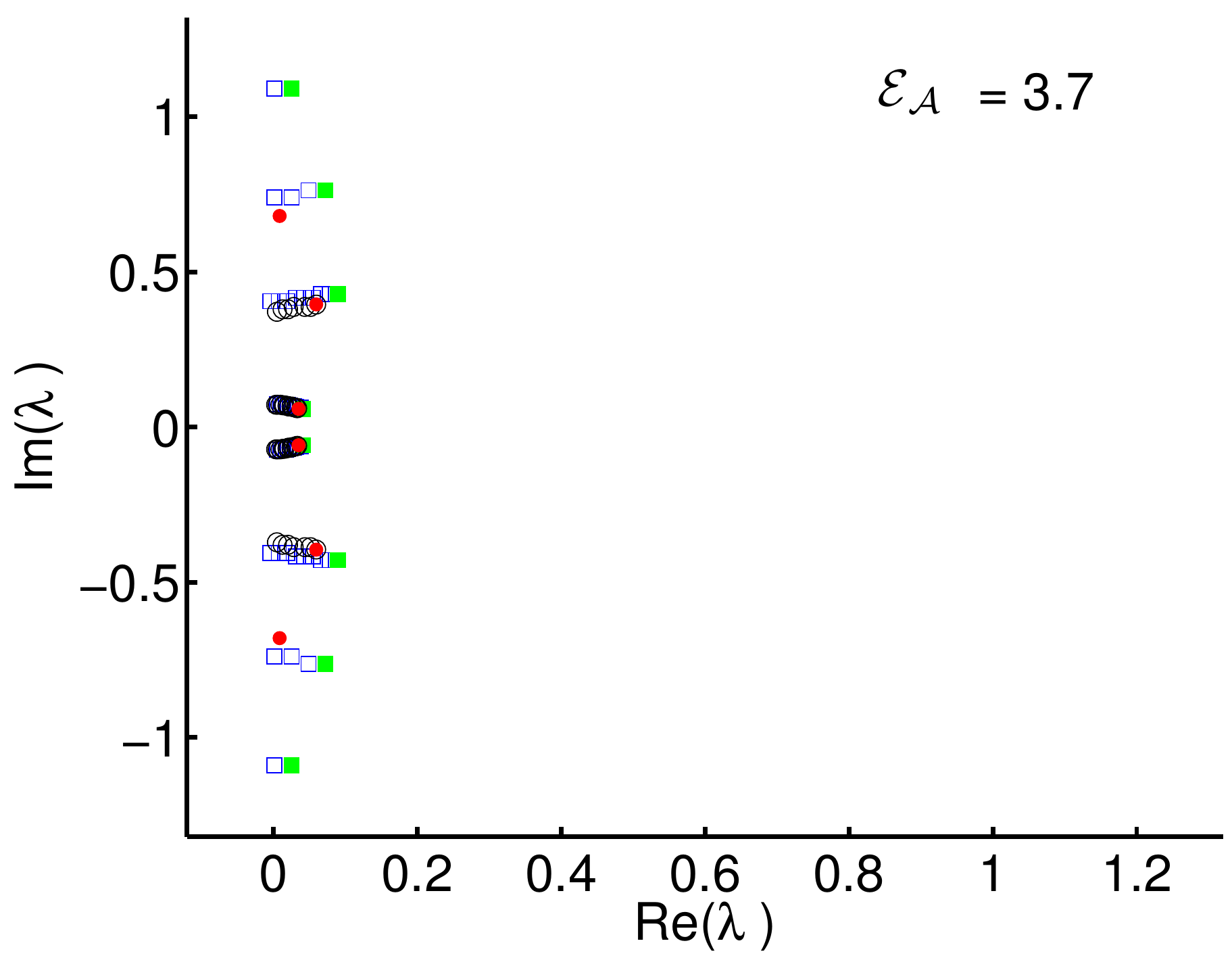} \\
   (d) \includegraphics[width=1.6in]{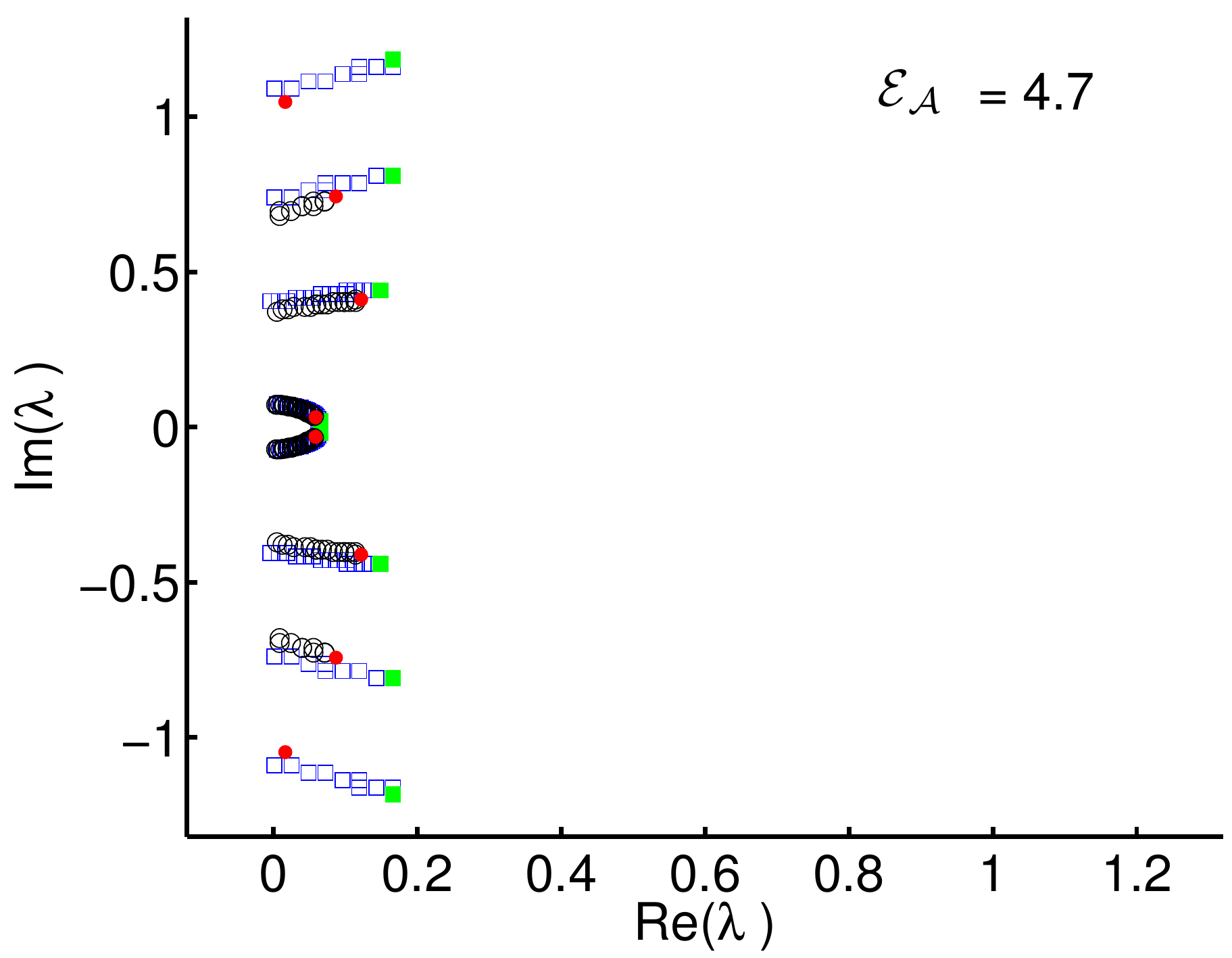} &  
   (e) \includegraphics[width=1.6in]{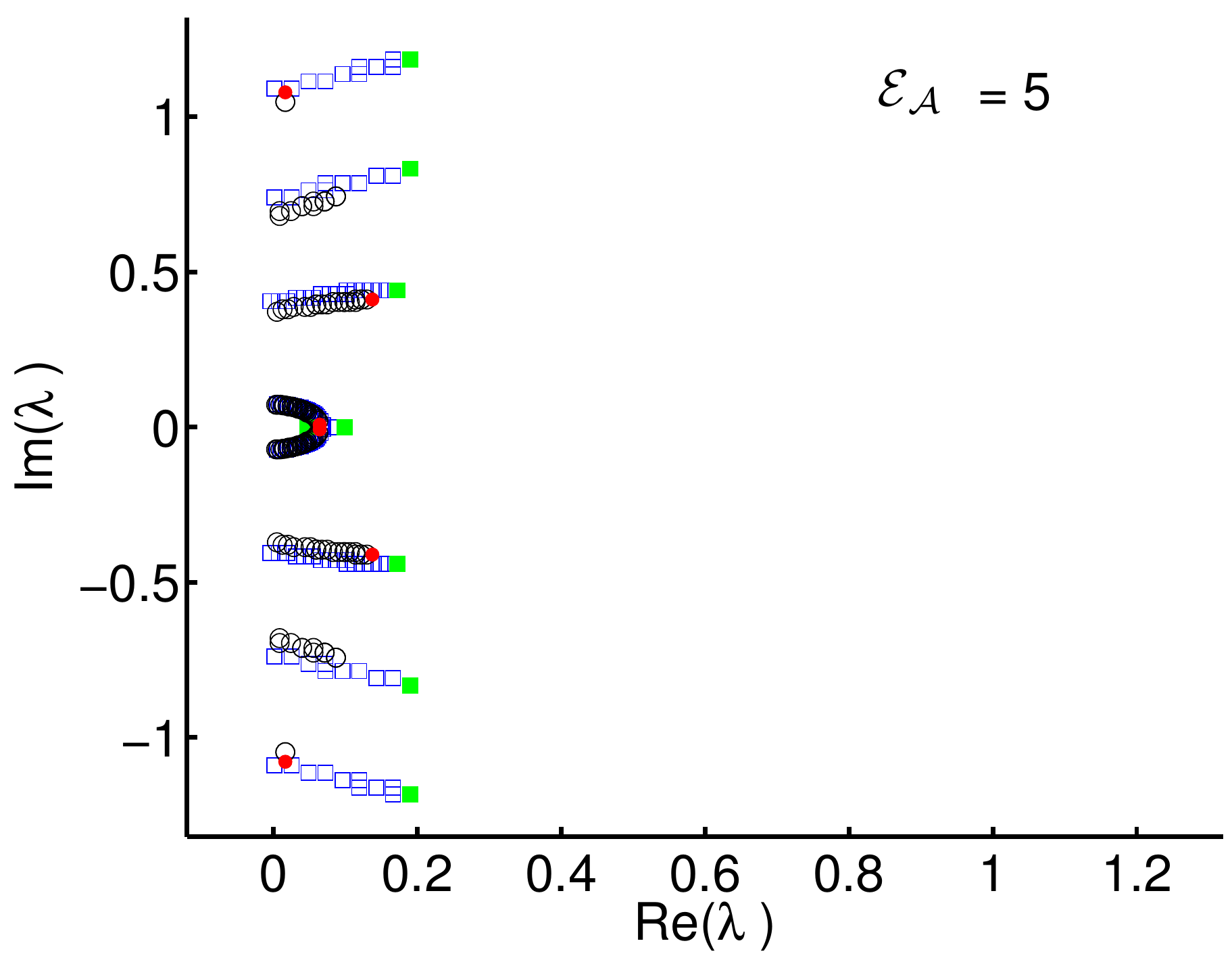} &  
   (f) \includegraphics[width=1.6in]{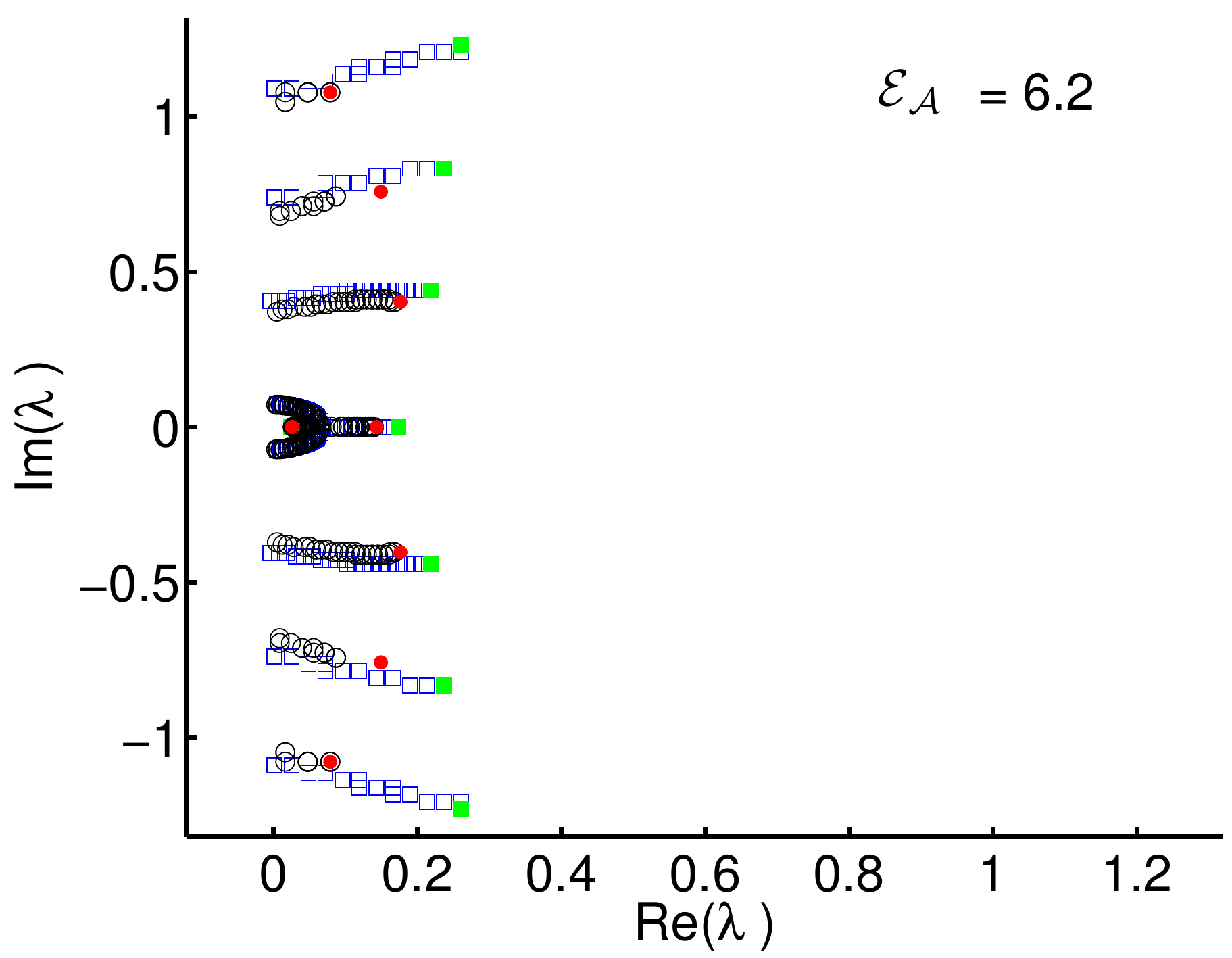} \\
   (g) \includegraphics[width=1.6in]{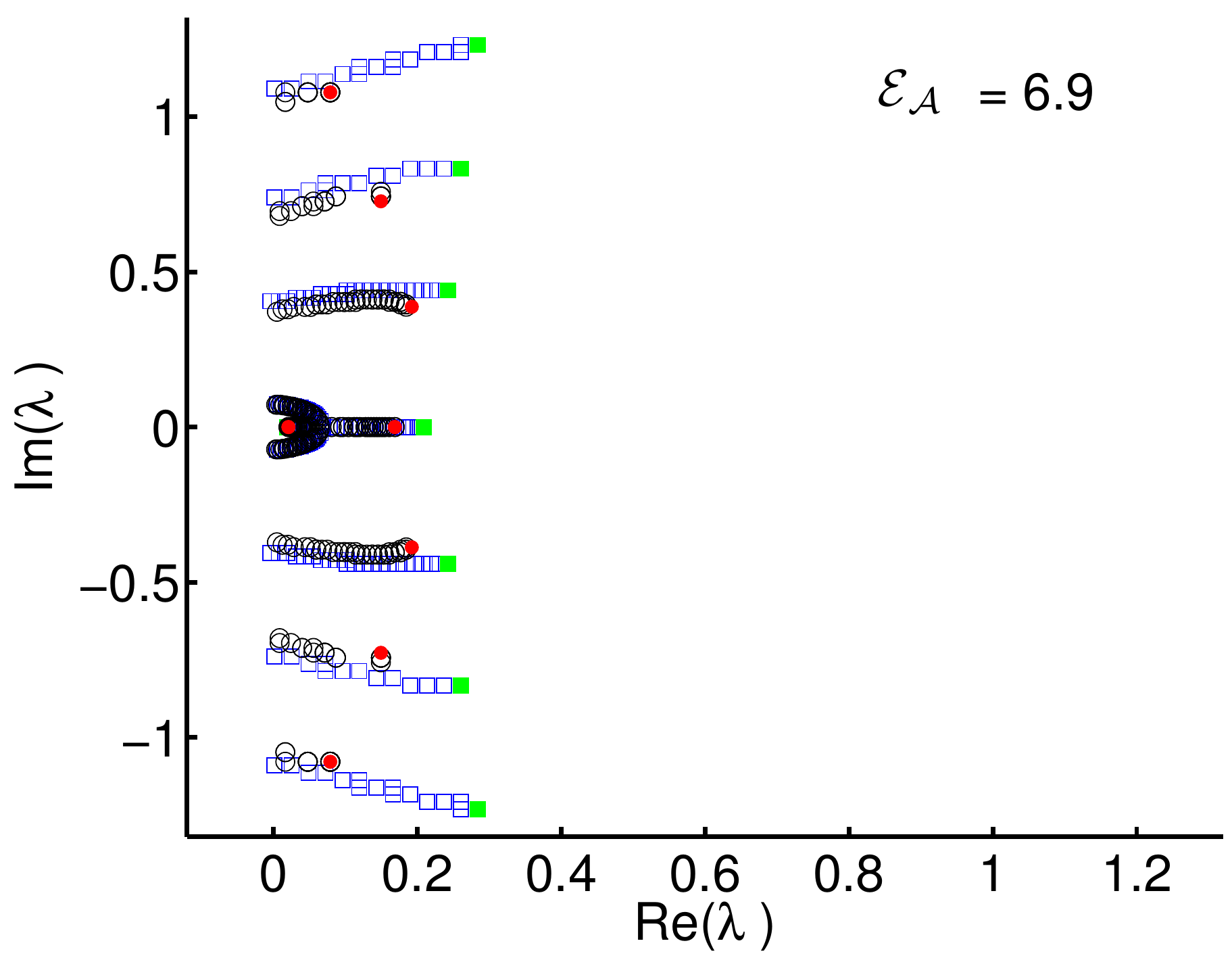} &  
   (h) \includegraphics[width=1.6in]{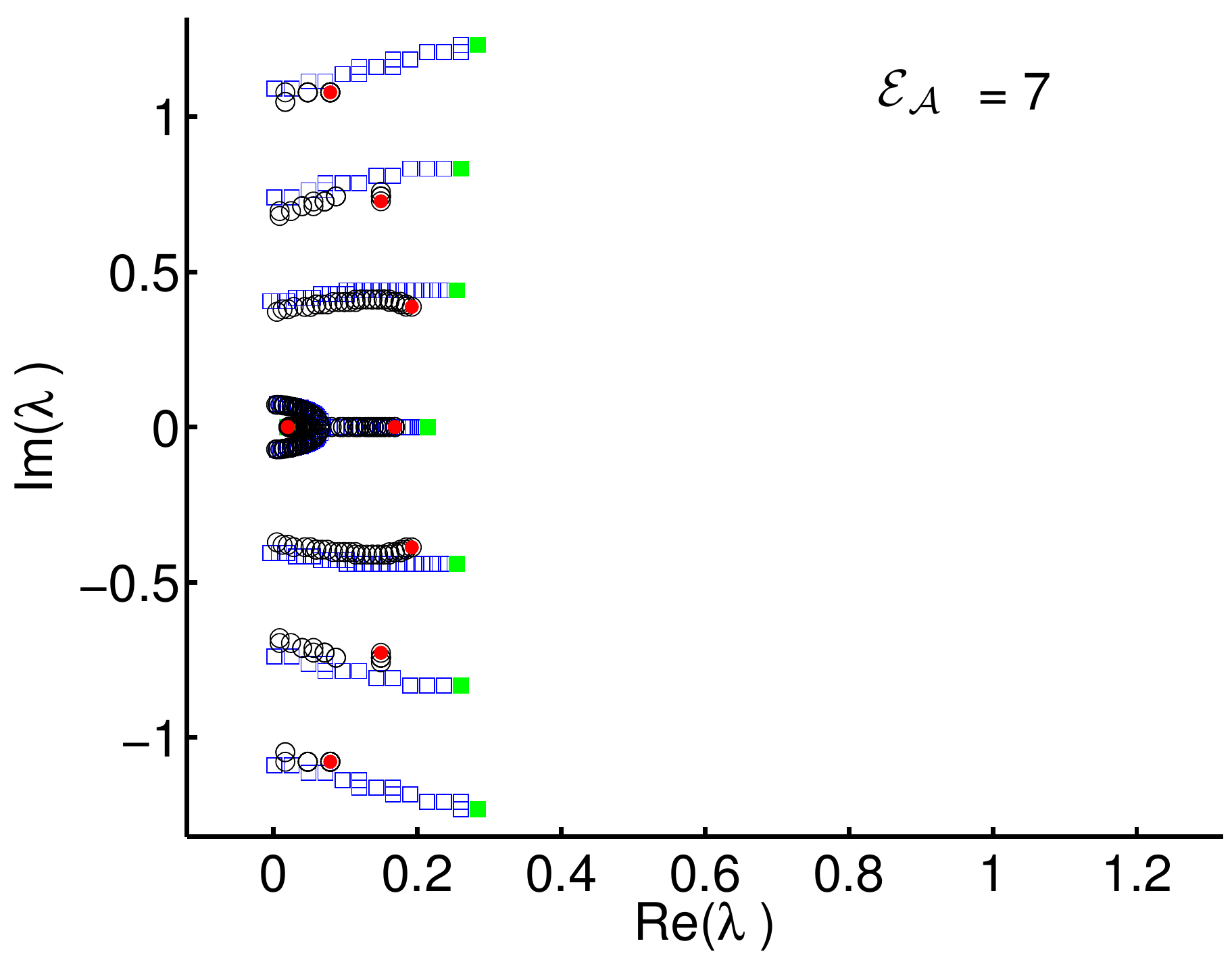} &  
   (i) \includegraphics[width=1.6in]{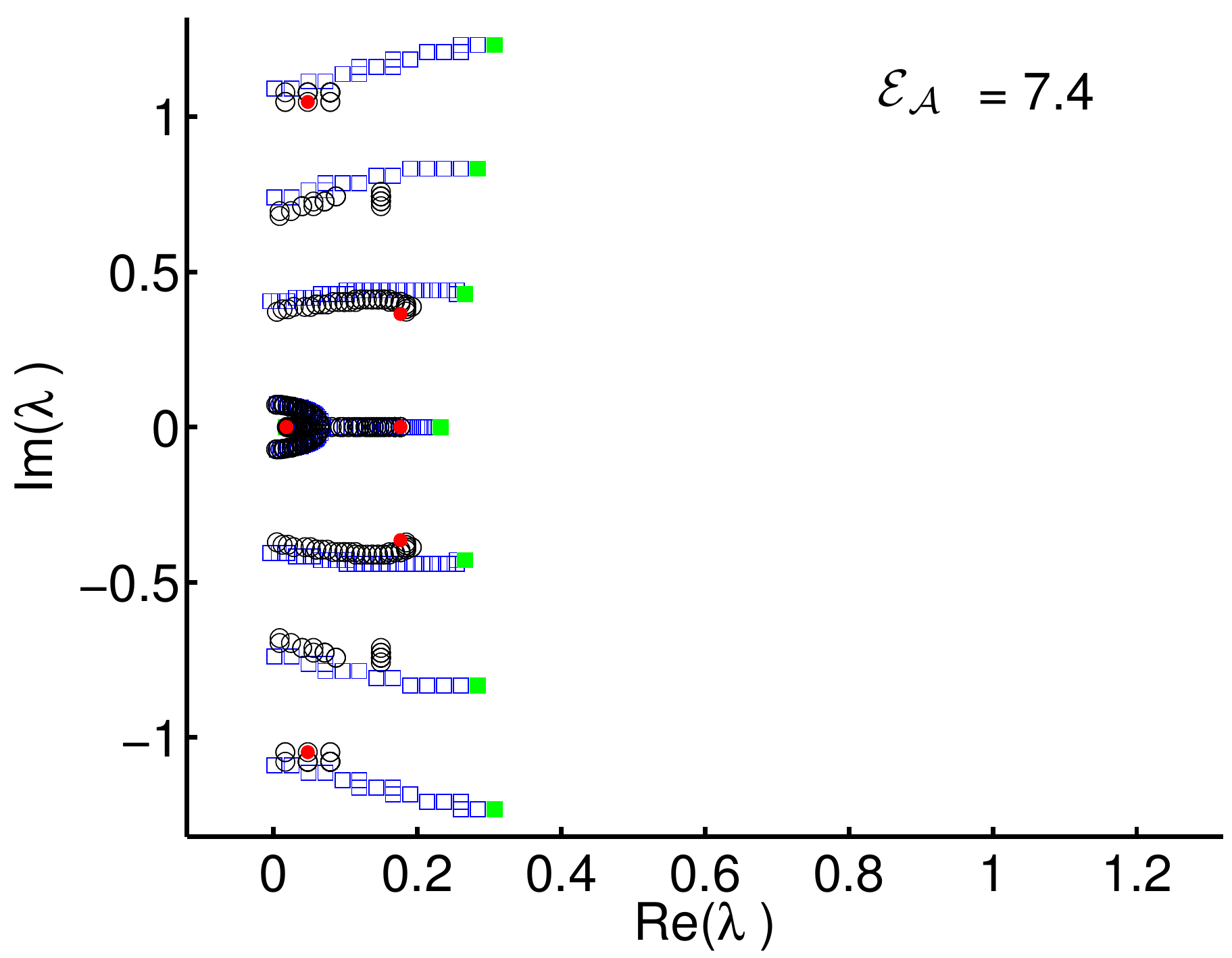} \\
   (j) \includegraphics[width=1.6in]{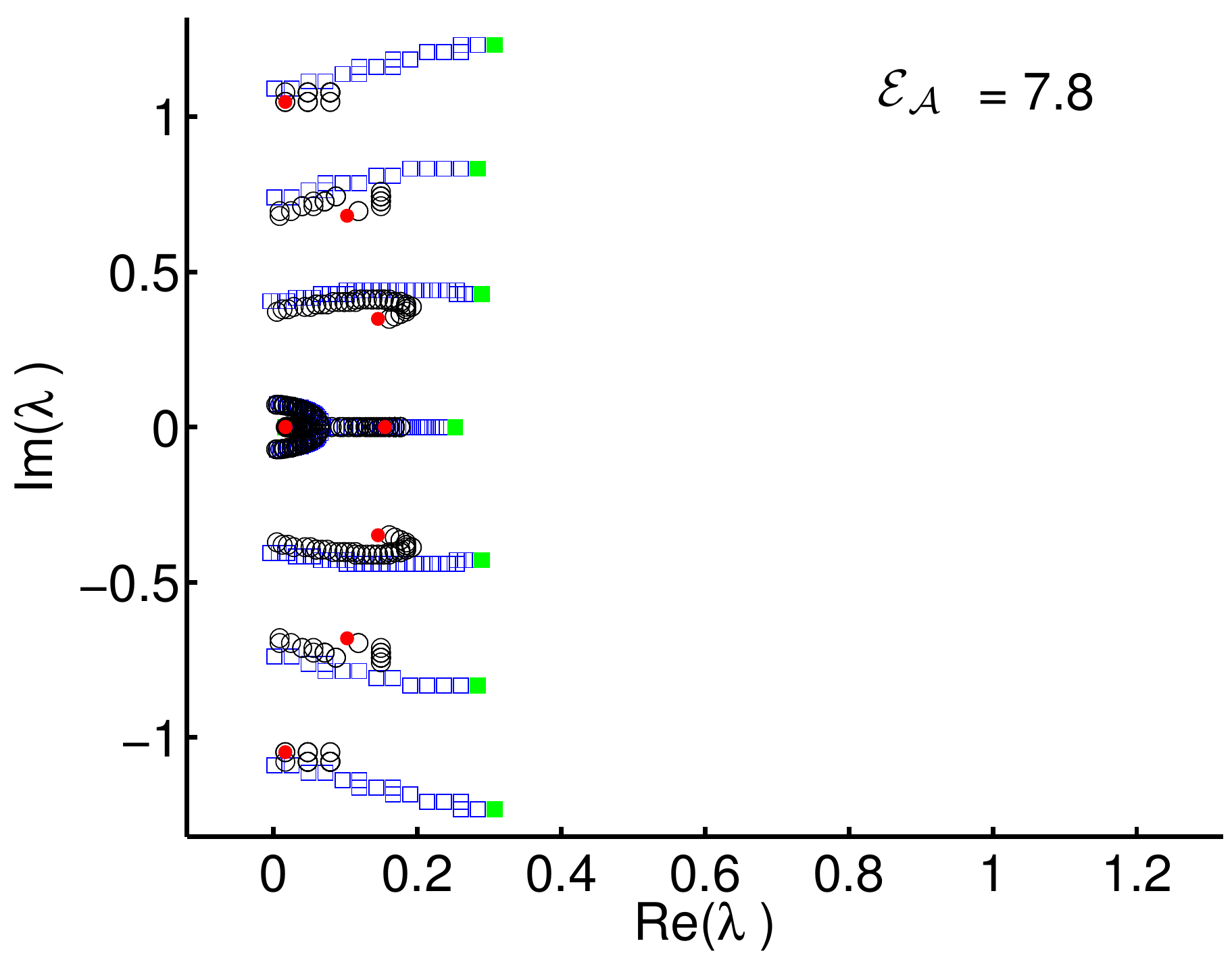} &  
   (k) \includegraphics[width=1.6in]{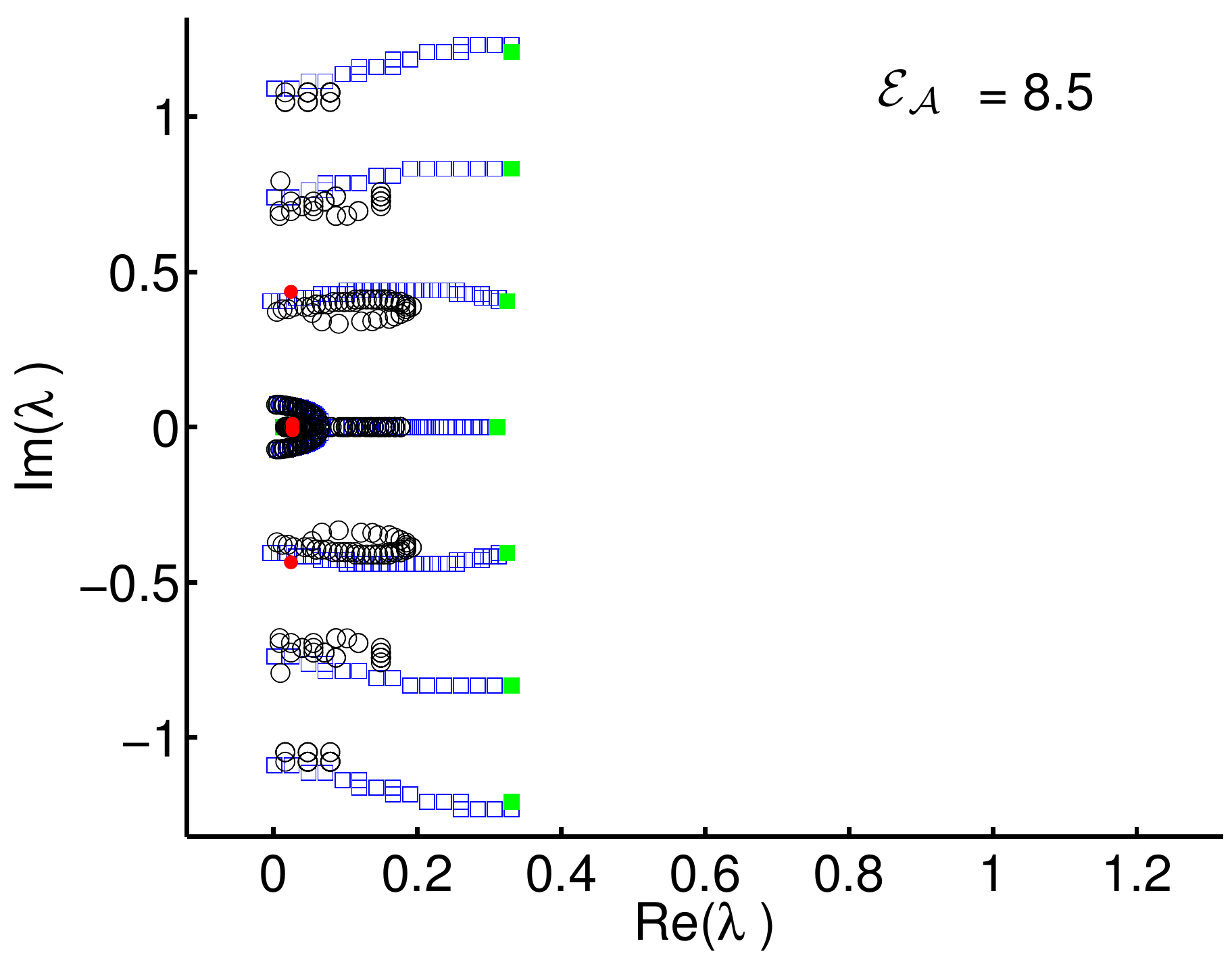} &  
   (l) \includegraphics[width=1.6in]{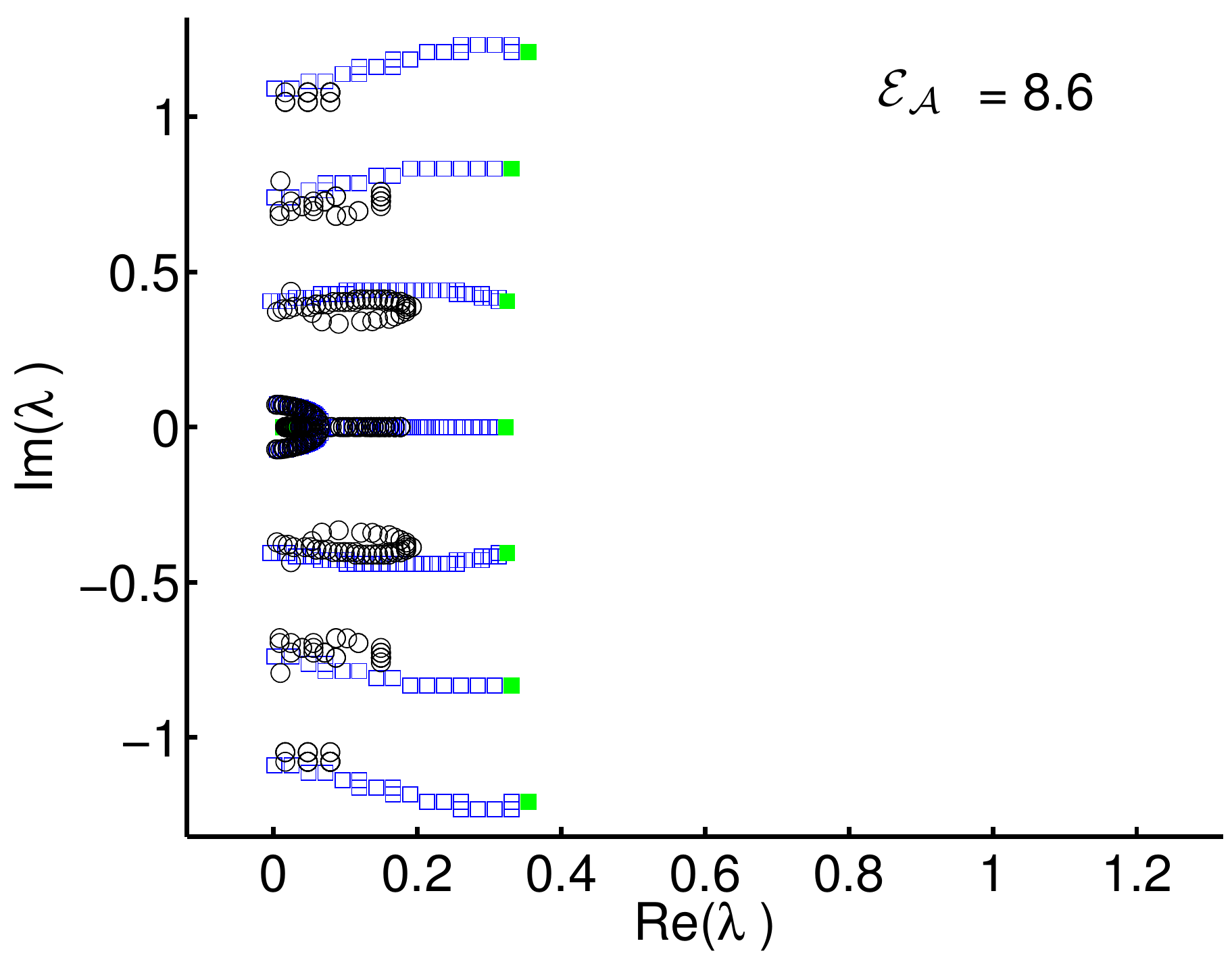} \\
   \end{tabular}
   \caption{The movement of unstable roots in the complex plane as $\Ea$ increases. Open circles (solid dot for current value of $\Ea$) mark the roots corresponding to the RNS model, while open squares (solid square for current value of $\Ea$) correspond to the ZND model. The roots with smallest modulus at the point they cross the imaginary axis appear at $\Ea \approx 2.6$, collide at $\Ea \approx 5$, then at $\Ea \approx 7.4$ the largest real roots turns back, and finally they cross back to the stable side of the complex plane at $\Ea \approx 8.6$. The roots with second largest modulus at initial time of crossing appear for $\Ea \approx 3$, turn at $\Ea \approx 7.0$ and cross back at $\Ea \approx 8.6$. The third largest modulus roots at time of initial crossing appear for $\Ea \approx  3.7$, turn at $\Ea \approx 6.2$, and cross again at $\Ea \approx 8.45$. Finally, the largest modulus roots cross at $\Ea \approx 4.7$, turn at $\Ea \approx 6.9$, and cross back at $\Ea \approx 7.85$.  The other parameters are  $\Gamma=0.2$; $e_\sp$=6.23e-2; $q$= 6.2e-1; $d=0.025$; $\kappa= 0.025$; $\nu= 0.025$; $k$ chosen as described in  \S\ref{ssec:parametrization}; $\tau_\sm$= 2.57e-1; $u_\sm$= 7.43e-1; $e_\sm$=9.71e-1; $c_v=1$; $z_\sp=1$; $z_\sm= 0$; $\tau_\sp=1$; $u_\sp=0$; $s=1$; $y_\spm=0$; $T_\sp$=6.2e-2; $T_\sm$= 9.71e-1; and $\Ti$=6.64e-2.}  
\label{fig:025-2}
\end{figure}

We see similar behavior even for smaller values of the viscous parameters; illustrative results are shown in Figures \ref{fig:025} and \ref{fig:025-2} which feature viscous parameters one-half and one-fourth the size of those described in the previous subsection. Here, we find three and four complex conjugate pairs of eigenvalues crossing into the unstable half plane, but their behavior is the same---a return to stability. The smaller viscosity increases the level of \Ea\ at which the return to stability occurs, but this growth seems to be slow; as noted below in \S\ref{sec:discuss}, we see a further investigation of the growth of the upper stability boundary to be an important direction for further research.

%

\begin{remark}[Dependence on $\Ti$]
Our experiments show that raising \Ti\ seems to make roots turn back sooner (in \Ea). One possible explanation for this behavior may be found in the simple observation that moving $\Ti$ and changing \Ea\ do similar things, roughly speaking.
\end{remark}

\subsection{Lower and upper neutral stability boundaries}\label{ssec:unsb}
\begin{figure}[ht]
\centering 
\begin{tabular}{cc}
(a) \includegraphics[width=0.4\textwidth]{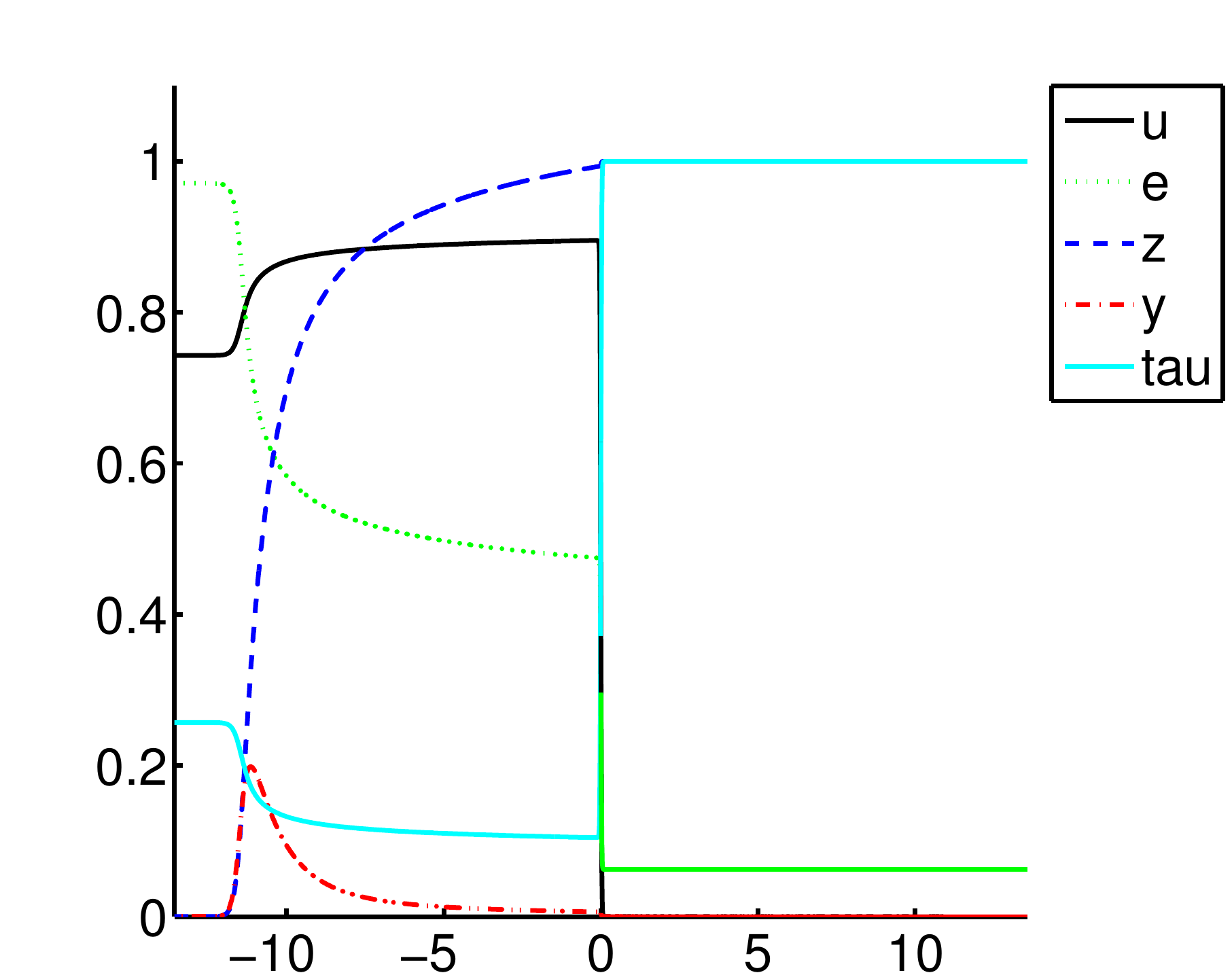} &
(b) \includegraphics[width=0.4\textwidth]{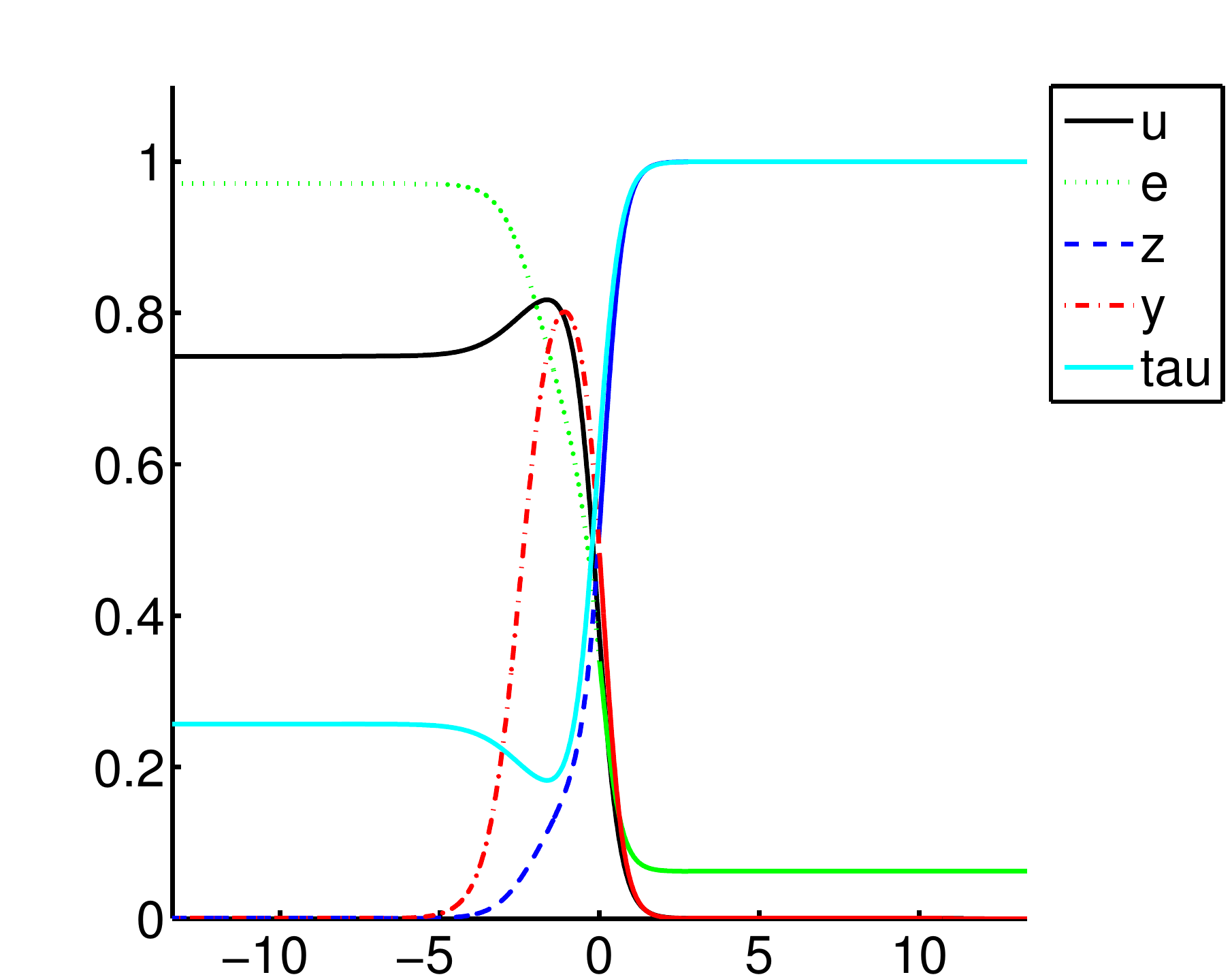}
\end{tabular}
\caption{Sample profiles illustrating diffusive effects. (a) $\nu=d=\kappa=0.01$. (b) $\nu=d=\kappa=0.3$.
In both cases the reaction zone structure is clearly visible, but in (b) the shock width is
of a similar order as the reaction zone width.  
%
For both plots, $\Ti =6.64$e-2, $e_\sp =6.23$e-2, $k = 1.53$e4, $q =6.23$e-1, $\Ea = 6$,
$\Gamma = 0.2$, $c_v= 1$. }
\label{fig:profiles}
\end{figure}

\begin{figure}[ht]
\centering
\includegraphics[scale=0.5]{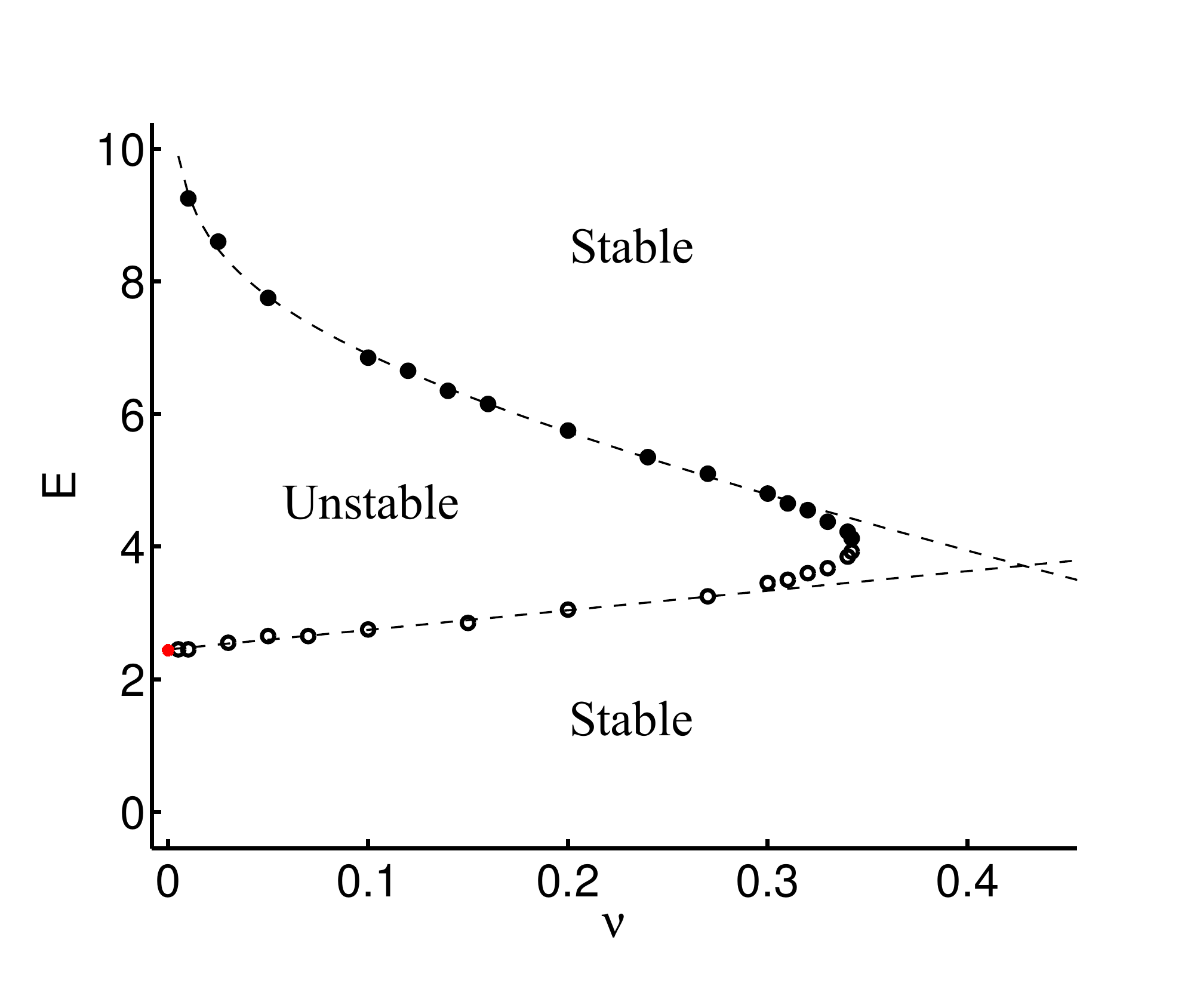}
\caption{Growth of the Upper Stability Boundary. We plot the neutral stability boundaries for RNS; the plot shows the activation energy ($\Ea$) against viscosity $\nu$; the lower stability boundary ($\Ea^\sm$) is marked by black circles and the ZND (inviscid) stability boundary with a red dot. The plot of the upper stability boundary (values of $\Ea^\sp$) is denoted by solid dots. We do curve fitting for the upper boundary with points  such that viscosity is less than 0.27, and less than or equal to 0.27 for the lower. The best-fit curve (dashed line) for the upper boundary is 
$$
\Ea^\sp(\nu)=5.67-6.16\nu -0.804 \ln(\nu)\,.
$$ 
For the lower boundary, the best-fit curve is linear:  $\Ea^\sm(\nu)=2.45+2.95\nu$. 
Here, $\nu  = d = \kappa$, $\Gamma = 0.2$, $e_\sp=6.23$e-2, $q$ = 6.23e-1, and $\Ti= 0.99$. The $(\nu,\Ea)$ values plotted here can be found in Tables \ref{tab:upper} and \ref{tab:lower} below. 
}
\label{fig:sb}
\end{figure}

We found the hyperstabilization phenomenon for the entire range of viscosity parameters in our numerical experiments. That is, for fixed viscosity parameters, we see, as \Ea\ increases, the emergence of complex-conjugate pairs of roots into the unstable half plane signaling instability, and we also find, as $\Ea$ continues to increase, that these unstable roots return to the closed stable half plane. We denote by $\Ea^\sm$ the value of the \emph{lower} neutral stability boundary (in advance of the appearance of unstable eigenvalues) and by $\Ea^\sp$ the value of the activation energy at which this return to neutral stability (no positive real part zeros of the Evans function) is achieved---the \emph{upper} neutral stability boundary. Figure \ref{fig:sb} shows a plot of the dependence of $\Ea^\spm$ on the strength of the viscosity. Notably, the figure shows that when the viscous effects are large enough, there is no instability. Indeed, figure \ref{fig:profiles} illustrates the effect of viscous strength on the shape of the traveling-wave profiles. The right-hand figure (panel (b)) is computed at a viscosity value close to the ``nose'' in Figure \ref{fig:sb} where the lower and upper boundaries coalesce. On the other hand, the left-hand figure (panel (a)) is computed near the inviscid ZND limit.  
As required, $\Ea^\sp$ increases  as viscosity decreases \cite{Z_ARMA11}. Notably, however, the growth appears to be quite slow; this we see as an intriguing avenue for further investigation. In fact, the best-fit curve plotted in Figure \ref{fig:sb} (dashed line) for the upper boundary includes a logarithmic term. By contrast, the behavior of the lower boundary seems to be much more regular. In this case, the best-fit curve (dashed line) for the relationship between $\Ea^\sm$ and $\nu$ is simply linear.  

\begin{table}[ht]
\caption{Upper boundary: $\Ea^\sp(\nu)$. Computing the Evans function on a semi-annulus of outer radius 10 and and inner radius $10^{-4}$, we find the presence of roots for $\Ea$ below the values in the table, and no roots for $\Ea$ above the table values. In most cases, we found the boundary with an absolute error of 0.05. For some values, we found the boundary with absolute error of 0.025; for those values we report more digits. 
}
\begin{tabular}{|c||c|c|c|c|c|c|c|c|c|c|c|c|c|c|c|c|}
\hline
 $\nu$ & 0.01 & 0.025 & 0.05 &  0.1 & 0.12 & 0.14 & 0.16 & 0.2\\
 \hline
 $\Ea^\sp$ & 9.25 &  8.6 & 7.75 & 6.85 &  6.65 & 6.35 & 6.15 &  5.75 \\ 
\hline
\hline
$\nu$ & 0.24 & 0.27 & 0.3 & 0.31 & 0.32 & 0.33 & 0.34 & 0.342 \\
\hline
$\Ea^\sp$ &  5.35 & 5.1 & 4.8 &  4.65 & 4.55 & 4.375 & 4.225 & 4.125\\
\hline
\end{tabular}
\label{tab:upper}
\end{table}

\begin{table}[ht]
\caption{Lower boundary: $\Ea^\sm(\nu)$. Computing the Evans function on a semi-annulus of outer radius 10 and and inner radius $10^{-4}$, we find the presence of roots for $\Ea$ above the values in the table, and no roots for $\Ea$ below the table values. In most cases, we found the boundary with an absolute error of 0.05. For some values, we found the boundary with absolute error of 0.025; for those values we report more digits. 
}
\begin{tabular}{|c||c|c|c|c|c|c|c|c|c|c|c|c|c|c|c|c|}
\hline
 $\nu$ & 0.005 & 0.01 & 0.03 &  0.05 & 0.07 & 0.1 & 0.15 & 0.2\\
 \hline
 $\Ea^\sm$ & 2.45&  2.45 & 2.55 & 2.65 &  2.65 & 2.75 & 2.85 &  3.05 \\ 
\hline
\hline
$\nu$ & 0.27 & 0.3 & 0.31 & 0.32 & 0.33 & 0.34 & 0.342& \\
\hline
$\Ea^\sm$ &  3.25 & 3.45 & 3.5 &  3.6 & 3.675 & 3.85 & 3.925& \\
\hline
\end{tabular}
\label{tab:lower}
\end{table}

\section{Discussion}\label{sec:discuss}
\subsection{Findings}\label{ssec:findings}
\begin{figure}[htb]
\centering
\includegraphics[width=0.4\textwidth]{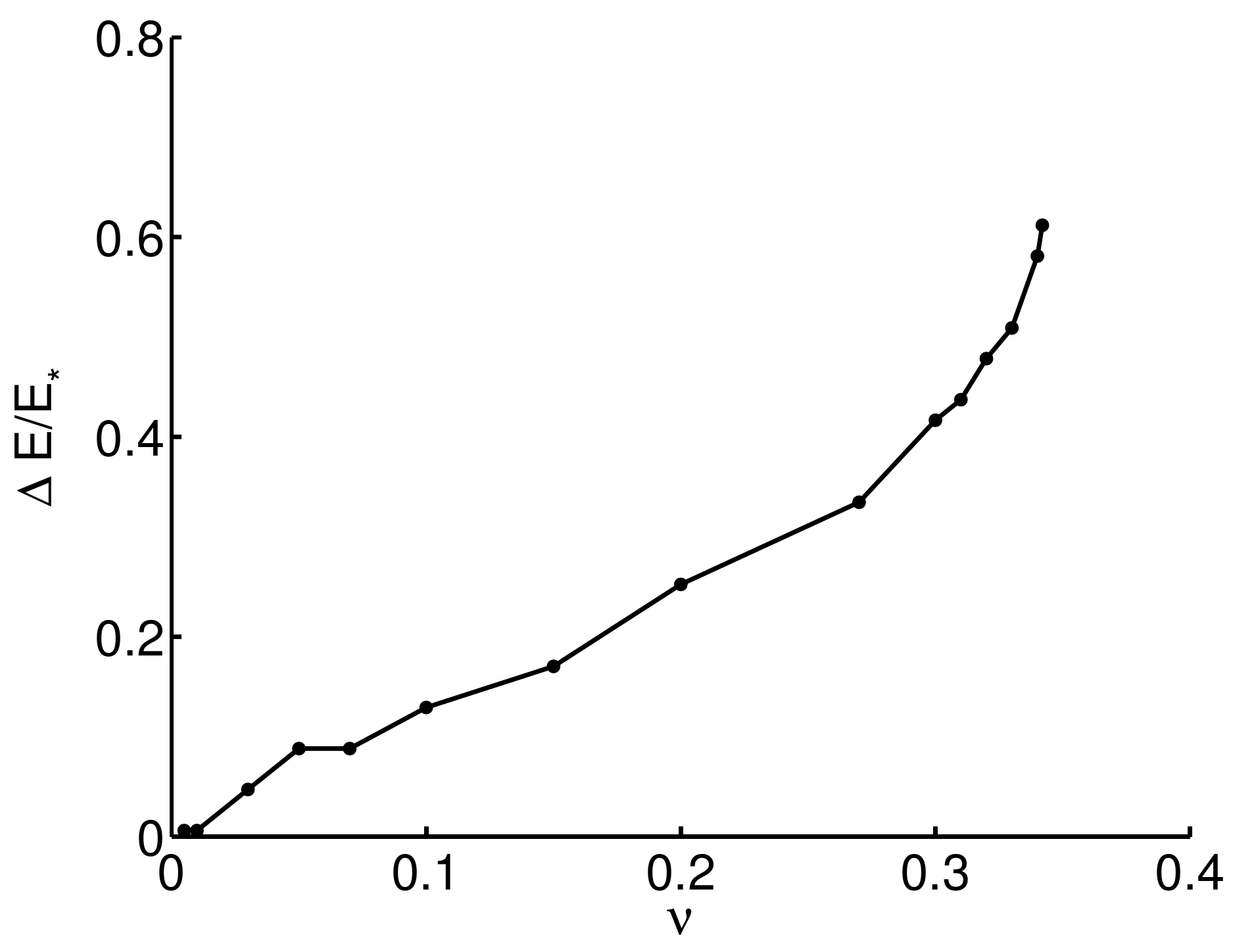}
\caption{Viscous Delay: We plot $\Delta E/\mathcal{E}_*=(\Ea^\sm(\nu)-\mathcal{E}_*)/\mathcal{E}_*$ against $\nu$, where $\mathcal{E_*}$ is the approximation to the ZND neutral boundary. 
Here, $\nu  = d = \kappa$, $\Gamma = 0.2$, $e_\sp=6.23$e-2, $q$ = 6.23e-1, and $\Ti= 0.99$. 
}
\label{fig:delay}
\end{figure}

Our experiments indicate that viscous effects \emph{can} be important in the stability analysis of detonation waves. First, we note that the inclusion of these effects delays the onset of instability; see Figure \ref{fig:delay}. 
We thus recover the results of Romick et al.\ \cite{RAP_JFM12}, and we confirm that the effect of viscosity on $\Ea^\sm$ is essentially perturbative with a linear relationship between the value of $\Ea^\sm$ and the strength of viscosity.  
On the other hand, for all of the viscosity values tested, the difference between RNS and ZND \emph{at the spectral level} is striking; see, e.g., Figure \ref{fig:025}. Moreover, our results reveal the unexpected phenomenon of \textbf{viscous hyperstabilization}---a return to stability---as the bifurcation parameter (activation energy) increases.
Now, as described above in \S\ref{ssec:vs}, the broad range of viscosities considered in our experiments includes a physical range, starting with the fairly large viscosity considered by Romick et al.\ \cite{RAP_JFM12} and decreasing. Thus, it is of significant interest to determine whether or not the activation energies needed for complete stabilization also lie within a physically relevant range. For the viscosity level considered by Romick et al.\ \cite{RAP_JFM12}, it is roughly twice the activation energy at which instability occurs for Fickett--Wood benchmark problem (7 vs. 3); for smaller viscosities, the level is successively larger, though the growth is surprisingly slow; see Figure \ref{fig:sb} and the surrounding discussion.
It is not clear whether this size of activation energy might arise in physical situations. Nonetheless, even for much smaller activation energies, our results show that the spectra of viscous and inviscid waves differ profoundly. Now, Powers \& Paolucci \cite{PP_AIAAJ05} have shown that the finest length scales in the reaction zones of realistic combustion models are of the same order as diffusive length scales. They conclude that a physically consistent, fully resolved model of detonation must include diffusive effects. In terms of the stability of such waves,  Bdzil \& Stewart have identified ``a better understanding of the hydrodynamic stability characteristics of our engineering-scale models''  as one of four major challenges that must be met before a practically useful predictive theory for detonation phenomena can be constructed \cite{BS_ARFM07}. Combining these observations, we find a call to action to develop a stability theory for detonation waves that incorporates diffusive effects. We see the present work, an Evans-based analysis, as a critical first step toward addressing this challenge. 

\subsection{Future directions}\label{ssec:future}
Several intriguing future directions for investigation have been identified in the body of this paper. Most notably, it would be of substantial mathematical and practical interest to carry out a singular perturbation analysis of the ZND limit and the phenomenon of return to stability seen here. It is known that the roots of the viscous Evans function converge to the roots of the limit ZND (inviscid) Evans--Lopatinski\u\i\ function \cite{Z_ARMA11}, but a better understanding of viscous hyperstabilization and the corresponding growth of the upper stability boundary, as in Figure \ref{fig:sb}, is clearly needed. 
Another possibility would be to extend these calculations to a more sophisticated model including, e.g., multiple space dimensions, multi-step chemistry, and more realistic equations of state; this is clearly necessary to answer the call of Bdzil \& Stewart \cite{BS_ARFM07}. 
We note that for the Navier--Stokes equations of (nonreacting) gas dynamics, the extension of the \textsc{StabLab} package to the multi-dimensional setting is well underway \cites{BHLZ-el,BHLZ-flux}. The incorporation of more realistic multi-step chemical reaction processes provides an excellent opportunity to test this developing computational suite on problems of substantial physical import.  
\appendix

\section{Linearization and Evans function: details}\label{sec:details}

In this section we provide the details of the linearization that is the basis for our computation of the Evans function. This calculation is an essential first step in the stability analysis, and there are a number of choices one must make to obtain a first-order eigenvalue ODE. Given these choices, a natural objective in this regard is to determine which choices lead to the simplest, most useful formulation of the eigenvalue equation for our computational needs. Here, we follow the strategy of using \emph{flux variables} (see \eqref{eq:Ydef} below). Notably, the use of flux variables gives the benefits of the oft-used integrated coordinates, and, moreover, has a natural extension to multidimensional problems (whereas the integrated coordinates do not). See \cite{BHLZ-flux} for an explanation and further discussion of this point. 

\subsection{Preliminaries}
We write the system \eqref{eq:rns} 
as
\beq\label{eq:rns_nc}
f^0(U)_t+\Big(f^1(U)-sf^0(U)\Big)_x=\big(B(U)U_x\big)_x+R(U)\,.
\eeq
where $U=(\t,u,e,z)^\tr$. The functions $f^0$ and $f^1$ are given by 
\beq\label{eq:fs}
f^0(U)=\bp \t \\ u \\ e+\frac{1}{2}u^2 \\ z \ep\,,\quad f^1(U)=\bp -u \\ p(\t,e,z) \\ up(\t,e,z) \\ 0 \ep\,.
\eeq
Also, 
\beq\label{eq:B}
B(U)= 
\left(\begin{array}{c|ccc}
0 & 0 & 0 & 0 \\ \hline 
0 & \frac{\nu}{\t} & 0 & 0 \\
0 & \frac{\nu u}{\t} & \frac{\kappa_v}{\t} & 0 \\
0 & 0 & 0 & \frac{d}{\t^2}
\end{array}\right)=\bp 0 & 0 \\ 0 & b(U) \ep\,,
\eeq
and
\beq\label{eq:R}
R(U)=\bp 0 \\ 0 \\ qk\phi(T(\t,e,z))z \\ -k\phi(T(\t,e,z))z \ep\,.
\eeq
\br
At this point we assume that the temperature $T$ is a function of the internal energy, $e$, alone and does not depend on $\t$ or $z$. This is certainly true in the fundamental case,
\(
T=c_v^{-1}e,
\)
that we consider in the calculations of this paper. 
We may thus write, as above,
\(
\check\phi (e)=\phi(T(e))\,.
\)
This simplifies some of the calculations. See, e.g., the form of \eqref{eq:E} below.
\er
Straightforward calculations show that 
\beq\label{eq:binverse}
b(U)^{-1}=\bp
\frac{\t}{\nu} & 0 & 0 \\
-\frac{\t u}{\kappa_v} & \frac{\t}{\kappa_v} & 0   \\
0 & 0 & \frac{\t^2}{d} 
\ep\,,
\eeq
and that 
\beq\label{eq:E}
E(U):=\dif R(U)=
\left(\begin{array}{c|ccc}
0 & 0 & 0 & 0 \\ \hline
0 & 0 & 0 & 0 \\
0 & 0 & qk\check\phi'(e)z & qk\check\phi(e) \\
0 & 0 & -k\check\phi'(e)z & -k\check\phi(e)
\end{array}\right)\,.
\eeq
We also record here, for future reference, the forms of the Jacobian matrices $a^j=\dif f^j$ for $j=0,1$. They are
\beq\label{eq:a0a1}
a^0(U):=
\left(\begin{array}{c|ccc} 
1 & 0 & 0 & 0 \\ \hline 
0 & 1 & 0 & 0 \\ 
0 & u & 1 & 0 \\
0 & 0 & 0 & 1
\end{array}\right)
\quad\text{and}\quad
a^1(U):=
\left(\begin{array}{c|ccc} 
0 & -1 & 0 & 0 \\ \hline 
p_\t & 0 & p_e & p_z \\ 
up_\t & p & up_e & up_z \\
0 & 0 & 0 & 0
\end{array}\right)\,.
\eeq

\br
We will utilize the block structure of the matrix $B$ in our calculations below. 
\er

\subsection{Linearization and block structure}
The linearized system of interest is given by 
\beq\label{eq:lin}
\lambda a^0(x) W+\big(A(x)W-B(x)W'\big)'=E(x)W\,.
\eeq
The unknown is 
\beq
W=\bp \t \\ \hline  u \\ e \\ z \ep\,,
\eeq
with $\t$, $u$, $e$, $z$ now representing (Laplace-transformed) perturbations. In \eqref{eq:lin}, by a slight abuse of notation, we have written $a^0(x)=a^0(\bar U(x))$, $B(x)=B(\bar U(x))$, etc. More,
\[
A(x)W:=[a^1(\bar U(x))-sa^0(\bar U(x))]W-\dif B(\bar U)(W,\bar U_x)\,.
\]
We will almost always suppress the dependence of the coefficients on $x$. 
Our aim is to write the eigenvalue problem \eqref{eq:lin} in the form of a first-order system 
\beq\label{eq:1order}
\calW'=G(x;\l)\calW\,.
\eeq
To that end, we write 
\beq
W=\bp W_1 \\ \hline W_2 \ep=\bp \t \\ \hline u \\ e \\ z \ep\,,
\eeq
with
\beq
W_1=\t\,\quad\text{and}\quad W_2=\bp u \\ e \\ z \ep\,,
\eeq
and we write the matrices $A$, $B$, and $E$ in corresponding block form. That is, we write
\beq\label{eq:ablock}
A=\bp A_{11} & A_{12} \\ A_{21} & A_{22} \ep\,,
\eeq
so that $A_{11}$ is $1\times 1$, $A_{12}$ is $1\times 3$, $A_{21}$ is $3\times 1$, and $A_{22}$ is $3\times 3$, and we use the same block decomposition and notational conventions for the matrices $B$ and $E$. 
%

\subsection{Flux variables}\label{ssec:flux}

We define $Y=(Y_1,Y_2)^\tr$ by 
\beq\label{eq:Ydef}
Y:=AW-BW'\,.
\eeq
When we rewrite \eqref{eq:Ydef} in terms of the blocks, we see
\beq\label{eq:blockY}
\bp Y_1 \\ Y_2 \ep = \bp A_{11} & A_{12} \\ A_{21} & A_{22} \ep\bp W_1 \\ W_2 \ep - \bp 0 & 0 \\ 0 & b \ep \bp W_1' \\ W_2' \ep\,.
\eeq
Then, the vector \calW\ in \eqref{eq:1order} will be defined by 
\beq\label{eq:calw}
\calW:=\bp Y \\ W_2 \ep\,.
\eeq
To extract the equations for $Y_1'$, $Y_2'$, and $W_2'$, we note that 
the system \eqref{eq:lin} takes the block form
\beq\label{eq:evalblock}
\l\bp a_{11}^0 & a^0_{12} \\ a^0_{21} & a^0_{22}\ep\bp W_1 \\ W_2\ep + \bp Y_1' \\ Y_2' \ep=\bp E_{11} & E_{12} \\ E_{21} & E_{22}\ep \bp W_1 \\ W_2 \ep\,.
\eeq

\subsubsection{Equation for $Y_1'$}

Evidently, from the first row of \eqref{eq:blockY}, 
\[
Y_1=A_{11}W_1+A_{12}W_2\,,
\]
whence 
\beq\label{eq:w1}
W_1=A_{11}^{-1}Y_1-A_{11}^{-1}A_{12}W_2\,.
\eeq
\br
We'll see shortly by direct calculation that the $1\times 1$ matrix $A_{11}$ is invertible. We shall use \eqref{eq:w1} to eliminate $W_1$ in favor of $Y_1$ and $W_2$ in the calculations that follow. 
\er

From \eqref{eq:evalblock} and the fact that $E_{11}=0$ and $E_{12}=0$, it is clear that 
\[
\l (a^0_{11}W_1 +a^0_{12}W_2)+Y_1'=0\,.
\]
More, inspection of \eqref{eq:a0a1} reveals that 
\[
a^0_{11}=1\,,\quad a^0_{12} =0\,.
\]
Therefore, upon using \eqref{eq:w1}, we find 
 \beq\label{eq:y1p}
 Y_1'=\big[-\l A_{11}^{-1}\big]Y_1+\big[0\big] Y_2+ \big[\l A_{11}^{-1}A_{12}\big]W_2\,.
 \eeq
 \subsubsection{Equation for $Y_2'$}
 Similarly, we obtain an equation for $Y_2'$. First, we note that the second row of \eqref{eq:evalblock} reads
 \[
\l (a^0_{21}W_1 +a^0_{22}W_2)+Y_2'=E_{21}W_1 +E_{22}W_2\,.
\]
Again, several of the terms are zero. From \eqref{eq:a0a1} and \eqref{eq:E}, we see
\beq
a^0_{21}=0\,,\quad E_{21}=0\,.
\eeq
Thus, we finally obtain 
 \beq\label{eq:y2p}
 Y_2'=\big[0\big]Y_1+\big[0\big]Y_2+\big[E_{22}-\l a^0_{22}\big]W_2\,.
 \eeq
 
 \subsubsection{Equation for $W_2'$}
 All that remains is to derive an ODE for $W_2$. However, from \eqref{eq:blockY} it is clear that 
 \begin{align} \label{eq:bw2p}
 bW_2'&=A_{21}W_1+A_{22}W_2-Y_2 \\ \nonumber
 	&=A_{21}(A_{11}^{-1}Y_1-A_{11}^{-1}A_{12}W_2)+A_{22}W_2-Y_2\\ \nonumber
	&=A_{21}A_{11}^{-1}Y_1-Y_2+\big(A_{21}(-A_{11}^{-1}A_{12})+A_{22}\big)W_2\,.
 \end{align}
 Thus, from \eqref{eq:bw2p}, we see that 
 \beq\label{eq:w2p}
 W_2'=\big[b^{-1}A_{21}A_{11}^{-1}\big] Y_1 +\big[-b^{-1}\big]Y_2+\big[b^{-1}(-A_{21}A_{11}^{-1}A_{12}+A_{22})\big]W_2\,.
 \eeq
 
 \subsection{System}\label{ssec:system}
 We write the entries of the coefficient matrix $G$ as 
 \beq
 G=\bp 
 g_{11} & g_{12} & g_{13} \\
  g_{21} & g_{22} & g_{23}\\
   g_{31} & g_{32} & g_{33}
   \ep\,,
 \eeq
 and we note that the form of the entries of $G$ can simply be read off from \eqref{eq:y1p}, \eqref{eq:y2p}, and \eqref{eq:w2p}. That is,
\begin{align*}
 g_{11} & = -\l A_{11}^{-1}\,, &
 g_{12} & = 0\,,&
 g_{13} & = \l A_{11}^{-1}A_{12}\,, \\
 g_{21} & = 0\,,&
 g_{22} & = 0 \,, &
 g_{23} & = E_{22}-\l a^0_{22}\,, \\
 g_{31} & = b^{-1}A_{21}A_{11}^{-1}\,,&
 g_{32} & = -b^{-1}\,,&
 g_{33} & = b^{-1}(-A_{21}A_{11}^{-1}A_{12}+A_{22})\,.
 \end{align*}
We now use the explicit forms of $a^0$, $a^1$, $b$, and $E$ above to expand these entries. First, we note that, since the first row of $B$ is zero, 
 \beq
 A_{11}=a^1_{11}-sa^0_{11}=-s
 \eeq
 so that 
 \beq
 A_{11}^{-1}=-\frac{1}{s}\,.
 \eeq
  \br
 Recall that we require $s\neq 0$. 
 \er
 Also, 
 \beq
 A_{12}=a^1_{12}-sa^0_{12}= \bp -1 & 0 & 0 \ep\,.
 \eeq
 Thus, the entries in the first row of $G$ are 
 \begin{align}
 g_{11} & =\l/s\,,\\
 g_{12} & =\bp 0 & 0 & 0 \ep\,, \\
 g_{13} & =\bp \l/s & 0 & 0 \ep\,.
 \end{align}
 
 Similarly, we see that the entries in the second row of $G$ are 
 \begin{align}
 g_{21} & = \bp 0 \\ 0 \\ 0 \ep\,, \\
 g_{22} & = \bp 0 & 0 & 0 \\ 0 & 0 & 0 \\ 0 & 0 & 0 \ep\,, \\
 g_{23} & = \bp -\l & 0 & 0 \\ 
 		-\l u & -\l + qk\check\phi'(e)z & qk\check\phi(e) \\
		0 & -k\check\phi'(e)z & -\l-k\check\phi(e)
		\ep\,.
 \end{align}
 Finally, it remains to compute the entries in the third row. We start by calculating the two lower blocks of $A$. They are
 \beq\label{eq:A21}
 A_{21} =
 \bp p_\t + \displaystyle{\frac{\nu u_x}{\t^2}} \\
 p_\t u+\displaystyle{\frac{\nu uu_x}{\t^2}+\frac{\kappa_v e_x}{\t^2}} \\
 2\displaystyle{\frac{d z_x}{\t^3}}
 \ep\,,
 \eeq
 and
 \beq\label{eq:A22}
 A_{22} = 
 \bp
 -s & p_e & p_z \\ 
 -su + p -\displaystyle{\frac{\nu u_x}{\t}} & -s+up_e & up_z \\
 0 & 0 & -s
 \ep\,.
 \eeq
 Then, 
 \begin{align}\label{eq:g31}
 g_{31} & = b^{-1}A_{21}A_{11}^{-1} \\
 	& = \bp
		\displaystyle{\frac{\t p_\t}{\nu}+\frac{u_x}{\t}} \\
		\displaystyle{\frac{e_x}{\t}} \\
		\displaystyle{\frac{2z_x}{\t}}
		\ep\big(-1/s\big)\,.
 \end{align}
Also, from \eqref{eq:binverse}, we see immediately that 
\beq\label{eq:g32}
g_{32}=-\bp
\displaystyle{\frac{\t}{\nu}} & 0 & 0 \\
-\displaystyle{\frac{\t u}{\kappa_v}} & \displaystyle{\frac{\t}{\kappa_v}} & 0   \\
0 & 0 & \displaystyle{\frac{\t^2}{d}}
\ep\,.
\eeq
Finally, we compute 
\begin{align}
g_{33} & = b^{-1}(-A_{21}A_{11}^{-1}A_{12}+A_{22}) \\
	&=\bp
	\frac{\t}{\nu}\mathbf{c}_1  & \frac{\t}{\nu}p_e & \frac{\t}{\nu}p_z \\
	-\frac{\t u}{\kappa_v}\mathbf{c}_1 + \frac{\t}{\kappa_v}\mathbf{c}_2 & -p_e\frac{\t u}{\kappa_v}-s\frac{\t}{\kappa_v}+\frac{\t}{\kappa_v}up_e & 0 \\
	\frac{\t^2}{d}\mathbf{c}_3 & 0  & -s\frac{\t^2}{d}
	\ep
\end{align}
with 
\begin{align*}
\mathbf{c}_1 & = -\frac{1}{s}\left(p_\t+\frac{\nu u_x}{\t^2}\right)-s \,,\\
\mathbf{c}_2 & =-\frac{1}{s}\left(p_\t u+\frac{\nu u u_x}{\t^2}+\frac{\kappa_v e_x}{\t^2}\right)-su + p -\frac{\nu u_x}{\t} \,,\\
\mathbf{c}_3 & = -\frac{1}{s}\left(2d\frac{z_x}{\t^3}\right)\,.\\
\end{align*}
This completes the derivation of the entries in the coefficient matrix $G$ in \eqref{eq:1order}.


\end{document}